\documentclass{preprint}

\usepackage{mathtools}
\usepackage{mathrsfs}
\usepackage{tikz}
\usetikzlibrary{math}

\usepackage{comment}

\title{Study of the finite element discretization of the iterated quasi-reversibility method \\ The Laplacian data assimilation case}

\addauthor{Jérémi}{Dardé}{jeremi.darde@math.univ-toulouse.fr}{Université de Toulouse, INSA Toulouse, CNRS, Institut de Mathématiques de Toulouse, Toulouse, France.}
\addauthor{Jérémy}{Heleine}{jeremy.heleine@math.univ-toulouse.fr}{Université de Toulouse, INSA Toulouse, CNRS, Institut de Mathématiques de Toulouse, Toulouse, France.}

\abstract{
	We investigate the finite element discretization of the iterated quasi-reversibility method, in the context of the data assimilation problem for the Laplacian.
	We provide error estimates explicit in all the parameters appearing in the model and the method (\textit{i.e.} the amplitude of noise, the size mesh, the regularization
	parameter and the number of iterations of the method). This allows to devise strategies to choose the aforementioned parameters, for example to guarantee convergence of the method when the level of noise goes to zero.
	Numerical experiments on a two-dimensional problem are provided, illustrating the relevance of the obtained estimates.
}

\begin{document}

\tableofcontents

	\section{Introduction}

This paper investigates the finite element discretization of the iterated quasi-reversibility method in the context of the data assimilation problem for the Laplace equation.
The problem under consideration is as follows: let $\Omega$ be a bounded open polygon ($d= 2$) or polyhedron ($d=3$), and $\omega$ be a non-empty open set compactly contained in $\Omega$.
	Our objective is to reconstruct a function $u$, belonging to $H^1(\Omega)$, that satisfies the following system:
	\begin{equation} \label{eq_dataassimpb} \tag{$\mathcal P_{DA}$}
		\Delta u = f^\delta + \nabla \cdot \mathbf{F}^\delta \text{ in } \Omega, \quad u = g^\delta \text{ in } \omega.
	\end{equation}
	In this system, $f^\delta \in L^2(\Omega)$ and $F^\delta \in L^2(\Omega)^d$ are source terms, while $g^\delta \in L^2(\omega)$ is a measurement within $\omega$ of the function $u^\delta$ that we seek to reconstruct.
	The triplet $(f^\delta, F^\delta, g^\delta)$ constitutes the problem data. The role of the parameter $\delta$  will be clarified in the subsequent discussion.

	Problem \eqref{eq_dataassimpb} is the classical data assimilation problem for the Laplacian operator, which is known to be a severely ill-posed inverse problem. While the unique continuation property satisfied by the Laplace operator \cite[Chapter 5]{lerousseaulebeaurobbiano22} ensures that Problem \eqref{eq_dataassimpb} has at most one solution, the existence of such a solution is not assured.
%

We introduce a triplet of \emph{exact}, \emph{admissible}, or \emph{idealized data}
$(f, \mathbf{F}, g) := (f^0, \mathbf{F}^0, g^0) \in L^2(\Omega) \times L^2(\Omega)^d \times L^2(\omega)$,
for which we assume the existence of a solution \( u_s \) to the inverse problem~\eqref{eq_dataassimpb}.
In practice, however, the exact data is inevitably contaminated by noise.
To account for this, we define a triplet of \emph{noisy} or \emph{perturbed data}
$(f^\delta, \mathbf{F}^\delta, g^\delta)$,
where the known parameter \(\delta \geq 0\) quantifies the magnitude of the perturbation:

	\begin{equation} \label{eq_def_noise_Laplace}
		\sqrt{\Vert f^\delta - f\Vert_{L^2(\Omega)}^2 + \Vert \mathbf{F}^\delta - \mathbf{F} \Vert^2_{L^2(\Omega)} + \Vert g^\delta - g \Vert^2_{L^2(\omega)}} \leq \delta.
	\end{equation}
	We do not assume that the data \((f^\delta, \mathbf{F}^\delta, g^\delta)\) is admissible; consequently, Problem~\eqref{eq_dataassimpb} may not admit a solution. Furthermore, due to the ill-posed nature of the problem, even if a solution \(u^\delta\) exists, it may deviate significantly from the target function \(u_s\) we seek to reconstruct. In fact, for any \(\delta > 0\), there exists an admissible triplet \((f^\delta, \mathbf{F}^\delta, g^\delta)\) satisfying~\eqref{eq_def_noise_Laplace} for which the corresponding solution \(u^\delta\) satisfies
\[
\| u^\delta - u_s \|_{H^1(\Omega)} \geq \frac{1}{\delta}.
\]
Therefore, the objective is not to reconstruct \(u^\delta\) itself, but rather to obtain an approximate solution \(\tilde{u}^\delta\), dependent on parameters such as noise amplitude, regularization, and discretization, while controlling the error \(\| \tilde{u}^\delta - u_s \|\) in terms of these parameters. This enables, for instance, the selection of parameters to guarantee convergence of \(\tilde{u}^\delta\) to \(u_s\) as the noise amplitude tends to zero.

	\subsection*{Iterated quasi-reversibility method for the data assimilation problem}

	A wide range of methods has been proposed to stabilize problem \eqref{eq_dataassimpb} and approximate its solution; a concise overview is provided in \autoref{sec_ScientifContext}. Here, we focus on the \emph{iterated quasi-reversibility method} \cite{darde2016}, applied to the $H^1 \times H^1$ mixed formulation initially introduced in \cite{bourgeois2005} for a closely related problem, the \emph{data completion problem},  and later revisited in \cite{bourgeoisrecoquillay2018}. To describe the method, we first introduce
	$\lambda^\delta$, the unique element of $H^1_0(\Omega)$ verifying
	\begin{equation} \label{eq_def_lambda}
	\int_\Omega \nabla \lambda^\delta \cdot \nabla \mu \, dx = \int_\Omega (f^\delta\, \mu - \mathbf F^\delta\cdot \nabla \mu) dx, \quad \forall \mu \in H^1_0(\Omega).
	\end{equation}
	The Lax-Milgram theorem  \cite[Corollary 5.8]{brezis2011} ensures the existence and uniqueness of $\lambda^\delta$, as well as its continuous dependence on $f^\delta$ and $\mathbf F^\delta$. Additionally, if Problem \eqref{eq_dataassimpb} has a solution $u^\delta$, then $\lambda^\delta$ is the orthogonal projection of $u^\delta$ onto $H^1_0(\Omega)$, meaning it is the unique element in $H^1_0(\Omega)$ such that
	$$
	\int_\Omega \nabla \lambda^\delta \cdot \nabla \mu dx = \int_\Omega \nabla u^\delta \cdot \nabla \mu dx, \quad \forall \mu \in H^1_0(\Omega).
	$$
	In particular, $\lambda^0$ is the orthogonal projection of $u_s$ onto $H^1_0(\Omega)$. Hereafter, we denote $\lambda^0$ as $\lambda_s$.

	At the continuous level (which means, in this study, prior to finite-element discretization), the \emph{iterated quasi-reversibility method} proceeds as follows.
	First, we choose a \emph{moderately small} regularization parameter $\varepsilon>0$.
	We then define $(u_\varepsilon^{M,\delta},\lambda_{\varepsilon}^{M, \delta})_{M\in \N}$, referred to as \emph{iterated quasi-reversibility sequence of approximate solutions to Problem \eqref{eq_dataassimpb}},  recursively. The initial pair is $(u_\varepsilon^{0,\delta},\lambda_\varepsilon^{0,\delta}) = (0,0)$ and, for each $M\in \N$, the pair $(u_\varepsilon^{M+1,\delta},\lambda_\varepsilon^{M+1,\delta}) \in H^1(\Omega)\times H^1_0(\Omega)$
	satisfies:
	\begin{equation} \label{eq_defQRiterated}
	\begin{cases}
	\displaystyle	\int_\Omega \nabla \lambda_{\varepsilon}^{M+1,\delta} \cdot \nabla v\, dx
	+ \int_\omega u_{\varepsilon}^{M+1,\delta}\, v\, dx + \varepsilon \int_\Omega \nabla u_{\varepsilon}^{M+1,\delta}\cdot \nabla v \, dx = \int_\Omega \nabla \lambda^\delta \cdot \nabla v \, dx + \int_\omega g^\delta\, v \, dx\\ \displaystyle \hfill + \varepsilon \int_\Omega \nabla u_{\varepsilon}^{M,\delta} \cdot \nabla v \, dx,  \quad
	\forall v\in H^1(\Omega), \\ \\
	\displaystyle	\int_\Omega \nabla \lambda_{\varepsilon}^{M+1,\delta} \cdot \nabla \mu \, dx - \int_\Omega \nabla u_{\varepsilon}^{M+1,\delta}\cdot \nabla \mu \, dx = 0, \quad \forall \mu \in H^1_0(\Omega).
	\end{cases}
	\end{equation}
	For simplicity, we denote  the sequence corresponding to the exact data (\textit{i.e.}
	$(u_\varepsilon^{M,0},\lambda_\varepsilon^{M,0})_{M\in \N}$) as $(u_\varepsilon^M , \lambda_\varepsilon^M)_{M\in \N}$.

	Although the primary focus of the article is the \emph{discretized quasi-reversibility sequence} introduced in the next section, we begin by presenting several results for the continuous setting. We start with the sequence corresponding to unperturbed data.
	We define the set
	$$
	H^1(\Omega,\Delta) := \left\lbrace v \in H^1(\Omega), \ \Delta v \in L^2(\Omega) \right\rbrace,
	$$
	which, equipped with the norm $\Vert v \Vert_{H^1(\Omega,\Delta)} := \sqrt{\Vert v \Vert_{H^1(\Omega)}^2 + \Vert \Delta v \Vert_{L^2(\Omega)}^2}$, is a Hilbert space.

	\begin{proposition} \label{prop_intro_first}
	For any $\varepsilon>0$, the sequence $(u_\varepsilon^M,\lambda_\varepsilon^M)_{M\in \N}$ is well-defined and converges strongly to $(u_s,\lambda_s)$ in $H^1(\Omega)\times H^1(\Omega)$ as $M \rightarrow \infty$.

	Moreover, there exists $\mu \in (0,1)$, and for every $\varepsilon_0>0$, a constant $c>0$, such that for all $u_s \in H^1(\Omega,\Delta)$,
	for all $\varepsilon\in (0,\varepsilon_0)$ and for all $M\in \N$,
	$$
	\Vert u_\varepsilon^M - u_s \Vert_{L^2(\Omega)} \leq c \ln \left(1+\frac{M}{\varepsilon}  \right)^{-\mu} \Vert u_s \Vert_{H^1(\Omega,\Delta)}.
	$$
	\end{proposition}

	The convergence of the sequence follows directly from the abstract framework developed in \autoref{sec_abstractQR}. The convergence estimate, meanwhile, is derived from the regularity estimates obtained in  \autoref{sec_contQRregres} and is proven in  \autoref{subsec_logconv}. Several observations are worth noting:

	\begin{itemize}
	\item Most regularization methods for Problem \eqref{eq_dataassimpb}, including the  \emph{standard quasi-reversibility method}, rely on a single parameter, usually called the \emph{regularization parameter}, which serves two conflicting purposes. On the one hand, it regularizes the problem, ensuring its numerical solvability, and thus should be chosen not too small. On the other hand, it acts as the convergence parameter, requiring it to be sufficiently small. Balancing these roles makes the selection of an appropriate parameter inherently challenging.

	An advantage of the iterated quasi-reversibility method is that these two roles are now decoupled into two separate parameters. The regularization parameter $\varepsilon$ solely governs regularization and is therefore chosen to be moderately small (typically between $0.01$ and $1$). The parameter $M$, in contrast, drives convergence, and is thus chosen to be large. This complementarity is evident in the convergence estimate, which highlights the importance of the ratio $\frac{M}{\varepsilon}$. This ratio will recur throughout the various convergence estimates in this study.

	\item We cannot directly apply the abstract framework for the iterated quasi-reversibility method developed in \cite{darde2016}. The theory in \cite{darde2016} was established for symmetric bilinear forms and a specific type of data. However, the bilinear form on the left-hand side of \eqref{eq_defQRiterated} is not symmetric, and the right-hand side data does not conform to the required form.
\autoref{sec_abstractQR} primarily focuses on adapting the known results to this new context, thereby generalizing the work presented in \cite{darde2016}.

	\item The relatively slow logarithmic convergence rate of the reconstruction method is to be expected, given the ill-posed nature of the data completion problem (see \cite[Section 6.1]{caubetdarde2020}, where the same logarithmic convergence rate emerges for a dual Kohn-Vogelius method applied to the Laplacian data completion problem). This is reminiscent of the logarithmic-type conditional stability estimate, classicaly obtained for elliptic data assimilation problems \cite{alessandrinirondirossetvessella2009,lerousseaulebeaurobbiano22}.

	Also, note that this convergence rate relies on a regularity assumption for the exact solution $u_s$, specifically that $\Delta u_s \in L^2(\Omega)$ , which is equivalent to $\mathbf{F}^0 = 0$.
	\end{itemize}

	Although not essential for the study of the discretization of the iterated quasi-reversibility method, we incidentally obtain
	an error estimate for the sequence with perturbed data,  $(u_\varepsilon^{M,\delta},\lambda_\varepsilon^{M,\delta})$.

	\begin{proposition} \label{prop_intro_second}
	The sequence $(u_\varepsilon^{M,\delta}, \lambda_\varepsilon^{M,\delta})_{M\in \N}$ is well-defined.

	There exists a constant $c>0$ such that, for all $(f^\delta,F^\delta,g^\delta)$ satisfying \eqref{eq_def_noise_Laplace}, for all $\varepsilon>0$ and all $M\in \N$, the following estimate holds
	$$
	\Vert u_\varepsilon^{M,\delta} - u_s\Vert_{H^1(\Omega)} \leq c \left(\sqrt{\frac{M}{\varepsilon}} \delta + \Vert u_\varepsilon^M - u_s \Vert_{H^1(\Omega)} + \Vert \lambda_\varepsilon^M - \lambda_s \Vert_{H^1(\Omega)} \right).
	$$
	Furthermore, there exists $\mu \in (0,1)$ and $c>0$ such that for all $u_s \in H^1(\Omega,\Delta)$, for all $(f^\delta,F^\delta,g^\delta)$ satisfying \eqref{eq_def_noise_Laplace}, for all $\varepsilon>0$ and all $M\in \N$,
	$$
	\Vert u_\varepsilon^{M,\delta} - u_s\Vert_{L^2(\Omega)} \leq c \left(\sqrt{\frac{M}{\varepsilon}} \delta + \ln \left(1+\frac{M}{\varepsilon} \right)^{-\mu} \Vert u_s \Vert_{H^1(\Omega,\Delta)}\right).
	$$
	\end{proposition}

	As suggested by the term~$\sqrt{\frac{M}{\varepsilon}} \delta$, the ratio~$\frac{M}{\varepsilon}$ should not be chosen too large when dealing with perturbed data, as this may amplify the effect of perturbations on the reconstruction quality. If $\varepsilon$ is fixed, selecting $M$ as a function of $\delta$ so that $M(\delta) \xrightarrow{\delta\rightarrow 0} +\infty$ and $\sqrt{M(\delta)} \delta \xrightarrow{\delta \rightarrow \infty} 0$ ensures the convergence of $u_\varepsilon^{M(\delta),\delta}$ to $u_s$.
	Given that this term will also appear in the estimate for the discretized iterated quasi-reversibility method, we will explore this issue in greater detail in the next section.

\subsection*{Discretization of the iterated quasi-reversibility method} \label{subsection_intro_discretization}

We now turn our attention to the numerical resolution of the iterated quasi-reversibility problem \eqref{eq_defQRiterated} using a finite element approach.

Specifically, for $h>0$, we  introduce $\mathcal T_h$, an admissible, shape regular and uniform family of partitions of $\Omega$
into triangular elements ($d=2$) or tetrahedral elements ($d=3$). We assume that each element $\tau_h \in \mathcal{T}_h$ has a diameter  bounded by $c\,h$, where $c>0$ is a fixed constant.
We define the finite element spaces as follows:
$$
V_h = \left\lbrace w \in H^1(\Omega), \ \forall \tau_h \in \mathcal{T}_h, \ w_{|\tau_h} \in \mathbb{R}_1[X] \right\rbrace, \quad V_{h,0} = V_h \cap H^1_0(\Omega).
$$
Here, $V_h$ is the standard \emph{conforming $P_1$-Lagrange finite element space}, while $V_{h,0}$ is the subspace of functions of $V_h$ that vanish on $\partial \Omega$.
For further details on $\mathcal{T}_h$, $V_h$ and $V_{h,0}$, we refer the reader to \cite[Section 9]{steinbach08}.

We now introduce the sequence of \emph{discretized iterated quasi-reversibility solutions}, denoted by $(u_{\varepsilon,h}^{M,\delta}, \lambda_{\varepsilon,h}^{M,\delta})_{M\in \N}$. Like its continuous counterpart, this sequence is defined recursively as follows: the initial pair is $(u_{\varepsilon,h}^{0,\delta},\lambda_{\varepsilon,h}^{0,\delta}) = (0,0)$, and for each $M\in \N$, the pair $(u_{\varepsilon,h}^{M+1,\delta},\lambda_{\varepsilon,h}^{M+1,\delta})$ is the unique element in $V_h\times V_{h,0}$
satisfying

\begin{equation} \label{eq_itQRDAapproxnoisy_core}
\begin{cases}
\displaystyle \int_\Omega \nabla \lambda_{\varepsilon,h}^{M+1,\delta} \cdot \nabla v_h dx + \int_\omega u_{\varepsilon,h}^{M+1,\delta} v_h dx + \varepsilon \int_\Omega \nabla u_{\varepsilon,h}^{M+1,\delta} \cdot \nabla v_h dx
=  \int_\Omega \nabla \lambda^\delta \cdot \nabla v_h dx + \int_\omega g^\delta\, v_h dx
\\ \hfill \displaystyle  + \varepsilon \int_\Omega \nabla u_{\varepsilon,h}^{M,\delta} \cdot \nabla v_h \, dx, \quad \forall v_h \in V_h,
\\
\\ \displaystyle \int_\Omega \nabla \lambda_{\varepsilon,h}^{M+1,\delta} \cdot \nabla \mu_h dx - \int_\Omega \nabla u_{\varepsilon,h}^{M+1,\delta} \cdot \nabla \mu_h dx = 0, \hfill \quad \forall \mu_h \in V_{h,0}.
\end{cases}
\end{equation}

Once again, the sequence is well-defined. The primary objective of this study is to derive an explicit bound on the discrepancy \(\Vert u_{\varepsilon,h}^{M,\delta} - u_s \Vert\), with a focus on clarifying its dependence on the key parameters: the regularization parameter~\(\varepsilon\), the sequence index~\(M\), the perturbation amplitude~\(\delta\), and the mesh size~\(h\).

As is typical when deriving finite element convergence estimates, we first need to establish \textit{a priori} regularity estimates at the continuous level, specifically for \( u_\varepsilon^M \) and \( \lambda_\varepsilon^M \). To the best of our knowledge, the only study addressing regularity estimates for the quasi-reversibility method is~\cite{bourgeoischesnel2020}. However, this work focuses on a different context: the authors investigate a data completion problem, which is arguably more complex from a regularity perspective due to intricate boundary conditions. Their analysis, however, is limited to the \emph{standard quasi-reversibility method}, a simpler case that corresponds to the iterated quasi-reversibility method terminated after a single iteration.

In  \autoref{sec_contQRregres}, we derive the following regularity estimates:

\begin{theorem} \label{thm_intro_first}
There exists $\sigma \in \left]\frac{1}{2},1\right]$ depending only on $\Omega$, and for every $\varepsilon_0>0$, a constant $c>0$
depending only on $\Omega$ and $\varepsilon_0$, such that for all $u_s \in H^1(\Omega,\Delta)$, all $\varepsilon \in (0,\varepsilon_0]$ and all
$M\in \N$, the functions $u_\varepsilon^M$ and $\lambda_\varepsilon^M$ belong to $H^{1+\sigma}(\Omega)\cap H^1(\Omega,\Delta)$. Moreover, the following estimates hold:
$$
\Vert \lambda_\varepsilon^M \Vert_{H^{1+\sigma}(\Omega)} \leq c \left(\left(1 - \frac{\varepsilon^M}{(1+\varepsilon)^M} \right)
 \Vert \Delta u_s \Vert_{L^2(\Omega)} + \sqrt{\frac{\varepsilon}{M}} \Vert u_s \Vert_{L^2(\Omega)} \right),
$$
$$
\Vert u_\varepsilon^M \Vert_{H^{1+\sigma}(\Omega)} \leq c \sqrt{\frac{M}{\varepsilon}} \Vert u_s \Vert_{H^1(\Omega,\Delta)}.
$$
In the case where $\Omega$ is convex, these estimates hold with $\sigma = 1$.
\end{theorem}

Before presenting the convergence estimates for the discretized quasi-reversibility solutions, we provide some observations regarding \autoref{thm_intro_first}.

\begin{itemize}
\item By setting $M=1$ in the estimates, we recover regularity estimates for the standard quasi-reversibility method, exhibiting a similar dependence on the regularization parameter $\varepsilon$ as observed in  \cite[Theorem 2.2]{bourgeoischesnel2020}.
\item The presence of \emph{reentrant corners} limits the gain in regularity. In two dimensions, the parameter~$\sigma$ can be expressed explicitly in terms of the angles of these reentrant corners~\cite{grisvard85,grisvard92}. When~$\Omega$ is convex (\textit{i.e.} in the absence of reentrant corners) $H^2$ regularity is guaranteed.
\item
There is an improvement in regularity from \( u_s \) to \( u_\varepsilon^M \):

\[
u_s \in H^1(\Omega), \ \Delta u_s \in L^2(\Omega) \implies u_\varepsilon^M \in H^{1+\sigma}(\Omega), \ \Delta u_\varepsilon^M \in L^2(\Omega).
\]

However, even if \( u_s \) is smoother (for instance, \( u_s \in H^p(\Omega) \) with \( p \geq 3 \)), this does not enhance the regularity of \( u_\varepsilon^M \), even when \( \Omega \) is convex. This limitation arises from the quasi-reversibility formulation, specifically due to the two integrals over \( \omega \), which inhibit any further gain in regularity.
This issue is discussed further in \autoref{sec_contQRregres}, particularly in \autoref{rk_limitation_regularity},  and in \autoref{section_further_comments}.
\item The limited regularity of $(u_\varepsilon^M, \lambda_\varepsilon^M)$, which is at most in $H^2$, makes it ineffective to increase the polynomial order of the finite element approximation. Doing so would not improve the convergence rate of the method. Therefore, we limit our analysis to $P1$-Lagrange finite elements.
\end{itemize}

These estimates are essential for deriving our main result: the error estimate between the $M$-th term of the discretized quasi-reversibility sequence with perturbed data, $u_{\varepsilon,h}^{M,\delta}$, and the exact solution, $u_s$. 

\begin{theorem} \label{thm_intro_second}
There exists  $\sigma \in \left(\frac{1}{2},1\right]$ and $\mu \in (0,1)$, depending only on $\Omega$, and for every $\varepsilon_0>0$, a constant $c>0$ depending only on $\Omega$ and $\varepsilon_0$, such that for all $(f^\delta,\mathbf{F}^\delta,g^\delta)$ satisfying \eqref{eq_def_noise_Laplace}, for all $\varepsilon \in (0,\varepsilon_0)$, for all $M\in \N$ and all $h>0$,
\begin{multline*}
\left\Vert u_{\varepsilon,h}^{M,\delta} - u_s \right\Vert_{H^1(\Omega)} \leq
c \Bigg( \sqrt{\frac{M}{\varepsilon}} \delta + \left(\frac{M}{\varepsilon}\right)^{\frac{3}{2}} h^{\sigma} \Vert u_s \Vert_{H^1(\Omega,\Delta)}
 \\  + \left\Vert u_\varepsilon^M - u_s \right\Vert_{H^1(\Omega)} + \left\Vert \lambda_\varepsilon^M - \lambda_s \right\Vert_{H^1(\Omega)} \Bigg),
\end{multline*}
$$
\left\Vert u_{\varepsilon,h}^{M,\delta} - u_s \right\Vert_{L^2(\Omega)} \leq  c\left( \sqrt{\frac{M}{\varepsilon}} \delta + \left(\frac{M}{\varepsilon}\right)^{\frac{3}{2}} h^{\sigma} \Vert u_s \Vert_{H^1(\Omega,\Delta)} + \ln \left(1+ \frac{M}{\varepsilon}\right)^{-\mu} \Vert u_s \Vert_{H^1(\Omega,\Delta)}\right).
$$
\end{theorem}

A few remarks are in order:

\begin{itemize}
	\item The right-hand side of each estimate consists of three terms, each corresponding to a distinct type of error.
The leftmost term is the \emph{noise error}, which quantifies the impact of perturbations in the data. This term provides an estimate of $\Vert u_{\varepsilon,h}^{M,\delta} - u_{\varepsilon,h}^M \Vert$, as derived in \autoref{prop_errorestimate_noise}.

The middle term is the \emph{discretization error}, which arises from the finite-element discretization process. It estimates $\Vert u_{\varepsilon,h}^M - u_\varepsilon^M \Vert$, as established in \autoref{thm_errorestimatediscretization}.

Finally, the rightmost term is the \emph{regularization error}, which measures the discrepancy between the quasi-reversibility solution and the target solution, $\Vert u_\varepsilon^M - u_s \Vert$. This error can only be quantified in the $L^2$ norm, as shown in \autoref{prop_logestimateQR}.
The global estimate follows directly from the triangle inequality.
    \item \autoref{thm_intro_second} provides a framework for selecting the parameters $M$, $\varepsilon$, and $h$ to ensure that the reconstruction $u_{\varepsilon,h}^{M,\delta}$ converges to the exact solution $u_s$ as the perturbation amplitude $\delta$ tends to zero. Convergence is guaranteed for any choice of these parameters that satisfies the following conditions:
\[
\frac{M}{\varepsilon} \xrightarrow{\delta \to 0} +\infty, \quad \frac{M}{\varepsilon} \delta^2 \xrightarrow{\delta \to 0} 0, \quad \left( \frac{M}{\varepsilon}\right)^3 h^{2\sigma} \xrightarrow{\delta \to 0} 0.
\]

The last condition depends on the parameter $\sigma$, whose value can be computed (see, e.g., \cite{grisvard85,grisvard92}). However, it can be replaced by the stronger condition
\[
\left( \frac{M}{\varepsilon}\right)^3 h \xrightarrow{\delta \to 0} 0,
\]
which guarantees convergence for any polygonal or polyhedral domain $\Omega$.

Finally, the even weaker condition
\[
\left( \frac{M}{\varepsilon}\right)^3 h^2 \xrightarrow{\delta \to 0} 0
\]
suffices to ensure convergence for any convex domain $\Omega$.

\item In practice, the partition $\mathcal{T}_h$ is given, which determines the value of $h$. The maximum noise amplitude is also specified. This leaves two parameters to be chosen: $\varepsilon$ and $M$. We first set $\varepsilon$ to a moderately small value, typically between $0.01$ and $1$, to ensure that the linear system to be solved is well-conditioned.

Next, we consider the parameter $M$. A natural approach is to balance the three terms in the second estimate, aiming to give them comparable influence. Such strategy is closely related to the balancing principle proposed in \cite{caoklibanovpereverzev09}. However, it is important to note that, by definition, the $H^1(\Omega,\Delta)$-norm of $u_s$ is unknown, since $u_s$ itself is unknown, and the convergence of such a strategy as $\delta$ and $h$ approach zero is not guaranteed.

\item When expanded in the classical nodal basis of $\mathcal{X}_h$, Equation~\eqref{eq_itQRDAapproxnoisy_core} yields a sequence of linear systems, which can be schematically written as
$\mathbb{M} \mathbf{U}^{M+1} = \mathbf{F}^\delta + \mathbb{B} \mathbf{U}^M$.

An important observation is that the linear system to be inverted remains unchanged across iterations: the matrix $\mathbb{M}$ depends only on $\mathcal{T}_h$ and $\varepsilon$, and not on the iteration index $M$. Consequently, $\mathbb{M}$ can be precomputed and factorized, rendering the computational cost of solving the system negligible as the number of iterations grows.

\item This study does not account for computational-level errors, including those stemming from the numerical computation of the integrals in \eqref{eq_itQRDAapproxnoisy_core} and the numerical inaccuracies introduced during the resolution of the associated linear system. The impact of these errors will be addressed in future research.
\end{itemize}

	\subsection{Scientific context} \label{sec_ScientifContext}

In data assimilation problems, the main objective is to  reconstruct a system’s full state from partial internal observations, and the knowledge of a mathematical model of the underlying physics of the problem (usually, a system of partial differential equations).
As such, data assimilation problems naturally appear in geophysics, such as weather forecasting and oceanography \cite{auroux17,Navon2009}, but also in biomedical applications  \cite{chapellefragumalletmoireau13,quarteronimanzonivergara17}.
Data assimilation problems are closely related to data completion problems, the key distinction being that the latter assumes data is available on the domain’s boundary.  Problem \eqref{eq_dataassimpb}
serves as a simplified (toy) example of a data assimilation problem, our primary goal being to derive explicit error estimates for our reconstruction method within this simplified framework.

Mathematically, Problem \eqref{eq_dataassimpb} constitutes an inverse and ill-posed problem. The uniqueness of its solution stems from the well-known unique continuation property, while its severe instability is captured by logarithmic-type conditional stability estimates, typically derived via Carleman estimates \cite{alessandrinirondirossetvessella2009,Holmgren1901,lerousseaulebeaurobbiano22,vessella25}.


The numerical reconstruction of solutions for data assimilation and data completion problems involving partial differential equations (PDEs) is a highly active research area, with a wide range of methods proposed. Among these, and without aiming to be exhaustive, we mention methods based on stabilized finite elements~\cite{burmanhansbolarson18,burmanoksanen18,burmangargpreuss24}, $L^2$-minimization with PDE constraints~\cite{cindeamunch15,cindeamunch16}, surface integral methods~\cite{boukarihaddar15,abidetal23}, fading parameter methods~\cite{Caille_Hanus_Delvare_Michaux_Leblonda_2018,caille19,abidboukaridelvare26}, convexification methods~\cite{beilinaklibanov12}, duality-based methods~\cite{caubetdarde2020,bourgeoisdarde22,bourgeoisetal25}, and Kohn-Vogelius approaches~\cite{andrieux05,benbelgacem05,Azaiez06,benfatma07,BenBelgacem22,BENBELGACEM2023}.

The quasi-reversibility method, introduced by Robert Latt\`es and Jacques-Louis Lions in the late 1960s~\cite{latteslions69}, represents another approach to addressing these problems, as does its variant, the iterated quasi-reversibility method. Both methods have been extensively studied and applied across a range of contexts (see, \emph{e.g.},~\cite{bourgeois2005,caoklibanovpereverzev09,darde2016,bourgeoisrecoquillay2018,bourgeoischesnel2020,darbasetal20,le21,heleine25}).

	\subsection{Outline of the article}

The article is organized as follows.
In \autoref{sec_abstractQR}, we develop an abstract framework for the iterated quasi-reversibility method, suited to our problem of interest.
\autoref{sec_contQRregres} is devoted to proving the regularity results presented in \autoref{thm_intro_first}.
The finite element discretization of the method is thoroughly analyzed in \autoref{sec_FEMdisc}.
Numerical results are discussed in \autoref{sec_NumIll}, while \autoref{section_further_comments} provides additional comments and insights on the study.
Finally, technical results are gathered in \autoref{appendix_miscres}.

Throughout this paper, the notation \(A \lesssim B\) signifies that there exists a constant \(c > 0\), independent of all relevant parameters (such as \(\delta\), \(\varepsilon\), \(M\), and \(h\)), such that \(A \leq cB\).

\section{Abstract iterated quasi-reversibility} \label{sec_abstractQR}

\subsection{Standard quasi-reversibility}

 Let $\mathcal X$ be a Hilbert space, with its scalar product denoted by $(\cdot,\cdot)_\mathcal X$ and the corresponding norm $\Vert\cdot\Vert_{\mathcal X}$.
Let  $a : \mathcal X \times \mathcal X \mapsto \R$ be a continuous positive bilinear form. For any $x \in \mathcal X$, we denote $\Vert x \Vert_a = \sqrt{a(x,x)}$. Note that $\Vert . \Vert_a$ is well-defined since $a$ is positive, but it is not necessarily a norm on $\mathcal X$.

Finally, let $L : \mathcal X \mapsto \R$ be a linear form such  that for all $x$ in $ \mathcal X$,
$
\vert L (x) \vert \leq c_a \Vert x \Vert_a,
$
for some nonnegative constant $c_a$.
Since $\Vert \cdot \Vert_a \lesssim \Vert \cdot \Vert_{\mathcal X}$, $L$ is automatically continuous. Our abstract problem of interest reads: find $x \in \mathcal X$ such that
\begin{equation} \label{prob_abstract_a} \tag{$\mathscr{P}$}
a(x,\tilde x)  = L(\tilde x),\quad \forall \tilde x \in X.
\end{equation}
Clearly, such problem is not necessarily well-posed. We make the following assumption:

\begin{assumption}[Identifiability assumption] \label{assumption_identifiability}
	If $a(x,\tilde x) = 0$ for all $\tilde x$, then $x = 0$.
\end{assumption}

Under this assumption, Problem \eqref{prob_abstract_a} admits at most one solution. However, the problem remains ill-posed and may still fail to have a solution. To address this, we introduce a bilinear, continuous, symmetric and positive form $b: \mathcal X \times \mathcal X \mapsto \R$. We set $\Vert x \Vert_b = \sqrt{b(x,x)}$, and assume that
$$
\Vert \cdot \Vert_{Q_1} : x \in \mathcal X \mapsto \sqrt{\Vert x \Vert_a^2 + \Vert x \Vert_b^2},
$$
is a norm on $\mathcal X$, equivalent to the norm $\Vert \cdot \Vert_{\mathcal X}$. Such a form \(b\) always exists; indeed, one can simply take \(b\) to be the scalar product on \(\mathcal{X}\).

We now introduce our regularized problem: for $\varepsilon>0$, find $x_\varepsilon \in \mathcal X$ such that
\begin{equation} \label{prob_abstract_Qe} \tag{$\mathscr P_\varepsilon$}
Q_\varepsilon(x_\varepsilon,\tilde x) := a(x_\varepsilon,\tilde x) + \varepsilon b(x_\varepsilon,\tilde x) = L(\tilde x),\quad \forall \tilde x \in X.
\end{equation}
By Lax-Milgram theorem, such problem is well-posed, hence $x_\varepsilon$ exists and is unique. Furthermore, it satisfies the following a priori estimates.

\begin{lemma} \label{lemma_apriori_estimates_standard_reg}
	For all $\varepsilon> 0$, we have
	$
	\displaystyle \Vert x_\varepsilon\Vert_a \leq c_a$, and $\displaystyle  \sqrt{\varepsilon} \Vert x_\varepsilon \Vert_b \leq c_a.
	$
\end{lemma}

\begin{proof}
	It is a direct consequence of the identity
	$$
	Q_\varepsilon(x_\varepsilon,x_\varepsilon) = \Vert x_\varepsilon\Vert_a^2 +\varepsilon \Vert x_\varepsilon\Vert_b^2 = L(x_\varepsilon) \leq c_a \Vert x_\varepsilon \Vert_a.
	$$
\end{proof}

%

\begin{lemma}
	Suppose there exists $I \subset (0,\infty)$ such that $0 \in \overline I$ and $(\Vert x_\varepsilon\Vert_b)_{\varepsilon \in I}$ is uniformly bounded. Then problem \eqref{prob_abstract_a} admits a unique solution $x$ and there exists a sequence of positive real numbers $(\varepsilon_n)_{n\in \N}$ converging towards zero such that $x_{\varepsilon_n}$
	weakly converges to $x$ in $\mathcal X$.
\end{lemma}

\begin{proof}
	Let $(\varepsilon_n)_{n\in \N} \in I^\N$ be converging to zero, and denote $x_n = x_{\varepsilon_n}$. By assumption, we have $\Vert x_n \Vert_b \lesssim 1$ for all integer $n$. But we also have $\Vert x_n \Vert_a \lesssim 1$ by \autoref{lemma_apriori_estimates_standard_reg}, which implies $\Vert x_n \Vert_{\mathcal X} \lesssim 1$. As a consequence, there exists a subsequence (which for convenience we still denote $x_n$), and $x\in \mathcal X$ such that $x_n$ weakly converges to $x$ in $\mathcal X$. In particular, for all $\tilde x  \in X$, we have
	$$
	a(x,\tilde x) \xleftarrow{n\rightarrow \infty} a(x_n,\tilde x) + \varepsilon_n b(x_n,\tilde x) = L(\tilde x).
	$$
	In other words, $x$ is a solution of \eqref{prob_abstract_a}. Uniqueness follows from \autoref{assumption_identifiability}.
\end{proof}

\begin{proposition}
	Suppose problem \eqref{prob_abstract_a} admits a unique solution $x_s$. Then $x_\varepsilon$ strongly converges to $x_s$ in $\mathcal X$ as $\varepsilon$ tends to zero, and we have
	$$
	  \Vert x_\varepsilon\Vert_b \leq \Vert x_s \Vert_b,\quad\Vert x_\varepsilon - x_s \Vert_a \leq \sqrt \varepsilon \Vert x_s \Vert_b,  \quad  \Vert x_\varepsilon-x_s \Vert_b \leq \Vert x_s \Vert_b.
	$$
\end{proposition}

\begin{proof}
	If $x_s$ is the solution of problem \eqref{prob_abstract_a}, then for all $\tilde x \in \mathcal X$, we have
	$$
	a(x_\varepsilon,\tilde x) + \varepsilon b(x_\varepsilon,\tilde x) = a(x_s,\tilde x),
	$$
	or equivalently
	$$
	a(x_\varepsilon - x_s, \tilde x) + \varepsilon b(x_\varepsilon,\tilde x) = 0.
	$$
	Setting $\tilde x = x_\varepsilon - x_s$, we obtain
	$$
	\Vert x_\varepsilon - x_s \Vert_a^2 + \varepsilon \Vert x_\varepsilon \Vert_b^2 =  \varepsilon b(x_\varepsilon,x_s) \leq \varepsilon \Vert x_\varepsilon \Vert_b \Vert x_s \Vert_b.
	$$
	The first estimate follows directly, and the second one is an immediate consequence. As furthermore
	$$
	\Vert x_\varepsilon - x_s \Vert_a^2 + \varepsilon \Vert x_\varepsilon - x_s \Vert_b^2 = - \varepsilon b(x_s, x_\varepsilon - x_s),
	$$
	the third estimate follows.

	As the family $(\Vert x_\varepsilon \Vert_b)_{\varepsilon>0}$ is uniformly bounded, the previous Lemma gives a sequence of positive numbers $(\varepsilon_n)_{n\in \N}$
	converging towards zero such that $x_{\varepsilon_n}$ weakly converges to $x$ as $n$ tends to infinity. As
	\begin{multline*}
	\Vert x_{\varepsilon_n} - x_s \Vert_{\mathcal X}^2 \lesssim \Vert x_{\varepsilon_n} - x_s \Vert_{Q_1}^2
	= \Vert x_{\varepsilon_n} - x_s \Vert_a^2 + \Vert x_{\varepsilon_n} - x_s \Vert_b^2
	= \Vert x_{\varepsilon_n} - x_s \Vert_a^2 + \Vert x_{\varepsilon_n}\Vert_b^2 - 2 b(x_{\varepsilon_n},x_s) +\Vert x_s\Vert_b^2 \\
	\leq {\varepsilon_n} \Vert x_s \Vert_b^2 + 2 \Vert x_s \Vert_b^2 - 2 b(x_{\varepsilon_n},x) \xrightarrow{n\rightarrow \infty} 0,
	\end{multline*}
	the sequence converges strongly. A standard proof by contradiction allows to conclude that the whole family converges as $\varepsilon$ goes to zero, which ends the proof.
\end{proof}

\subsubsection*{Differentiability}

All this section is adapted from \cite[Sections 3.1 and 3.2]{darde2016}.

For \(\varepsilon > 0\), we inductively define a sequence \((x_\varepsilon^m)_{m \in \mathbb{N}}\) in \(\mathcal{X}^\mathbb{N}\) as follows: set \(x_\varepsilon^0 = x_\varepsilon\), and for each \(m \in \mathbb{N}\), let \(x_\varepsilon^{m+1}\) be the unique element in \(\mathcal{X}\) satisfying
\begin{equation} \label{eq_definition_abstractderivatives}
Q_\varepsilon(x_\varepsilon^{m+1}, \tilde{x}) = - (m+1) \, b(x_\varepsilon^m, \tilde{x}), \quad \forall \tilde{x} \in \mathcal{X}.
\end{equation}
Lax-Milgram theorem ensures that the sequence is well-defined.

We now define a mapping \(F : (0, \infty) \to \mathcal{X}\) by \(F(\varepsilon) = x_\varepsilon\).

\begin{proposition}
The mapping \(F\) is in \(C^\infty((0, \infty); \mathcal{X})\), and for all \(\varepsilon > 0\) and \(m \in \mathbb{N}\), the \(m\)-th derivative of \(F\) at \(\varepsilon\) is given by \(F^{(m)}(\varepsilon) = x_\varepsilon^m\).
\end{proposition}

\begin{proof}
	The proof is done by induction. By definition, $F^{(0)}(\varepsilon) = F(\varepsilon) = x_\varepsilon$. It is therefore sufficient to prove that $F$ is continuous to initialize the induction. Let $\varepsilon_1 >0$ and $\varepsilon_2>0$. We have
	\begin{multline*}
	\Vert x_{\varepsilon_1} - x_{\varepsilon_2} \Vert_{Q_{\varepsilon_1}}^2 = Q_{\varepsilon_1}(x_{\varepsilon_1} - x_{\varepsilon_2},x_{\varepsilon_1} - x_{\varepsilon_2}) = (\varepsilon_2 - \varepsilon_1) b(x_{\varepsilon_2}, x_{\varepsilon_1} - x_{\varepsilon_2}) \\ \leq \vert \varepsilon_2 - \varepsilon_1 \vert
	 \Vert x_{\varepsilon_2} \Vert_b \Vert x_{\varepsilon_1} - x_{\varepsilon_2} \Vert_b \lesssim \frac{\vert \varepsilon_2 - \varepsilon_1 \vert}{\sqrt{\varepsilon_1 \varepsilon_2}} \Vert x_{\varepsilon_1} - x_{\varepsilon_2}\Vert_{Q_{\varepsilon_1}}.
	\end{multline*}
	By interchanging the roles of $\varepsilon_1$ and $\varepsilon_2$, we obtain
	$$
	\Vert  x_{\varepsilon_1} - x_{\varepsilon_2} \Vert_{\mathcal X} \lesssim \frac{\vert \varepsilon_2 - \varepsilon_1\vert}{\sqrt{\varepsilon_1 \varepsilon_2} \, \min(1,\sqrt{ \varepsilon_1},\sqrt {\varepsilon_2})}.
	$$
	This inequality directly implies the continuity of $F$.

	We now proceed with the induction. First of all, assuming that $\varepsilon >0 \mapsto x_\varepsilon^m \in \mathcal X$ is continuous, the same lines of reasoning gives that $\varepsilon > 0 \mapsto x^{m+1}_\varepsilon$ is also continuous. Now, setting $W^m = x^m_{\varepsilon_2} - x^m_{\varepsilon_1} - (\varepsilon_2 - \varepsilon_1) x_{\varepsilon_1}^{m+1}$, we get, for all $x$ in $\mathcal X$,
	$$
	Q_{\varepsilon_1}(W^m, x) = \left\lbrace
	\begin{array}{l}
	(\varepsilon_2 - \varepsilon_1) b(x^m_{\varepsilon_1}- x^m_{\varepsilon_2},x) - m b(W^{m-1},x), \text{ if } m\geq 1, \\
	(\varepsilon_2 - \varepsilon_1) b(x^m_{\varepsilon_1}- x^m_{\varepsilon_2},x), \text{ if } m = 0.
	\end{array}
	\right.
	$$
	The result easily follows.
\end{proof}

Up to that point, our only assumption on the bilinear form $a$ has been its positivity. We now introduce a second assumption, which, while arguably less natural, is necessary for the subsequent development and is verified in the problem under consideration.

\begin{assumption} \label{assumption_sym_of_a_on_solutions}
For all $\varepsilon>0$, for all $(m,k) \in \N \times \N$, we have $a(x_\varepsilon^m, x_\varepsilon^k) = a(x_\varepsilon^k, x_\varepsilon^m)$. Furthermore, if the Problem \eqref{prob_abstract_a} admits a (necessarily unique) solution $x_s$, we also have $a(x_\varepsilon^m, x_s) = a(x_s, x_\varepsilon^m)$.
\end{assumption}

\autoref{assumption_sym_of_a_on_solutions} corresponds to a sort of symmetry of the bilinear form $a$ on the quasi-reversibility solutions and its derivatives. In fact, it is more fundamentally an assumption on the linear form L and the bilinear form b, but it is simpler to express and utilize in the aforementioned form.

Under such assumption, we obtain the following lemma.

\begin{lemma} \label{lemma_sign_bxx}
	For all integers $m$, $p$, we have $(-1)^{m+p} b(x_\varepsilon^m,x_\varepsilon^p) \geq 0$.
\end{lemma}

\begin{proof}
	We proceed by induction. For $p\in \N$, we set the following statement:
	\begin{center}
		$P_p$ : for all $m\in \N$, $m\leq p$, $(-1)^{m+p} b(x_\varepsilon^m,x_\varepsilon^p) \geq 0$.
	\end{center}
	$P_0$ is trivially true. Now suppose that $P_p$ holds for some integer $p$. Let $m\in \N$, $m\leq p+1$. If $m=p+1$, we have $(-1)^{m+p+1} b(x_\varepsilon^{p+1}, x_\varepsilon^{m}) = \Vert x_{\varepsilon}^{p+1} \Vert_b^2 \geq 0$. If $m = p$, we have
	$$
	(-1)^{m+p+1} b(x_\varepsilon^m,x_\varepsilon^{p+1}) = - b(x_\varepsilon^p, x_\varepsilon^{p+1})  = \frac{\Vert x_\varepsilon^{p+1}\Vert_{Q_\varepsilon}^2}{p+1} \geq 0.
	$$
	Finally, if $m \in \left\lbrace 0, \ldots, p-1 \right\rbrace$, we have
	$$
	(-1)^{m+p+1} b(x_\varepsilon^m,x_\varepsilon^{p+1}) = (-1)^{m+p+2}\frac{Q_\varepsilon(x_\varepsilon^{m+1},x_\varepsilon^{p+1})}{m+1}
	= (-1)^{m+p+1} \frac{m+1}{p+1} b(x_\varepsilon^{m+1},x_\varepsilon^p).
	$$
	Here, we use  \autoref{assumption_sym_of_a_on_solutions}, which implies $Q_\varepsilon(x_\varepsilon^{m+1},x_\varepsilon^{p+1})  = Q_\varepsilon(x_\varepsilon^{p+1},x_\varepsilon^{m+1})$.
	As $m+1 \leq p$, $P_p$ implies $(-1)^{m+p+1} b(x_\varepsilon^m,x_\varepsilon^{p+1}) \geq 0$, which ends the proof.
\end{proof}

\subsection{Iterated quasi-reversibility}

We are now ready to introduce and analyze the abstract iterated quasi-reversibility method, which is the primary focus of this section.


Let \(\varepsilon > 0\). We define a sequence \((X_\varepsilon^M)_{M \in \mathbb{N}}\) in \(\mathcal{X}\) recursively as follows: set \(X_\varepsilon^0 = 0_{\mathcal{X}}\), and for each \(M \in \mathbb{N}\), let \(X_\varepsilon^{M+1}\) be the unique element of \(\mathcal{X}\) satisfying
\begin{equation} \label{eq_abstract_itQR_continuous}
Q_\varepsilon(X_\varepsilon^{M+1},x) = L(x) + \varepsilon b(X_\varepsilon^M,x), \quad \forall x \in \mathcal X.
\end{equation}
By repeated application of the Lax-Milgram theorem, the sequence is well-defined. A straightforward induction further shows that
$$
X_\varepsilon^M = \sum_{p=0}^{M-1} (-1)^p \frac{\varepsilon^p}{p!} x_\varepsilon^p.
$$
As a consequence, \autoref{assumption_sym_of_a_on_solutions} implies that for all $\varepsilon>0$, for all $(M,K) \in \N \times \N$, $a(X_\varepsilon^M, X_\varepsilon^K) = a(X_\varepsilon^K, X_\varepsilon^M)$, and if there exists $x_s$ satisfying \eqref{prob_abstract_a}, $a(X_\varepsilon^M,x_s) = a(x_s, X_\varepsilon^M)$.

Let us start with some a priori estimates. From now on, we fix $\varepsilon >0$.

\begin{proposition} \label{prop_aprioriestimates_bis}
	For all $M\in \N$, we have
	\begin{enumerate}
		\item $\Vert X_\varepsilon^{M} \Vert_{b} \leq \Vert X_\varepsilon^{M+1} \Vert_{b} $.
		\item $\Vert X_\varepsilon^M \Vert_a \leq c_a$.
	\end{enumerate}
\end{proposition}

\begin{proof}
	\begin{enumerate}
		\item 
		Given that $X_\varepsilon^{M+1} = X_\varepsilon^M + (-1)^M \frac{\varepsilon^M}{M!} x_\varepsilon^M$, we obtain
		$$
		\Vert X_\varepsilon^{M+1} \Vert_b^2 = \Vert X_\varepsilon^M \Vert_b^2 + \frac{\varepsilon^{2M}}{(M!)^2} \Vert x_\varepsilon^M \Vert_b^2
		 + 2 (-1)^M \frac{\varepsilon^M}{M! } b\left( x_\varepsilon^M, X_\varepsilon^M \right).
		$$
		By applying \autoref{lemma_sign_bxx}, we obtain
		$$
		(-1)^M \frac{\varepsilon^M}{M! } b\left( x_\varepsilon^M, X_\varepsilon^M \right) = \sum_{p=0}^M \frac{\varepsilon^{M+p}}{M!p!} (-1)^{p+M} b(x_\varepsilon^M , x_\varepsilon^p) \geq 0,
		$$
		which gives the estimate.
		\item It is a slight adaptation of \cite[Proposition~8, item~(b)]{darde2016}: obviously the estimate holds for $M=0$. Now,
		as
		$$
		\Vert X_\varepsilon^{M+1}\Vert_a^2 + \varepsilon \Vert X_\varepsilon^{M+1}\Vert_b^2 = L(X_\varepsilon^{M+1}) + \varepsilon b(X_\varepsilon^M,X_\varepsilon^{M+1}),
		$$
		we have
		$$
		\Vert X_\varepsilon^{M+1}\Vert_a^2 \leq c_a \Vert X_\varepsilon^{M+1} \Vert_a + \varepsilon \left( b(X_\varepsilon^M,X_\varepsilon^{M+1}) - \Vert X_\varepsilon^{M+1} \Vert_b^2 \right).
		$$
		The first estimate gives $b(X_\varepsilon^M,X_\varepsilon^{M+1}) - \Vert X_\varepsilon^{M+1} \Vert_b^2 \leq 0$, which ends the proof.
	\end{enumerate}
\end{proof}

\begin{proposition} \label{prop_aprioriestimates_ter}
	For all $M\in \N$, we have $\Vert X_\varepsilon^M \Vert_b \leq c_a \sqrt{\frac{2\,M}{\varepsilon}}.$
\end{proposition}

\begin{proof}
	We proceed as in \cite[Proposition 9]{darde2016}. To begin with, we note that if $c_a = 0$, then $X_\varepsilon^M = 0$ for all $M$, hence the result holds. Therefore we focus on the case $c_a>0$.

	Consider a sequence   $(x_p)_{p\in \N_*}$ of positive real numbers, where $x_1=c_a > 0$ and, for each
	$p \in \N_*$, $x_{p+1}$ is defined as the positive root of the polynomial  $P_p(x) = x^2 - x x_p - x_1^2$.
	As $P_p(0) = -x_1^2$, the sequence is well-defined. We claim that $x_p < \sqrt{2p}\, x_1$. Indeed, suppose it is true for some integer $p$. Then, as
	\begin{multline*}
	P_p\left(\sqrt{2 (p+1)}x_1\right) = 2 (p+1) x_1^2 - \sqrt{2 (p+1)} x_1 x_p - x_1^2
	> 2 (p+1)x_1^2 - \sqrt{2(p+1)} \sqrt{2p} x_1^2 - x_1^2  \\ = \left(\sqrt{p +1} - \sqrt{p}\right)^2 x_1^2 > 0,
	\end{multline*}
	we have $x_{p+1} < \sqrt{2(p+1)}x_1 $.
	\begin{figure}[h!]
	\begin{center}
	\begin{tikzpicture}[scale=0.5]
	\draw [domain=0:15.5,samples=500, smooth, color=blue,thick] plot(\x,{sqrt(2*\x)}) node[above left,scale=0.7]{$\sqrt{2x}x_1$};
	\draw [->,thick] (0,0) -- (15.5,0);
	\draw [->,thick] (0,0) -- (0,6);
	\pgfmathsetmacro{\XX}{1};
	\foreach \i in {1,...,15} {
	\coordinate (A) at (\i,\XX);
	\draw (A) node{$\bullet$}  node[below right, scale=0.7]{$x_{\i}$};
	\draw (\i,0) node{+};
	\pgfmathsetmacro{\XXnext}{\XX + sqrt(\XX*\XX+4))/2};
	\global\let\XX=\XXnext;
	}
	\end{tikzpicture}
	\end{center}
	\caption{The terms of the sequence $(x_p)_{p\in \N}$  lie below the curve $f(x) = \sqrt{2x} \, x_1$.}
	\end{figure}
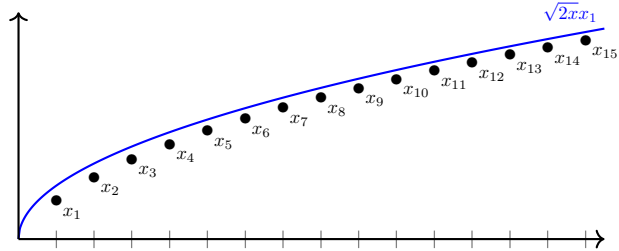

	We now return to the sequence \((X_\varepsilon^M)_{M \in \mathbb{N}}\). Since \(X_\varepsilon^0 = 0\), the result trivially holds for \(M = 0\). For any \(M \in \mathbb{N}\), we have
\[
\| X_\varepsilon^{M+1} \|_a^2 + \varepsilon \| X_\varepsilon^{M+1} \|_b^2 = L(X_\varepsilon^{M+1}) + \varepsilon \, b(X_\varepsilon^M, X_\varepsilon^{M+1}).
\]
This implies, in particular,
\[
\varepsilon \| X_\varepsilon^{M+1} \|_b^2 \leq c_a \| X_\varepsilon^{M+1} \|_a + \varepsilon \| X_\varepsilon^M \|_b \| X_\varepsilon^{M+1} \|_b \leq c_a^2 + \varepsilon \| X_\varepsilon^M \|_b^2.
\]
Let \(w_p = \sqrt{\varepsilon} \| X_\varepsilon^p \|_b\). Then \(w_p\) satisfies
\[
w_{p+1}^2 - w_{p+1} w_p - c_a^2 = w_{p+1}^2 - w_{p+1} w_p - x_1^2 \leq 0.
\]
We claim that \(w_p \leq x_p\) for all \(p \in \mathbb{N}_*\). To prove this, suppose the inequality holds for some \(p \in \mathbb{N}_*\). Then,
\[
P_p(w_{p+1}) = w_{p+1}^2 - w_{p+1} x_p - x_1^2 \leq w_{p+1}^2 - w_{p+1} x_p - x_1^2 \leq 0,
\]
which implies \(w_{p+1} \leq x_{p+1}\). Since \(w_1 = \sqrt{\varepsilon} \| x_\varepsilon \|_b \leq c_a = x_1\), the result follows by induction.

\end{proof}

The following Proposition can be easily deduced by adapting the proof of Proposition 10 in \cite{darde2016}.

\begin{proposition} \label{prop_aprioriestimates_abstractsetting}
Suppose Problem \eqref{prob_abstract_a} admits a solution $x_s$. Then we have
\begin{enumerate}
\item $\Vert X_\varepsilon^{M+1} - x_s \Vert_b \leq \Vert X_\varepsilon^M - x_s \Vert_b$.
\item $\Vert X_\varepsilon^M \Vert_b \leq \Vert x_s \Vert_b$.
\end{enumerate}
\end{proposition}

We can now state the following convergence result for the abstract iterated quasi-reversibility method:

\begin{proposition} \label{prop_residualestimates_abstractsetting}
Suppose Problem \eqref{prob_abstract_a} admits a solution $x_s$. Set $\varepsilon>0$. For all $M\in \N$,  we have
$$
 \Vert X_\varepsilon^{M+1} - x_s \Vert_a \leq \Vert X_\varepsilon^M - x_s \Vert_a \leq \sqrt{\frac{\varepsilon}{M}} \Vert x_s \Vert_b.
$$
As a consequence, we have $X_\varepsilon^M \xrightarrow{M\rightarrow \infty} x_s$.
\end{proposition}

Let us first make some comments on this Proposition before proceeding to the proof. First, the right-hand inequality, which provides a convergence estimate for $\Vert X_\varepsilon^M - x_s \Vert_a$, is novel: while all the necessary tools were available in \cite{darde2016}, it was not derived there.  Secondly, the left-hand inequality corresponds to inequality (d) of \cite[Proposition 8]{darde2016}. However, the proof of this inequality in \cite{darde2016} contains an error, and we will use this opportunity to rectify it.

\begin{proof}[Proof of \autoref{prop_residualestimates_abstractsetting}]
 We start by proving the left-hand inequality of the Proposition.

Set $\varepsilon>0$. Recalling that for all $M\in \N$, $\displaystyle X_\varepsilon^M = \sum_{p=0}^{M-1} (-1)^p \frac{\varepsilon^p}{p!}  x_\varepsilon^p$, and  \autoref{assumption_sym_of_a_on_solutions}, we see that
$$
\Vert X_\varepsilon^{M+1} - x_s \Vert_a^2 = \Vert X_\varepsilon^M - x_s \Vert_a^2 + 2\, (-1)^M \frac{\varepsilon^M}{M!} a(X_\varepsilon^M - x_s, x_\varepsilon^M)
+ \frac{\varepsilon^{2M}}{(M!)^2} \Vert x_\varepsilon^M \Vert_a^2.
$$
We first apply this identity for $M= 0$, which gives
\begin{multline*}
\Vert x_\varepsilon - x_s \Vert_a^2 = \Vert x_s \Vert_a^2  - 2 a(x_s,x_\varepsilon) + a(x_\varepsilon, x_\varepsilon)
= \Vert x_a\Vert_a^2 -2 L(x_\varepsilon) + a(x_\varepsilon,x_\varepsilon) \\ = \Vert x_s\Vert_a^2 - 2 \varepsilon \Vert x_\varepsilon\Vert_b^2 - \Vert x_\varepsilon\Vert_a^2
\leq \Vert x_s\Vert_a^2.
\end{multline*}
We now deal with the case $M\geq 1$. We note that $\displaystyle X_\varepsilon^M - x_s = x_\varepsilon^M - x_s + \sum_{p=1}^{M-1} (-1)^p \frac{\varepsilon^p}{p!}x_\varepsilon^p$.
Using the equation defining $x_\varepsilon^p$, we obtain
$$
a(x_\varepsilon^p, x_\varepsilon^M) = -\varepsilon b(x_\varepsilon^p,x_\varepsilon^M) - p b(x_\varepsilon^{p-1},x_\varepsilon^M).
$$
Multiplying this identity by $(-1)^p \frac{\varepsilon^p}{p!}$ and summing for $p$ in $\left\lbrace 1, \ldots, M-1\right\rbrace$ gives
$$
a\left( \sum_{p=1}^{M-1} (-1)^p \frac{\varepsilon^p}{p!}  x_\varepsilon^p, x_\varepsilon^M \right) = (-1)^M \frac{\varepsilon^M}{(M-1)!} b(x_\varepsilon^{M-1},x_\varepsilon^M) + \varepsilon b(x_\varepsilon,x_\varepsilon^M).
$$
As by definition of $x_\varepsilon$ and $x_s$, we have $a(x_\varepsilon - x_s, x_\varepsilon^M) = -\varepsilon b(x_\varepsilon,x_\varepsilon^M)$, we finally obtain  that
$a(X_\varepsilon^M - x_s, x_\varepsilon^M) = (-1)^M \frac{\varepsilon^M}{(M-1)!} b(x_\varepsilon^{M-1},x_\varepsilon^M) $, which yields
$$
\Vert X_\varepsilon^{M+1} - x_s \Vert_a^2 = \Vert X_\varepsilon^M - x_s \Vert_a^2 + 2 \frac{\varepsilon^{2M}}{M! (M-1)!} b(x_\varepsilon^{M-1},x_\varepsilon^M) + \frac{\varepsilon^{2M}}{(M!)^2} \Vert x_\varepsilon^M \Vert_a^2.
$$
Finally, the equation defining $x_\varepsilon^M$ gives $M b(x_\varepsilon^{M-1},x_\varepsilon^M) = -\Vert x_\varepsilon^M \Vert_a^2 - \varepsilon \Vert x_\varepsilon^M\Vert_b^2$, leading to
$$
\Vert X_\varepsilon^{M+1} - x_s \Vert_a^2 = \Vert X_\varepsilon^M - x_s \Vert_a^2 - 2 \frac{\varepsilon^{2M+1}}{(M!)^2} \Vert x_\varepsilon^M \Vert_b^2 - \frac{\varepsilon^{2M}}{(M!)^2} \Vert x_\varepsilon^M \Vert_a^2.
$$
As a result, the left-hand inequality of  \autoref{prop_residualestimates_abstractsetting} is proven.
\smallskip

Using the equations defining $X_\varepsilon^{M+1}$ and $x_s$ again, and the first inequality of  \autoref{prop_aprioriestimates_abstractsetting}, we obtain, for $p \geq 1$,
\begin{multline*}
\Vert X_\varepsilon^{p+1} - x_s \Vert_a^2 \leq \Vert X_\varepsilon^{p+1} -x_s \Vert_{\mathcal{Q}_\varepsilon}^2 = \varepsilon b(X_\varepsilon^p - x_s, X_\varepsilon^{p+1} - x_s)
 \leq \varepsilon \Vert X_\varepsilon^p - x_s \Vert_b \Vert X_\varepsilon^{p+1} - x_s \Vert_b \\
\leq \varepsilon \Vert X_\varepsilon^p - x_s \Vert_b^2 = \varepsilon b(X_\varepsilon^p - x_s, X_\varepsilon^{p-1} - x_s) - \Vert X_\varepsilon^p - x_s \Vert_a^2.
\end{multline*}
In consequence, we obtain
$$
\Vert X_\varepsilon^{p+1} - x_s \Vert_a^2  + \Vert X_\varepsilon^p - x_s \Vert_a^2 \leq \varepsilon b(X_\varepsilon^p - x_s, X_\varepsilon^{p-1} - x_s) .
$$
Reiterating this process starting from $p = M$ to $p=1$ gives
$$
\sum_{p=1}^{M+1} \Vert X_\varepsilon^p - x_s \Vert_a^2  \leq \varepsilon b(X_\varepsilon^1 - x_s,X_\varepsilon^0 - x_s) = - \varepsilon b(x_\varepsilon - x_s, x_s) \leq \varepsilon \Vert x_s \Vert_b^2.
$$
Consequently, the series $\sum_{p\geq 1} \Vert X_\varepsilon^p - x_s \Vert_a^2$ converges, and its sum is smaller than $\varepsilon \Vert x_s \Vert_b^2$. The right-hand inequality of  \autoref{prop_residualestimates_abstractsetting} is then a direct consequence of  \autoref{lemma_estimate_series}.
\smallskip

Finally, we prove the convergence of the sequence $X_\varepsilon^M$ to $x_s$. We note that by  \autoref{prop_aprioriestimates_bis}, inequality 2, and  \autoref{prop_aprioriestimates_abstractsetting}, inequality 2, the sequence $X_\varepsilon^M$ is uniformly bounded in $\mathcal{X}$. As a consequence, there exists an increasing sequence of non-negative integers $(M_n)_{n\in\N}$ and $x_\infty$ in $\mathcal X$ such that $X_\varepsilon^{M_n}$ converges weakly to $x_\infty$ as $n$ goes to infinity. In particular, for all $x$ in $\mathcal X$, $a(X_\varepsilon^{M_n}, x ) \xrightarrow{n\rightarrow \infty} a(x_\infty,x)$. But  \autoref{prop_residualestimates_abstractsetting} immediately implies that
$a(X_\varepsilon^{M_n}, x) \xrightarrow{n\rightarrow \infty} a(x_s,x)$, which leads thanks to  \autoref{assumption_identifiability} to $x_\infty = x_s$.

Then, combining  \autoref{prop_aprioriestimates_abstractsetting} and \autoref{prop_residualestimates_abstractsetting} gives
$$
\Vert X_\varepsilon^{M_n} - x_s \Vert_{\mathcal{Q}_1}^2 = \Vert X_\varepsilon^{M_n} - x_s \Vert_a^2 + \varepsilon \Vert X_\varepsilon^{M_n} - x_s \Vert_b^2
\leq \frac{\varepsilon}{M_n} \Vert x_s\Vert_b^2 + 2\Vert x_s \Vert_b^2 - 2 b(x_s, X_\varepsilon^{M_n}) \xrightarrow{n\rightarrow \infty} 0.
$$
Hence $X_\varepsilon^{M_n}$ strongly converges to $x_s$. We obtain the convergence of the entire sequence by arguing by contradiction.
\end{proof}

\subsubsection*{Galerkin Approximation}

Let \(\mathcal{X}_h\) be a closed subspace of \(\mathcal{X}\). The Galerkin approximation of the abstract iterated quasi-reversibility method consists in defining the sequence \((X_{\varepsilon,h}^M)_{M \in \mathbb{N}}\) in \(\mathcal{X}_h^\mathbb{N}\), which approximates the sequence \((X_\varepsilon^M)_{M \in \mathbb{N}}\). We set \(X_{\varepsilon,h}^0 = X_\varepsilon^0 = 0_{\mathcal{X}}\), and for each integer \(M\), we define \(X_{\varepsilon,h}^{M+1}\) as the unique element in \(\mathcal{X}_h\) satisfying
\[
\mathcal{Q}_\varepsilon\left(X_{\varepsilon,h}^{M+1}, x_h\right) = L(x_h) + \varepsilon \, b\left(X_{\varepsilon,h}^M, x_h\right), \quad \forall x_h \in \mathcal{X}_h.
\]
Since \((\mathcal{X}_h, \|\cdot\|_{\mathcal{X}})\) is a Hilbert space, the Lax-Milgram theorem and a straightforward induction argument ensure that the sequence \((X_{\varepsilon,h}^M)_{M \in \mathbb{N}}\) is well-defined.

Our key result concerning the Galerkin approximation is as follows:

\begin{proposition} \label{prop_abstract_Galerkin_estimation}
There exists a constant $c\in \R$ such that, for all $\varepsilon>0$, for all $M\in \N$,
$$
\Vert X_{\varepsilon,h}^M - X_\varepsilon^M \Vert_{\mathcal{Q}_\varepsilon} \leq c  \frac{\max(1,\varepsilon)}{\sqrt \varepsilon} \sum_{p=0}^M \inf_{x_h \in \mathcal{X}_h} \Vert X_\varepsilon^p - x_h \Vert_{\mathcal X}.
$$
\end{proposition}

\begin{proof}
Set $\varepsilon>0$. We define the auxiliary sequence $(X_{\varepsilon,h}^{M,\star})_{M\in \N}$ of elements of $\mathcal{X}_h$ by $X_{\varepsilon,h}^{0,\star} = 0_\mathcal{X}$, and for all $M\in\N$, $X_{\varepsilon,h}^{M+1,\star}$ as the unique element of $\mathcal{X}_h$ solution of
$$
\mathcal{Q}_\varepsilon\left(X_{\varepsilon,h}^{M+1,\star},x_h\right) = L(x_h) + \varepsilon b\left(X_\varepsilon^M,x_h\right), \quad \forall x_h \in \mathcal{X}_h.
$$
By Galerkin orthogonality (see, \textit{e.g.}, \cite[Lemma 2.17]{ernguermond2004}), the difference $ X_{\varepsilon,h}^{M,\star} - X_\varepsilon^M$ satisfies
$$
\mathcal{Q}_\varepsilon(X_{\varepsilon,h}^{M,\star} - X_{\varepsilon}^{M} , x_h) = 0, \quad \forall x_h \in \mathcal X_h,
$$
which classically implies
$$
\Vert X_{\varepsilon,h}^{M,\star} - X_{\varepsilon}^{M} \Vert^2_{\mathcal{Q}_\varepsilon} = \mathcal{Q}_\varepsilon \left(X_{\varepsilon,h}^{M,\star} - X_{\varepsilon}^{M} , x_h - X_\varepsilon^M\right), \quad \forall x_h \in \mathcal{X}_h.
$$
As for all $(x,\tilde x)$ in $\mathcal{X} \times \mathcal{X}$,
\begin{multline} \label{eq_estimate_Qe(x,tildex)_nonsym}
\mathcal{Q}_\varepsilon(x,\tilde x) = a(x,\tilde x) + \varepsilon b(x,\tilde x) \leq \Vert a \Vert \, \Vert x\Vert_{\mathcal{X}} \, \Vert \tilde x \Vert_{\mathcal X} + \varepsilon \Vert x \Vert_b\,  \Vert \tilde x \Vert_b \lesssim \Vert x \Vert_{\mathcal{Q}_1} \, \Vert \tilde x \Vert_{\mathcal{X}} + \sqrt \varepsilon \Vert x \Vert_{\mathcal{Q}_\varepsilon}\, \Vert \tilde x \Vert_{\mathcal{Q}_1} \\
\lesssim \max \left(1, \frac{1}{\sqrt{\varepsilon}}\right) \Vert x \Vert_{\mathcal{Q}_\varepsilon} \, \Vert \tilde x \Vert_{\mathcal{X}} + \sqrt \varepsilon \Vert x \Vert_{\mathcal{Q}_\varepsilon}\, \Vert \tilde x \Vert_{\mathcal{X}} \lesssim \frac{\max(1,\varepsilon)}{\sqrt \varepsilon} \Vert x \Vert_{\mathcal{Q}_\varepsilon} \, \Vert \tilde x \Vert_{\mathcal X},
\end{multline}
we deduce
\begin{equation} \label{eq_XepsMminusXepshMstar_estimate}
\Vert X_{\varepsilon,h}^{M,\star} - X_{\varepsilon}^{M} \Vert_{\mathcal{Q}_\varepsilon} \lesssim \frac{\max(1,\varepsilon)}{\sqrt \varepsilon} \inf_{x_h \in \mathcal{X}_h} \Vert X_\varepsilon^M - x_h \Vert_{\mathcal{X}}.
\end{equation}
In addition, using the definition of both $X_{\varepsilon,h}^M$ and $X_{\varepsilon,h}^{M,\star}$, we have, for all integer $M$,
\begin{multline*}
	\Vert X_{\varepsilon,h}^{M+1,\star} - X_{\varepsilon,h}^{M+1} \Vert_{\mathcal{Q}\varepsilon}^2 =  \varepsilon \, b(X_\varepsilon^{M} - X_{\varepsilon,h}^{M}, X_{\varepsilon,h}^{M+1,\star} - X_{\varepsilon,h}^{M+1})
 \\	= \varepsilon b(X_\varepsilon^M - X_{\varepsilon,h}^{M,\star} , X_{\varepsilon,h}^{M+1,\star} - X_{\varepsilon,h}^{M+1} )  + \varepsilon b(X_{\varepsilon,h}^{M,\star} - X_{\varepsilon,h}^M,X_{\varepsilon,h}^{M+1,\star} - X_{\varepsilon,h}^{M+1}),
\end{multline*}
and we conclude that
$$
\Vert X_{\varepsilon,h}^{M+1,\star} - X_{\varepsilon,h}^{M+1} \Vert_{\mathcal{Q}\varepsilon} \leq \Vert X_\varepsilon^M - X_{\varepsilon,h}^{M,\star}\Vert_{\mathcal{Q}_\varepsilon}
+ \Vert X_{\varepsilon,h}^{M,\star} - X_{\varepsilon,h}^M \Vert_{\mathcal{Q}_\varepsilon}.
$$
The result follows from a simple induction.
\end{proof}

\subsubsection*{Data perturbation}

We now focus on the effects of data perturbation, analyzing this at the discretization level---that is, within the space \( X_h \). We introduce a family of linear continuous forms $(L^\delta : \mathcal X \mapsto \R)_{\delta \geq 0}$, such that $\Vert L^\delta - L \Vert \leq  c \delta$ for some fixed constant $c$. Here, $\delta$ is understood as a \emph{noise amplitude}, and
$L^\delta$ as the \emph{perturbed data}, $L^0 = L$ being the \emph{exact or unperturbed data}.

We introduce the discretized iterated quasi-reversibility sequence of solutions, denoted by \(\big(X_{\varepsilon,h}^{M,\delta}\big)_{M \in \mathbb{N}}\), defined by setting \(X_{\varepsilon,h}^{0,\delta} = 0\) and, for each \(M \in \mathbb{N}\), defining \(X_{\varepsilon,h}^{M+1,\delta} \in \mathcal{X}_h\) such that
\begin{equation} \label{eq_def_itQRabs_discr_noise}
\mathcal{Q}_\varepsilon\big(X_{\varepsilon,h}^{M+1,\delta},x_h \big)  = L^\delta(x_h) + \varepsilon b\big( X_{\varepsilon,h}^{M,\delta},x_h\big),\quad \forall x_h \in \mathcal X_h.
\end{equation}
Note that the sequence \((X_{\varepsilon,h}^{M,0})_{M \in \mathbb{N}}\) coincides with the sequence \((X_{\varepsilon,h}^M)_{M \in \mathbb{N}}\) introduced in the previous section.

\begin{proposition} \label{prop_abstracterrorestimatenoise}
There exists $c\in \R$,  such that for all $\mathcal X_h$ closed subspace of $\mathcal X$, for all $\varepsilon>0$, for all $M\in \N$ and for all $\delta \geq 0$,
$$
\Vert X_{\varepsilon,h}^{M,\delta} - X_{\varepsilon,h}^M \Vert_{\mathcal{Q}_\varepsilon} \leq c \sqrt M \delta.
$$
\end{proposition}

\begin{proof}
By linearity, $W_{\varepsilon,h}^{M,\delta} = X_{\varepsilon,h}^{M,\delta} - X_{\varepsilon,h}^M \in \mathcal{X}_h$ satisfies
$$
\mathcal{Q}_\varepsilon\big(W_{\varepsilon,h}^{M+1,\delta},x_h \big)  = \mathcal{L}^\delta(x_h) + \varepsilon b\big( W_{\varepsilon,h}^{M,\delta},x_h\big), \quad \forall x_h \in \mathcal X_h,
$$
with $\mathcal{L}^\delta = L^\delta - L$. The result follows directly by applying the second estimate of \autoref{prop_aprioriestimates_bis} together with the estimate of   \autoref{prop_aprioriestimates_ter}, since these results, obtained at the continuous level (\textit{i.e.} in $\mathcal{X}$), are valid at the discrete level (\textit{i.e.} in $\mathcal X_h$).
\end{proof}

\subsubsection*{Final estimate}

We are now in position to state a global convergence estimate for the abstract quasi-reversibility method, which reads as follows.

\begin{theorem}  Let $\mathcal{X}$ be a Hilbert space, $a : \mathcal X \times \mathcal X$ a continuous, positive bilinear form and $L : \mathcal{X} \mapsto \R$ a continuous linear form.
Assume \eqref{prob_abstract_a} admits a unique solution $x_s$, and that Assumptions \eqref{assumption_identifiability}
and \eqref{assumption_sym_of_a_on_solutions} are verified.



There exists a constant $c\in \R$ such that, for all $\varepsilon>0$, $M\in \N$, $\mathcal X_h$ closed subspace of $\mathcal X$ and $\delta \geq 0$,
$$
\Vert X_{\varepsilon,h}^{M,\delta} - x_s \Vert_{\mathcal X} \leq c \left( \sqrt{\frac{M}{\min(1,\varepsilon)}} \delta + \max \left(\frac{1}{\varepsilon},\sqrt \varepsilon \right) \sum_{p=0}^M
\inf_{x_h \in \mathcal{X}_h} \Vert X_\varepsilon^p - x_h \Vert_\mathcal{X} + \Vert X_\varepsilon^M - x_s \Vert_{\mathcal X} \right),
$$
where $X_\varepsilon^M$ is defined by \eqref{eq_abstract_itQR_continuous}, and $X_{\varepsilon,h}^{M,\delta}$ is defined by \eqref{eq_def_itQRabs_discr_noise}.
\end{theorem}

\begin{proof}
By the triangular inequality, we obtain
$$
\Vert X_{\varepsilon,h}^{M,\delta} - x_s \Vert_{\mathcal X} \leq \Vert X_{\varepsilon,h}^{M,\delta} - X_{\varepsilon,h}^M \Vert_{\mathcal{X}}
+ \Vert X_{\varepsilon,h}^M - X_\varepsilon^M \Vert_{\mathcal{X}} + \Vert X_\varepsilon^M - x_s \Vert_{\mathcal{X}}.
$$
As the  $\mathcal{X}$-norm is equivalent to the $\mathcal{Q}_1$-norm, we have
$$
\Vert X_{\varepsilon,h}^{M,\delta} - X_{\varepsilon,h}^M \Vert_{\mathcal{X}}
+ \Vert X_{\varepsilon,h}^M - X_\varepsilon^M \Vert_{\mathcal{X}} \lesssim \frac{1}{\min(1,\sqrt{\varepsilon})}
\left(\Vert X_{\varepsilon,h}^{M,\delta} - X_{\varepsilon,h}^M \Vert_{\mathcal{Q}_\varepsilon}
+ \Vert X_{\varepsilon,h}^M - X_\varepsilon^M \Vert_{\mathcal{Q}_\varepsilon} \right).
$$
The result is then a direct application of  \autoref{prop_abstract_Galerkin_estimation} and  \autoref{prop_abstracterrorestimatenoise}.
\end{proof}

\section{Data assimilation problem: continuous quasi-reversibility formulation and regularity results} \label{sec_contQRregres}

\subsection{From the abstract setting to the data assimilation problem}

Our goal is now to apply the abstract method developed in the previous Section to the data assimilation problem \eqref{eq_dataassimpb}. To do so, we set $\mathcal{X} = H^1(\Omega) \times H^1_0(\Omega)$, and define
$$
L : (v,\mu) \in \mathcal X \mapsto \int_\Omega \nabla \lambda^\delta \cdot \nabla v dx + \int_\omega g^\delta v  dx \in \R,
$$
and
$$
a : \big( (u,\lambda), (v,\mu) \big) \in \mathcal X \times \mathcal X \mapsto \int_\Omega \nabla \lambda \cdot \nabla v dx + \int_\omega u v dx + \int_\Omega \nabla \lambda \cdot \nabla \mu dx - \int_\Omega \nabla u \cdot \nabla \mu dx.
$$
Clearly, $L$ is a linear continuous form on $\mathcal X$, and $a$ is a continuous positive bilinear form on $\mathcal X$. Furthermore, a straightforward computation shows that
$\vert L(u,\lambda) \vert \leq \sqrt{\Vert \nabla \lambda^\delta \Vert_{L^2(\Omega)}^2 + \Vert g^\delta \Vert_{L^2(\omega)}^2} \Vert u,\lambda \Vert_a$.

We now show that the bilinear form $a$ satisfies  \autoref{assumption_identifiability}, and a solution of Problem \eqref{eq_dataassimpb}, if it exists, is necessarily a solution of Problem \eqref{prob_abstract_a}:
\begin{proposition} \label{prop_asatisfiesA1}
The bilinear form \( a \) satisfies {\autoref{assumption_identifiability}}. Moreover, assume that \textbf{Problem~\eqref{eq_dataassimpb}} admits a (necessarily unique) solution \(\tilde{u}\). Define \(\tilde{\lambda}\) as its orthogonal projection onto \( H^1_0(\Omega) \), characterized by the identity:
\[
\int_\Omega \nabla \tilde{\lambda} \cdot \nabla \mu \, \mathrm{d}x = \int_\Omega \nabla \tilde{u} \cdot \nabla \mu \, \mathrm{d}x, \quad \forall \mu \in H^1_0(\Omega).
\]
Then, the pair \((\tilde{u}, \tilde{\lambda})\) solves \textbf{Problem~\eqref{prob_abstract_a}}.
%
\end{proposition}

\begin{proof}
Let $(u,\lambda) \in \mathcal{X}$ be such that $a\big( (u,\lambda), (v,\mu) \big) = 0$ for all $(v,\mu) \in \mathcal X$. Choosing $(v,\mu) = (u,\lambda)$ immediately implies $u = 0$ in $\omega$ and $\lambda = 0$ in $\Omega$. Now, choosing $(v,\mu) = (0,\phi)$ with $\phi \in C^\infty_c(\Omega)$, we get
$$
\int_\Omega \nabla u \cdot \nabla \phi\, dx = 0, \quad \forall \varphi \in C^\infty_c(\Omega).
$$
In consequence, $\Delta u = 0$ in $\Omega$. As $u = 0$ in $\omega$,  unique continuation  immediately implies $u = 0$ in $\Omega$. In conclusion, \autoref{assumption_identifiability} holds true.

Now, for a solution $\tilde u$ of Problem \eqref{eq_dataassimpb} and $\tilde \lambda$ its orthogonal projection on $H^1_0(\Omega)$, we have
$$
a\big((\tilde u,\tilde \lambda),(v,\mu) \big) = \int_\Omega \nabla \tilde \lambda\cdot \nabla v dx + \int_\omega \tilde  u v dx
 = \int_\Omega \nabla \tilde \lambda\cdot \nabla v dx + \int_\omega g^\delta v dx, \quad \forall (v,\mu) \in \mathcal X.
$$
But for all $\phi \in C^\infty_c(\Omega)$, we have by definition of $u$ and $\lambda^\delta$,
$$
\int_\Omega \nabla \tilde u \cdot \nabla \phi\, dx = - \langle \Delta \tilde u , \phi \rangle  = - \langle f^\delta + \nabla \cdot \mathbf{F}^\delta , \phi \rangle
 = \int_\Omega (\mathbf{F}^\delta \cdot \nabla \phi  - f^\delta \phi ) dx = \int_\Omega \nabla \lambda^\delta \cdot \nabla \phi\, dx.
$$
The density of $C^\infty_c(\Omega)$ in $H^1_0(\Omega)$ and the definition of $\tilde \lambda$ give $\tilde \lambda = \lambda^\delta$, hence $(\tilde u,\tilde \lambda)$ is solution of \eqref{prob_abstract_a}.
\end{proof}

We now set
$$
b : \big((u,\lambda),(v,\mu) \big) \in \mathcal X \times \mathcal X \mapsto \int_\Omega \nabla u \cdot \nabla v\, dx.
$$
Clearly, $b$ is a bilinear continuous positive and symmetric form on $\mathcal X$.

\begin{proposition}
The application $\mathcal N : (u,\lambda) \in \mathcal{X} \mapsto \sqrt{\Vert (u,\lambda) \Vert_a^2 + \Vert (u,\lambda) \Vert_b^2}$ is a norm on $\mathcal{X}$, equivalent to the standard $H^1\times H_0^1$ norm.
\end{proposition}

\begin{proof}
This is straightforward since for all $(u,\lambda) \in \mathcal{X}$,
$$
\Vert (u,\lambda) \Vert_a^2 + \Vert (u,\lambda) \Vert_b^2  = \int_\Omega \vert \nabla \lambda\vert^2 dx + \int_\omega \vert u \vert^2 dx + \int_\Omega \vert \nabla u \vert^2 dx.
$$
\end{proof}

As a result, the abstract setting developed in  \autoref{sec_abstractQR} for the \emph{standard quasi-reversibility method} applies.
In our context, the quasi-reversibility problem is formulated as follows: for $\varepsilon>0$, find $(u_\varepsilon,\lambda_\varepsilon) \in H^1(\Omega) \times H_0^1(\Omega)$ such that, for all $(v,\mu) \in H^1(\Omega) \times H^1_0(\Omega)$,
\begin{multline*}
\int_\Omega \nabla \lambda_\varepsilon \cdot \nabla v dx + \int_\omega u_\varepsilon v dx + \int_\Omega \nabla \lambda_\varepsilon\cdot \nabla \mu dx - \int_\Omega \nabla u_\varepsilon\cdot \nabla \mu dx + \varepsilon \int_\Omega \nabla u_\varepsilon\cdot \nabla v dx
\\ = \int_\Omega  \nabla \lambda^\delta \cdot \nabla v dx + \int_\omega g^\delta v dx.
\end{multline*}
Alternatively, the quasi-reversibility equation can be expressed in the equivalent system
$$
\begin{cases}
\displaystyle \int_\Omega \nabla \lambda_\varepsilon \cdot \nabla v dx + \int_\omega u_\varepsilon v dx + \varepsilon \int_\Omega \nabla u_\varepsilon\cdot \nabla v dx
 = \int_\Omega  \nabla \lambda^\delta \cdot \nabla v dx + \int_\omega g^\delta v dx, \quad \forall v \in H^1(\Omega), \\
 \\
 \displaystyle \int_\Omega \nabla \lambda_\varepsilon\cdot \nabla \mu dx - \int_\Omega \nabla u_\varepsilon\cdot \nabla \mu dx = 0, \quad \forall \mu \in H^1_0(\Omega).
\end{cases}
$$
In this form, it becomes clear that $\lambda_\varepsilon$ is the orthogonal projection of $u_\varepsilon$ onto $H^1_0(\Omega)$.


We now turn our attention to the sequence $(x_\varepsilon^m)_{m\in \N}$, defined by \eqref{eq_definition_abstractderivatives}. By definition, $x_\varepsilon^0 = x_\varepsilon = (u_\varepsilon,\lambda_\varepsilon)$.
For $m\in \N_*$, we set $x_\varepsilon^m = (v_\varepsilon^m, \mu_\varepsilon^m) \in H^1(\Omega) \times H_0^1(\Omega)$. The sequence satisfies
\begin{multline*}
\int_\Omega \nabla \mu_\varepsilon^{m+1} \cdot \nabla v dx + \int_\omega v_\varepsilon^{m+1} v dx + \int_\Omega \nabla \mu_\varepsilon^{m+1}\cdot \nabla \mu dx - \int_\Omega \nabla u_\varepsilon^{m+1} \cdot \nabla \mu dx + \varepsilon \int_\Omega \nabla u_\varepsilon^{m+1}\cdot \nabla v dx
\\ = -(m+1) \int_\Omega \nabla v_\varepsilon^m \cdot \nabla v dx, \quad \forall (v,\mu) \in H^1(\Omega) \times H^1_0(\Omega).
\end{multline*}
Note that here again $\mu_\varepsilon^m$ is the orthogonal projection of $v_\varepsilon^m$ onto $H^1_0(\Omega)$.

\begin{proposition} \label{prop_asatisfiesA2}
The bilinear form $a$ satisfies  \autoref{assumption_sym_of_a_on_solutions}.
\end{proposition}

\begin{proof}
We focus on the identity \( a(x_\varepsilon^m, x_\varepsilon^k) = a(x_\varepsilon^k, x_\varepsilon^m) \), as the identity \( a(x_\varepsilon^m, x_s) = a(x_s, x_\varepsilon^m) \) follows by an analogous argument.

Let \( m \) and \( k \) be two integers. A direct computation shows that the required identity holds if
\[
\int_\Omega \nabla v_\varepsilon^m \cdot \nabla \mu_\varepsilon^k \, \mathrm{d}x = \int_\Omega \nabla v_\varepsilon^k \cdot \nabla \mu_\varepsilon^m \, \mathrm{d}x.
\]
This equality is satisfied because \(\mu_\varepsilon^m\) and \(\mu_\varepsilon^k\) are the orthogonal projections of \(v_\varepsilon^m\) and \(v_\varepsilon^k\) onto \(H^1_0(\Omega)\), respectively. Consequently, we obtain
\[
\int_\Omega \nabla v_\varepsilon^m \cdot \nabla \mu_\varepsilon^k \, \mathrm{d}x = \int_\Omega \nabla \mu_\varepsilon^m \cdot \nabla \mu_\varepsilon^k \, \mathrm{d}x = \int_\Omega \nabla \mu_\varepsilon^m \cdot \nabla v_\varepsilon^k \, \mathrm{d}x.
\]
\end{proof}

We are now in position to study the iterated quasi-reversibility sequence adapted to the data assimilitation problem \eqref{eq_dataassimpb}. We defer the discussion on discretization and perturbation aspects to the next section and, for now, focus on the continuous framework without perturbation. In other words, we consider the \emph{exact data} $(f,\mathbf F,g) \in L^2(\Omega) \times L^2(\Omega)^d\times L^2(\omega)$, for which we assume that problem \eqref{eq_dataassimpb} admits an unique solution $u_s \in H^1(\Omega)$. Then, defining $\lambda^0$ by \eqref{eq_def_lambda} with $f^0 = f$ and $\mathbf{F}^0 = \mathbf{F}$, and setting $\lambda_s \in H^1_0(\Omega)$ as the orthogonal projection of $u$ into $H^1_0(\Omega)$, we have by  \autoref{prop_asatisfiesA1} that $x_s = (u_s,\lambda_s)$ is the unique solution of \eqref{prob_abstract_a} -- and $\lambda^0 = \lambda_s$.

In this context, for \(\varepsilon > 0\), the iterated quasi-reversibility sequence \(X_\varepsilon^M = (u_\varepsilon^m, \lambda_\varepsilon^m) \in H^1(\Omega) \times H^1_0(\Omega)\) is initialized by setting \(X_\varepsilon^0 = (0, 0)\). Then, for any \(M \in \mathbb{N}\), the pair \(X_\varepsilon^{M+1} = (u_\varepsilon^{M+1}, \lambda_\varepsilon^{M+1})\) satisfies \eqref{eq_abstract_itQR_continuous}, which reads as:
\begin{multline*}
\int_\Omega \nabla \lambda_\varepsilon^{M+1} \cdot \nabla v \, dx + \int_\omega u_\varepsilon^{M+1} v \, dx + \int_\Omega \nabla \lambda_\varepsilon^{M+1} \cdot \nabla \mu \, dx - \int_\Omega \nabla u_\varepsilon^{M+1} \cdot \nabla \mu \, dx + \varepsilon \int_\Omega \nabla u_\varepsilon^{M+1} \cdot \nabla v \, dx \\
= \int_\Omega \nabla \lambda_s \cdot \nabla v \, dx + \int_\omega g v \, dx + \varepsilon \int_\Omega \nabla u_\varepsilon^M \cdot \nabla v \, dx, \quad \forall (v, \mu) \in H^1(\Omega) \times H^1_0(\Omega).
\end{multline*}
Alternatively, this can be expressed in the equivalent system:
\begin{equation} \label{eq_itQRDA_core}
	\left\{
		\begin{aligned}
			&\begin{multlined}[0.9\linewidth]
				\int_\Omega \nabla \lambda_\varepsilon^{M+1} \cdot \nabla v \, dx + \int_\omega u_\varepsilon^{M+1} v \, dx + \varepsilon \int_\Omega \nabla u_\varepsilon^{M+1} \cdot \nabla v \, dx \\
				= \int_\Omega \nabla \lambda_s \cdot \nabla v \, dx + \int_\omega g v \, dx + \varepsilon \int_\Omega \nabla u_\varepsilon^M \cdot \nabla v \, dx, \quad \forall v \in H^1(\Omega),
			\end{multlined} \\
			&\int_\Omega \nabla \lambda_\varepsilon^{M+1} \cdot \nabla \mu \, dx - \int_\Omega \nabla u_\varepsilon^{M+1} \cdot \nabla \mu \, dx = 0, \quad \forall \mu \in H^1_0(\Omega).
		\end{aligned}
	\right.
\end{equation}
It is worth noting that in the right-hand side of \eqref{eq_itQRDA_core}$_1$, \(g\) may be replaced by \(u_s\). Additionally, we observe that \(\lambda_\varepsilon^M\) is once again the orthogonal projection of \(u_\varepsilon^M\) onto \(H^1_0(\Omega)\).

\begin{theorem} \label{Thm_itQRDA_aprioriestimates}
For all \(\varepsilon > 0\), the sequence \((u_\varepsilon^M, \lambda_\varepsilon^M)\) is well-defined. Moreover, it converges strongly to \((u_s, \lambda_s)\) in \(H^1(\Omega) \times H^1_0(\Omega)\) and satisfies the following estimates:
\begin{equation} \label{eq_itQRDA_aprioriestimates}
\begin{cases}
\displaystyle \Vert \nabla u_\varepsilon^M \Vert_{L^2(\Omega)} \leq \Vert \nabla u_{\varepsilon}^{M+1} \Vert_{L^2(\Omega)} \leq \Vert \nabla u_s \Vert_{L^2(\Omega)}, \\ \\
\displaystyle \Vert \nabla (u_\varepsilon^{M+1} - u_s) \Vert_{L^2(\Omega)} \leq \Vert \nabla (u_\varepsilon^M - u_s ) \Vert_{L^2(\Omega)}, \\ \\
\displaystyle \sqrt{\Vert \lambda_\varepsilon^M \Vert_{L^2(\Omega)}^2 + \Vert u_\varepsilon^M \Vert_{L^2(\omega)}^2} \leq \sqrt{\Vert \nabla \lambda_s \Vert_{L^2(\Omega)}^2 + \Vert u_s \Vert_{L^2(\omega)}^2}, \\ \\
\displaystyle  \sqrt{\Vert \nabla (\lambda_\varepsilon^{M+1} - \lambda_s) \Vert_{L^2(\Omega)}^2 + \Vert u_\varepsilon^{M+1} - u_s \Vert_{L^2(\omega)}^2}
\leq \sqrt{\Vert \nabla (\lambda_\varepsilon^{M} - \lambda_s) \Vert_{L^2(\Omega)}^2 + \Vert u_\varepsilon^{M} - u_s \Vert_{L^2(\omega)}^2} \\
 \hfill \displaystyle \leq  \sqrt{\frac{\varepsilon}{M}} \Vert \nabla u_s \Vert_{L^2(\Omega)}.
\end{cases}
\end{equation}
\end{theorem}

\begin{remark}
The right-hand side of \eqref{eq_itQRDA_aprioriestimates}$_3$ can be bounded from above by $c \Vert u_s \Vert_{H^1(\Omega)}$, for some constant $c \in \R$. Indeed, as $\lambda_s$ is the orthogonal projection of $u_s$ onto $H^1_0(\Omega)$, we have $\Vert \nabla \lambda_s \Vert_{L^2(\Omega)} \leq \Vert \nabla u_s \Vert_{L^2(\Omega)}$, which immediately implies the result.
\end{remark}

\begin{proof}[Proof of  \autoref{Thm_itQRDA_aprioriestimates}]
We simply note that, as for all $(u,\lambda) \in H^1(\Omega) \times H^1_0(\Omega)$, $$\vert L(u,\lambda) \vert \leq \sqrt{\Vert \lambda_s\Vert_{L^2(\Omega)}^2 + \Vert u_s \Vert_{L^2(\omega)}^2 } \Vert u,\lambda\Vert_a,$$ we have $c_a \leq \sqrt{\Vert \lambda_s\Vert_{L^2(\Omega)}^2 + \Vert u_s \Vert_{L^2(\omega)}^2 }$.
As $a$, $b$ and $L$ satisfy all the required assumptions,
the result is a direct application of  \autoref{prop_aprioriestimates_bis},  \autoref{prop_aprioriestimates_abstractsetting} and  \autoref{prop_residualestimates_abstractsetting}.
\end{proof}

\subsection{Regularity estimates} \label{subsection_regularity_estimates}

The main objective of this section is to examine how increased regularity in the exact solution $u_s$ translates into improved regularity for the iterated quasi-reversibility solution
$(u_\varepsilon^M,\lambda_\varepsilon^M)$. These regularity results will be used in the next section to obtain convergence estimates of the finite elements discretization.
Note that the primary challenge here is not merely establishing regularity results, but deriving regularity estimates with an explicit dependence on the parameters $\varepsilon$ and $M$.

\medskip

For the reader convenience, the elliptic regularity results underlying our analysis are summarized in the  \autoref{appendix_miscres},
specifically in  \autoref{thm_regularityresults_Diri} and  \autoref{thm_regularityresults_Neum}.

\medskip

From now on, we assume that \( u_s \in H^1(\Omega, \Delta) \). We will demonstrate that, under this additional regularity assumption on \( u_s \), both \( u_\varepsilon^M \) and \( \lambda_\varepsilon^M \) enjoy enhanced regularity properties.

We recall that the norm in \( H^1(\Omega, \Delta) \) is given by:
\[
\| u_s \|_{H^1(\Omega, \Delta)} = \sqrt{\| u_s \|_{H^1(\Omega)}^2 + \| \Delta u_s \|_{L^2(\Omega)}^2}.
\]
We begin by noting that \(\lambda_s\) exhibits higher regularity.

\begin{proposition}
There exists $\sigma \in \left] \frac{1}{2},1\right]$ and $c$, depending only on $\Omega$, such that $\lambda_s$ belongs to $H^1(\Omega,\Delta) \cap H^{1+\sigma}(\Omega)$ and
$$
\Delta \lambda_s = \Delta u_s, \quad \Vert \lambda_s \Vert_{H^{1+\sigma}(\Omega)} \leq c \Vert \Delta u_s \Vert_{L^2(\Omega)}.
$$
\end{proposition}

\begin{proof}
By definition of $\lambda_s$ and regularity of $u_s$, we have
$$
\int_\Omega \nabla \lambda_s \cdot \nabla \mu dx = \int_\Omega \nabla u_s \cdot\nabla \mu dx = - \int_\Omega \Delta u_s \ \mu dx, \quad \forall \mu \in H^1_0(\Omega).
$$
Choosing $\mu$ in $C^\infty_c(\Omega)$ implies $\Delta \lambda_s = \Delta u_s$. The result is then a direct application of  \autoref{thm_regularityresults_Diri}.
\end{proof}

We can now turn our attention to $\lambda_\varepsilon^M$ and $u_\varepsilon^M$.

\begin{proposition} \label{prop_estimate_diffLapla}
Set $\varepsilon_0>0$. There exists a constant $c>0$ such that for any $\varepsilon \in (0,\varepsilon_0]$, for any $M\in \N$,
$\lambda_\varepsilon^M$ and $u_\varepsilon^M$ belong to $H^1(\Omega,\Delta)$, and
$$
\Vert \Delta (\lambda_\varepsilon^M - \lambda_s) \Vert_{L^2(\Omega)} + \Vert \Delta (u_\varepsilon^M - u_s ) \Vert_{L^2(\Omega)}
\leq c \sqrt{\frac{\varepsilon}{M}} \Vert u_s \Vert_{H^1(\Omega,\Delta)}.
$$
\end{proposition}

\begin{proof}
We know that $u_\varepsilon^M$ and $\lambda_\varepsilon^M$ belong to $H^1(\Omega)$, so we only need to focus on the estimate on the Laplacian. Furthermore, as by definition of $\lambda_s$ we have
$$
\int_\Omega \nabla \lambda_s \cdot \nabla \phi\, dx = \int_\Omega \nabla u_s \cdot \nabla \phi\, dx, \quad  \forall \phi \in C^\infty_c(\Omega),
$$
we have $\Delta \lambda_s = \Delta u_s \in L^2(\Omega)$. Similarly, from \eqref{eq_itQRDA_core}$_2$ we obtain $\Delta \lambda_\varepsilon^M = \Delta u_\varepsilon^M$. For our purpose, it is therefore sufficient to prove that $\Delta \lambda_\varepsilon^M$ belongs to $L^2(\Omega)$, and to derive an $L^2$-estimate for $\Delta (\lambda_\varepsilon^M - \lambda_s)$.

Choosing $v = \mu \in C^\infty_c$ in \eqref{eq_itQRDA_core} easily leads to
$$
\Delta (\lambda_\varepsilon^{M+1}  - \lambda_s ) = \frac{1}{1+\varepsilon} (u_\varepsilon^{M+1} - u_s ) 1_\omega + \frac{\varepsilon}{1+\varepsilon} \Delta (\lambda_\varepsilon^M - \lambda_s).
$$
  \autoref{lemma_wm_explicit} with $q = \varepsilon (1+\varepsilon)^{-1}$, $w_M = \Delta (\lambda_\varepsilon^M - \lambda_s) $, $\eta_M =  \frac{1}{1+\varepsilon} (u_\varepsilon^{M+1} - u_s ) 1_\omega  $ and $\kappa = 0$ gives
 $$
 \Delta (\lambda_\varepsilon^M - \lambda_s)  = - \frac{\varepsilon^M}{(1+\varepsilon)^M} \Delta \lambda_s + \frac{1}{1+\varepsilon}\sum_{k=0}^{M-1} \frac{\varepsilon^k}{(1+\varepsilon)^k}
 (u_\varepsilon^{M-k} - u_s) 1_\omega.
 $$
 In particular, $\Delta \lambda_\varepsilon^M$ belongs to $L^2(\Omega)$.
Combining \eqref{eq_itQRDA_aprioriestimates}$_4$ and  \autoref{lemma_bound_sum}, we have, for some constant $c>0$ depending only on $\varepsilon_0$,
\begin{multline} \label{eq_proof_estimateDelta}
\sum_{k=0}^{M-1} \frac{\varepsilon^k}{(1+\varepsilon)^k} \Vert u_\varepsilon^{M-k} - u_s \Vert_{L^2(\omega)}
\leq \sqrt{\varepsilon}  \Vert \nabla u_s \Vert_{L^2(\Omega)} \sum_{k=0}^{M-1} \frac{\varepsilon^k}{(1+\varepsilon)^k} \frac{1}{\sqrt{M-k}} \\
\leq c (1+\varepsilon)  \sqrt{\frac{\varepsilon}{M}} \Vert \nabla u_s \Vert_{L^2(\Omega)}.
\end{multline}
The desired estimate follows directly.
\end{proof}

\begin{remark}
In the final step of the proof, it could be tempting to simply bound $\varepsilon (1+\varepsilon)^{-1}$ by 1. However, since
$$
\sum_{k=0}^{M-1} \Vert u_\varepsilon^{M-k} - u_s \Vert_{L^2(\omega)}\lesssim \sqrt \varepsilon \Vert \nabla u_s \Vert_{L^2(\Omega)} \sum_{k=1}^M \frac{1}{\sqrt k} \lesssim
\sqrt{\varepsilon\, M} \Vert \nabla u_s\Vert_{L^2(\Omega)},
$$
this naive bound would degrade the estimate in terms of powers of $M$.

We recover the correct asymptotic by introducing the upper bound $\varepsilon_0$ and using  \autoref{lemma_bound_sum}. 
It is worth noting that the constant $c$ appearing in \eqref{eq_proof_estimateDelta} diverges as $\varepsilon_0 \rightarrow \infty$.
\end{remark}

\begin{proposition} \label{prop_ellipticestimate_lambdaepsm}
Set $\varepsilon_0 >0$.
There exists $\sigma \in \left] \frac{1}{2},1 \right]$, depending only on $\Omega$,  and $c>0$, depending only on $\varepsilon_0$ and $\Omega$, such that
for any $\varepsilon \in (0,\varepsilon_0]$ and any $M\in \N_*$, $\lambda_\varepsilon^M$ belongs to $H^{1+\sigma}(\Omega)$ and
$$
\Vert \lambda_\varepsilon^M \Vert_{H^{1+\sigma}(\Omega)} \leq c \left( \left( 1 - \frac{\varepsilon^M}{(1+\varepsilon)^M} \right)\Vert \Delta u_s \Vert_{L^2(\Omega)} + \sqrt{\frac{\varepsilon}{M}} \Vert \nabla u_s \Vert_{L^2(\Omega)} \right).
$$
\end{proposition}

\begin{proof}
	Let $\mu \in H^1_0(\Omega)$. Using simultaneously \eqref{eq_itQRDA_core}$_1$ and \eqref{eq_itQRDA_core}$_2$  gives
	$$
	(1+\varepsilon) \int_\Omega \nabla \lambda_{\varepsilon}^{M+1} \cdot \nabla \mu dx =
	 \int_\Omega \nabla \lambda_s \cdot \nabla \mu dx + \int_\omega (u_s- u_{\varepsilon}^{M+1}) \mu dx + \varepsilon \int_\Omega \nabla \lambda_{\varepsilon}^{M} \cdot \nabla \mu dx,
	$$
	which can be rewritten $w_{M+1} = q w_M + \kappa + \eta_M$ with $q = \varepsilon (1+\varepsilon)^{-1} \not\in\left\lbrace 0,1 \right\rbrace$ and, for all $M\in \N$,
	$$
	w_M = \int_\Omega \nabla \lambda_{\varepsilon}^{M} \cdot \nabla \mu dx, \ \kappa = \frac{1}{1+\varepsilon}  \int_\Omega \nabla \lambda_s \cdot \nabla \mu dx, \
	\eta_M = \frac{1}{1+\varepsilon} \int_\omega (u_s- u_{\varepsilon}^{M+1}) \mu dx .
	$$
	Note that by definition $w_0 = 0$.
	 \autoref{lemma_wm_explicit} gives
	\begin{multline*}
	\int_\Omega \nabla \lambda_{\varepsilon}^{M} \cdot \nabla \mu dx = \frac{1-q^M}{1-q} \frac{1}{1+\varepsilon}  \int_\Omega \nabla \lambda_s \cdot \nabla \mu dx
	+ \sum_{k=0}^{M-1} \frac{q^k}{1+\varepsilon}  \int_\omega (u_s- u_{\varepsilon}^{M-p}) \mu dx\\
	= - \left(1 - \frac{\varepsilon^M}{(1+\varepsilon)^M}\right) \int_\Omega \Delta \lambda \, \mu dx + \frac{1}{1+\varepsilon} \sum_{k=0}^{M-1} q^k \int_\omega (u_s- u_{\varepsilon}^{M-k}) \mu dx.
	\end{multline*}
	As this holds for any $\mu$ in $H^1_0(\Omega)$,  \autoref{thm_regularityresults_Diri} implies that there exists $\sigma\in \left]\frac{1}{2}, 1 \right]$ such that  $\lambda_{\varepsilon}^{M}$ belongs to $H^{1+\sigma}(\Omega)$ and
	$$
	\Vert \lambda_{\varepsilon}^{M} \Vert_{H^{1+\sigma}(\Omega)} \lesssim  \left(1 - \frac{\varepsilon^M}{(1+\varepsilon)^M}\right) \Vert \Delta \lambda_s \Vert_{L^2(\Omega)}
	+ \frac{1}{1+\varepsilon} \sum_{k=0}^{M-1} q^k \Vert u_{\varepsilon}^{M-k} - u_s \Vert_{L^2(\omega)}.
	$$
	Estimate \eqref{eq_itQRDA_aprioriestimates}$_4$ combined with  \autoref{lemma_bound_sum} gives, for some constant $c$ depending only on $\varepsilon_0$,
	$$
	 \sum_{k=0}^{M-1} q^k \Vert u_{\varepsilon}^{M-k} - u_s \Vert_{L^2(\omega)} \leq \sqrt \varepsilon \Vert \nabla u_s \Vert_{L^2(\Omega)} \sum_{k=0}^{M-1}
	 	\frac{q^k}{\sqrt{M-k}} \leq c (1+\varepsilon) \sqrt{\frac{\varepsilon}{M}} \Vert \nabla u_s \Vert_{L^2(\Omega)}.
	$$
	The result follows.
\end{proof}

\begin{proposition} \label{prop_lambdaemminuslambdaH2}
Set $\varepsilon_0 >0$.
There exists $\sigma \in \left] \frac{1}{2},1 \right]$, depending only on $\Omega$,  and $c>0$, depending only on $\varepsilon_0$ and $\Omega$, such that
for any $\varepsilon \in (0,\varepsilon_0]$ and any $M\in \N_*$,
$$
\Vert \lambda_{\varepsilon}^{M} - \lambda_s \Vert_{H^{1+\sigma}(\Omega)}
\leq c \sqrt{\frac{\varepsilon}{M}}  \left(  \Vert \Delta u_s \Vert_{L^2(\Omega)} + \Vert \nabla u_s\Vert_{L^2(\Omega)} )\right..
$$
\end{proposition}

\begin{proof}
We proceed along the same lines as before. For $\mu$ in $H^1_0(\Omega)$,  we set
$$w_M = \int_\Omega \nabla (\lambda_{\varepsilon}^{M} - \lambda_s) \cdot \nabla \mu dx,$$ and
$$\eta_M = (1+\varepsilon)^{-1} \int_\omega (g-u_{\varepsilon}^{M+1}) \mu dx.$$
We see that
$w_{M+1} = q w_M + \eta_M$ with $q = \varepsilon (1+\varepsilon)^{-1}$, which immediately gives ( \autoref{lemma_wm_explicit})
$$
w_M = q^M \, w_0 + \sum_{k=0}^{M-1} q^{k} \eta_{M-1-k},
$$
or equivalently
$$
\int_\Omega \nabla (\lambda_\varepsilon^M - \lambda_s) \cdot \nabla \mu dx = q^M \int_\Omega \Delta \lambda_s\, \mu dx
+ \frac{1}{1+\varepsilon} \sum_{k=0}^{M-1} q^k \int_\omega (u_s - u_{\varepsilon}^{M-k}) \mu dx.
$$
As this holds for any $\mu$ in $H^1_0(\Omega)$, we deduce from  \autoref{thm_regularityresults_Diri} the existence of $\sigma \in \left]\frac{1}{2},1\right]$  such that
$$
\Vert  \lambda_{\varepsilon}^{M} - \lambda_s \Vert_{H^{1+\sigma}(\Omega)}
\lesssim q^M \Vert \Delta \lambda_s \Vert_{L^2(\Omega)} +  \sum_{k=0}^{M-1} q^{k} \Vert u_s- u_{\varepsilon}^{M-k} \Vert_{L^2(\omega)}.
$$
As previously, we have, for some constant $c>0$ depending only on $\varepsilon_0$,
$$
\sum_{k=0}^{M-1} q^{k} \Vert u_s- u_{\varepsilon}^{M-k} \Vert_{L^2(\omega)} \leq c \sqrt{\frac{\varepsilon}{M}} \Vert \nabla u_s \Vert_{L^2(\Omega)}.
$$
As furthermore for all $\varepsilon >0$ and for all $M\geq 1$,
$$
q^M \sqrt{\frac{M}{\varepsilon}} = \frac{\varepsilon^{M-\frac{1}{2}}}{(1+\varepsilon)^M} \sqrt M \leq \frac{1}{2},
$$
the result follows.
\end{proof}

\begin{proposition} \label{prop_ellipticestimate_uepsm}
Set $\varepsilon_0 >0$.
There exists $\sigma \in \left] \frac{1}{2},1 \right]$, depending only on $\Omega$,  and $c>0$, depending only on $\varepsilon_0$ and $\Omega$, such that  for any $\varepsilon \in (0,\varepsilon_0]$ and any $M\in \N_*$,
	 $u_{\varepsilon}^{M}$ belongs to $H^{1+\sigma}(\Omega)$, and
	$$
	\Vert u_{\varepsilon}^{M} \Vert_{H^{1+\sigma}(\Omega)} \leq c \sqrt{\frac{M}{\varepsilon}} \Vert u_s \Vert_{H^1(\Omega,\Delta)}.
	$$
\end{proposition}

\begin{proof}
	Let $v$ in $H^1(\Omega)$ be fixed. Equation \eqref{eq_itQRDA_core}$_1$ gives, for all $k \in \N$,
	$$
	\int_\Omega \nabla (\lambda_{\varepsilon}^{k+1}-\lambda_s + \varepsilon (u_{\varepsilon}^{k+1} - u_{\varepsilon}^{k})) \cdot \nabla v dx =   \int_\omega (u_s-u_{\varepsilon}^{k+1}) v dx.
	$$
	Summing these equalities from $k=0$ to $k= M-1$ yields
	$$
	\int_\Omega \nabla \left( \sum_{k=0}^{M-1} (\lambda_{\varepsilon}^{k+1} - \lambda_s)  + \varepsilon u_{\varepsilon}^{M}\right) \cdot \nabla v dx = \int_\omega \left( \sum_{k=0}^{M-1} \left(u_s-u_{\varepsilon}^{k+1}\right) \right) v dx.
	$$
	Then  \autoref{thm_regularityresults_Neum} gives that $\displaystyle \sum_{k=0}^{M-1} (\lambda_{\varepsilon}^{k+1} - \lambda_s)  + \varepsilon u_{\varepsilon}^{M}$ belongs to $H^{1+\sigma}(\Omega)$ and verifies
	$$
	\left\Vert \sum_{k=0}^{M-1} (\lambda_{\varepsilon}^{k+1} - \lambda_s)  + \varepsilon u_{\varepsilon}^{M} \right\Vert_{H^{1+\sigma}(\Omega)}
	\lesssim \left( \left\Vert \sum_{k=0}^{M-1} (\lambda_{\varepsilon}^{k+1} - \lambda_s)  + \varepsilon u_{\varepsilon}^{M} \right\Vert_{L^2(\Omega)}
	+ \left\Vert \sum_{k=0}^{M-1}\left(u_s-u_{\varepsilon}^{k+1}\right) \right\Vert_{L^2(\omega)}\right).
	$$
	As we know by the previous Proposition that $\lambda_{\varepsilon}^{k+1} - \lambda_s$  belongs to $H^{1+\sigma}(\Omega)$ for all $k$, we deduce that $u_{\varepsilon}^{M}$ belongs
	to $H^{1+\sigma}(\Omega)$ and verifies
	$$
	\varepsilon \Vert u_{\varepsilon}^{M} \Vert_{H^{1+\sigma}(\Omega)} \lesssim \left(\varepsilon \Vert u_{\varepsilon}^{M} \Vert_{L^2(\Omega)} + \sum_{k=0}^{M-1}
	 \left( \Vert \lambda_{\varepsilon}^{k+1} - \lambda_s \Vert_{H^{1+\sigma}(\Omega)} + \Vert u_s - u_{\varepsilon}^{k+1} \Vert_{L^2(\omega)} \right)\right).
	$$
	On one hand, as  for all $w\in H^1(\Omega)$, $\Vert w\Vert_{H^1(\Omega)} \lesssim\left(\Vert \nabla w\Vert_{L^2(\Omega)} + \Vert w \Vert_{L^2(\omega)}\right)$, we deduce from \eqref{eq_itQRDA_aprioriestimates}$_1$ and \eqref{eq_itQRDA_aprioriestimates}$_3$ that $\Vert u_\varepsilon^M \Vert_{L^2(\Omega)} \lesssim \Vert u_s \Vert_{H^1(\Omega)}$.
	On the other hand, the previous Proposition combined with \eqref{eq_itQRDA_aprioriestimates}$_4$ gives
	$$
	\sum_{k=0}^{m-1}
	\left( \Vert \lambda_{\varepsilon}^{k+1} - \lambda_s \Vert_{H^{1+\sigma}(\Omega)} + \Vert u_s - u_{\varepsilon}^{k+1} \Vert_{L^2(\omega)} \right) \lesssim
	\sqrt{\varepsilon} \Vert u_s \Vert_{H^1(\Omega,\Delta)} \sum_{k=0}^{M-1} \frac{1}{\sqrt{k+1}} \lesssim  \sqrt{\varepsilon\, M} \Vert u_s \Vert_{H^1(\Omega,\Delta)}.
	$$
	The result follows.
\end{proof}

\begin{remark} \label{rk_increased_bdn_regularity}
In this study, the domain $\Omega$ is polygonal or polyhedric. As a consequence, the domain may present entring corners, which limit the regularity of the solutions.
Nevertheless, we have the following properties.
\begin{enumerate}
\item For any $\tilde{\Omega}$ open set compactly supported in $\Omega$, $u_\varepsilon^M$ and $\lambda_\varepsilon^M$ belong to $H^2(\tilde \Omega)$. This follows directly from cutoff arguments and analogous reasoning as above.
\item If $\Omega$ is convex, than one may take $\sigma= 1$ in all above estimates. The same applies if $\Omega$ is a domain with a smooth boundary.
\end{enumerate}
\end{remark}

\begin{remark} \label{rk_limitation_regularity}
If $\Omega$ is convex and $u_s$ belongs to $H^1(\Omega,\Delta)$, we find that both $u_\varepsilon^M$ and $\lambda_\varepsilon^M$ lie in $H^2(\Omega)$. Unfortunately, assuming higher regularity for $u_s$ does not yield improved regularity for $u_\varepsilon^M$ and $\lambda_\varepsilon^M$. To illustrate this phenomenon, consider the first term of the sequence -- the solution to the standard quasi-reversibility method $(u_\varepsilon,\lambda_\varepsilon) = (u_\varepsilon^1,\lambda_\varepsilon^1)$. A direct computation yields
$$
\Delta \lambda_\varepsilon = \frac{1}{1+\varepsilon} \Delta \lambda_s + \frac{1}{1+\varepsilon} (u_\varepsilon - u_s) 1_\omega \text{ in } \Omega.
$$
If $u_s$ smooth, then $\lambda_s$ is also smooth, hence the first term on the right-hand side is smooth. However, the presence of the characteristic function of $\omega$ in the second term limits its regularity to at best $L^2$. Consequently, $\lambda_\varepsilon$ is at best $H^2$, as is $u_\varepsilon$.
\end{remark}

\subsection{Convergence estimate for the iterated quasi-reversibility method} \label{subsec_logconv}

As already stated in  \autoref{Thm_itQRDA_aprioriestimates}, 
we  know that, for any choice of the positive regularization parameter $\varepsilon$, $(u_\varepsilon^M,\lambda_\varepsilon^M)$ converges to $(u_s,\lambda_s)$ strongly in $H^1(\Omega)\times H^1_0(\Omega)$ as $M$ goes to infinity. However, the abstract study does not provide any rate of convergence.

\medskip

Still, the derived regularity estimates yield convergence estimate for sufficiently smooth exact solution, \text{i.e.} $u_s \in H^1(\Omega,\Delta)$. More precisely,  \autoref{prop_lambdaemminuslambdaH2}, which gives
$$
\Vert \lambda_{\varepsilon}^{M} - \lambda_s \Vert_{H^{1+s}(\Omega)}
\leq c \sqrt{\frac{\varepsilon}{M}}  \left(  \Vert \Delta u_s \Vert_{L^2(\Omega)} + \Vert \nabla u_s\Vert_{L^2(\Omega)} \right),
$$
is a quantified convergence estimate, in a  stronger norm  than the $H^1$-norm.

\medskip

Interestingly, the estimates obtained in  \autoref{Thm_itQRDA_aprioriestimates} and  \autoref{prop_estimate_diffLapla}, combined with known conditional stability estimates for the data assimilation problem, yield a convergence rate for $u_\varepsilon^M - u_s$, albeit in a weaker $L^2$-norm, as seen in the following result.

\begin{proposition} \label{prop_logestimateQR}
There exists $\mu \in (0,1)$ and, for all $\varepsilon_0>0$, a constant $c>0$, such that for all $u_s$ in $H^1(\Omega,\Delta)$, all $\varepsilon \in (0,\varepsilon_0)$ and all $M \in \N_*$,
$$
\Vert u_\varepsilon^M - u_s \Vert_{L^2(\Omega)} \leq c \ln \left( 1+ \frac{M}{\varepsilon} \right)^{-\mu} \Vert u_s \Vert_{H^1(\Omega,\Delta)}.
$$
\end{proposition}

\begin{proof}
Combining the estimates \eqref{eq_itQRDA_aprioriestimates}$_1$ and \eqref{eq_itQRDA_aprioriestimates}$_4$, and Proposition \eqref{prop_estimate_diffLapla} gives a constant $c$, depending only on $\varepsilon_0$ such that for all $\varepsilon \in (0,\varepsilon_0]$ and $M\in\N_*$,
$$
\sqrt{\frac{M}{\varepsilon}} \Vert \Delta (u_\varepsilon^M - u_s ) \Vert_{L^2(\Omega)} + \sqrt{\frac{M}{\varepsilon}} \Vert u_\varepsilon^M - u_s \Vert_{L^2(\omega)}
+ \Vert u_\varepsilon^M - u_s \Vert_{H^1(\Omega)} \leq c \Vert u_s \Vert_{H^1(\Omega,\Delta)}.
$$
We proceed by considering distinct cases: first, we suppose that $\varepsilon$ and $M$ are such that $\sqrt \frac{\varepsilon}{M}\leq \frac{1}{2}$, which implies $\frac{2 \sqrt{\frac{\varepsilon}{M}}}{1+ \sqrt{\frac{\varepsilon}{M}}} < 1$. Applying  \autoref{thm_GlobalStability_Alessandrini} with
$$
\eta_g = \eta_f = c \sqrt{\frac{\varepsilon}{M}} \Vert u_s \Vert_{H^1(\Omega,\Delta)}, \quad \mathcal{K} = c \Vert u_s \Vert_{H^1(\Omega,\Delta)},
$$
gives a constant $ c >0$ and $\mu \in (0,1)$, independent of $u_s$, $\varepsilon$ and $M$, such that
$$
\Vert u_\varepsilon^M - u_s \Vert_{L^2(\Omega)} \leq  c \left( 1+ \sqrt{\frac{\varepsilon}{M}} \right) \ln \left(\frac{1+\sqrt{\frac{\varepsilon}{M}}}{2\sqrt{\frac{\varepsilon}{M}}} \right)^{-\mu} \Vert u_s \Vert_{H^1(\Omega,\Delta)}.
$$
Since the function $f:t \in \left(0,\frac{1}{2}\right] \mapsto \ln \left(1+ \frac{1}{t^2}\right) \ln \left( \frac{1+t}{2t}\right)^{-1}$ is positive, continuous, and satisfies
$f(t) \xrightarrow[]{t\rightarrow 0^+}2$, it is bounded. Consequently, we obtain a new constant $c>0$ such that, for all $\varepsilon$ and $M$ verifying $\sqrt \frac{\varepsilon}{M}\leq \frac{1}{2}$,
$$
\Vert u_\varepsilon^M - u_s \Vert_{L^2(\Omega)} \leq  c \ln \left( 1+ \frac{M}{\varepsilon} \right)^{-\mu} \Vert u_s \Vert_{H^1(\Omega,\Delta)}.
$$
Now, in the case $\sqrt \frac{\varepsilon}{M} \geq \frac{1}{2}$, we simply note that $\ln\left( 1+\frac{M}{\varepsilon}\right) \leq \ln (5)$, which immediatly gives the result as $\Vert u_\varepsilon^M - u_s \Vert_{L^2(\Omega)} \lesssim\Vert u_s \Vert_{H^1(\Omega,\Delta)}$.
\end{proof}

\section{Finite element discretization} \label{sec_FEMdisc}

We now turn our attention to convergence estimates for the finite element discretization of the iterated quasi-reversibility variational formulation.
 We recall that in our study, 
 for $h>0$, we denote by $\mathcal{T}_h$ an admissible, shape regular and uniform family of partitions of the polytope $\Omega$
into triangular elements ($d=2$) or tetrahedral elements ($d=3$), such that for every $\tau_h$ in $\mathcal{T}_h$, the diameter of $\tau_h$ is bounded by $h$.
As stated in the Introduction, we denote
$$
V_h = \left\lbrace w \in H^1(\Omega), \ \forall \tau_h \in \mathcal{T}_h, \ w_{|\tau_h} \in \mathbb{R}_1[X] \right\rbrace, \quad V_{h,0} = V_h \cap H^1_0(\Omega).
$$
$V_h$ is the standard \emph{conforming $P_1$-Lagrange finite element}, while $V_{h,0}$ is the subspace of functions of $V_h$ that vanish on $\partial \Omega$.

%

\subsection{\texorpdfstring{$H^1$}{H1}-convergence} \label{subsection_FEM_H1conv}

Since $V_h$ (resp. $V_{h,0}$) is a finite-dimensional subspace of $H^1(\Omega)$ (resp. $H^1_0(\Omega)$),
$\mathcal X_h = V_h \times V_{h,0}$ is a closed subspace of $\mathcal X = H^1(\Omega) \times H^1_0(\Omega)$, and the Galerkin approximation of the iterated quasi-reversibility sequence within the space $\mathcal X_h$ can be analysed using the abstract results developed in  \autoref{sec_abstractQR}.

\medskip We first consider the case of unperturbed data. In other words, for $\varepsilon > 0$ fixed, we approximate the quasi-reversibility sequence $X_\varepsilon^M = (u_\varepsilon^M, \lambda_\varepsilon^M)$ defined by \eqref{eq_itQRDA_core} by the sequence
$
X_{\varepsilon,h}^M = (u_{\varepsilon,h}^M, \lambda_{\varepsilon,h}^M) \in \mathcal{X}_h
$
defined by $X_{\varepsilon,h}^0 = (0,0)$ and, for all $M\in \N$,

\begin{equation} \label{eq_itQRDAapprox_core}
\begin{cases}
\displaystyle \int_\Omega \nabla \lambda_{\varepsilon,h}^{M+1} \cdot \nabla v_h dx + \int_\omega u_{\varepsilon,h}^{M+1} v_h dx + \varepsilon \int_\Omega \nabla u_{\varepsilon,h}^{M+1} \cdot \nabla v_h dx
=  \int_\Omega \nabla \lambda_s \cdot \nabla v_h dx + \int_\omega g\, v_h dx
\\ \hfill \displaystyle  + \varepsilon \int_\Omega \nabla u_{\varepsilon,h}^M \cdot \nabla v_h \, dx, \quad \forall v_h \in V_h,
\\
\\ \displaystyle \int_\Omega \nabla \lambda_{\varepsilon,h}^{M+1} \cdot \nabla \mu_h dx - \int_\Omega \nabla u_{\varepsilon,h}^{M+1} \cdot \nabla \mu_h dx = 0, \hfill \quad \forall \mu_h \in V_{h,0}.
\end{cases}
\end{equation}
or equivalently $a(X_{\varepsilon,h}^{M+1},x_h) + \varepsilon b(X_{\varepsilon,h}^{M+1},x_h) = L(x_h) + \varepsilon b(X_{\varepsilon,h}^M,x_h)$, for all $x_h$ in $\mathcal{X}_h$.

\begin{theorem} \label{thm_errorestimatediscretization}
For all $\varepsilon >0$ and $h>0$, the sequence $X_{\varepsilon,h}^M = (u_{\varepsilon,h}^M, \lambda_{\varepsilon,h}^M)$ is well-defined. Furthermore, there exists $\sigma \in \left(\frac{1}{2},1\right]$, depending only on $\Omega$, and, for all $\varepsilon_0>0$, $c>0$ depending only on $\Omega$ and $\varepsilon_0$, such that if $u_s$ belongs to $H^1(\Omega,\Delta)$, for all $h>0$, $\varepsilon \in (0,\varepsilon_0]$ and $M\in \N_*$,
$$
\Vert u_{\varepsilon,h}^M - u_\varepsilon^M \Vert_{H^1(\Omega)} + \Vert \lambda_{\varepsilon,h}^M - \lambda_\varepsilon^M \Vert_{H^1(\Omega)}
\leq c h^\sigma \left(\frac{M}{\varepsilon}\right)^{\frac{3}{2}} \Vert u_s \Vert_{H^1(\Omega,\Delta)}.
$$
\end{theorem}

\begin{proof}
From  \autoref{prop_abstract_Galerkin_estimation} follows the existence of a constant $c$, depending only on $\varepsilon_0$, such that
\begin{multline*}
\Vert u_{\varepsilon,h}^M - u_\varepsilon^M \Vert_{H^1(\Omega)} + \Vert \lambda_{\varepsilon,h}^M - \lambda_\varepsilon^M \Vert_{H^1(\Omega)}
\\ \leq \frac{c}{\varepsilon} \sum_{p=0}^{M} \inf_{(u_h,\lambda_h) \in V_h \times V_{h,0}} \left(\Vert u_\varepsilon^p - u_h \Vert_{H^1(\Omega)} + \Vert \lambda_\varepsilon^p - \lambda_h \Vert_{H^1(\Omega)}\right)
\\ = \frac{c}{\varepsilon} \sum_{p=0}^{M} \left( \inf_{u_h \in V_h} \Vert u_\varepsilon^p - u_h \Vert_{H^1(\Omega)} +  \inf_{\lambda_h \in V_{h,0}}  \Vert \lambda_\varepsilon^p - \lambda_h \Vert_{H^1(\Omega)}\right)
\end{multline*}
 \autoref{thm_globalapprox_FE} combined with  \autoref{prop_ellipticestimate_lambdaepsm} and  \autoref{prop_ellipticestimate_uepsm} gives $\sigma \in \left( \frac{1}{2},1\right]$
depending only on $\Omega$, and $c>0$, depending only on $\varepsilon_0$ and $\Omega$, such that
$$
 \inf_{u_h \in V_h} \Vert u_\varepsilon^p - u_h \Vert_{H^1(\Omega)} +  \inf_{\lambda_h \in V_{h,0}}  \Vert \lambda_\varepsilon^p - \lambda_h \Vert_{H^1(\Omega)}
 \leq c \sqrt{\frac{p}{\varepsilon}} h^\sigma \Vert u_s \Vert_{H^1(\Omega,\Delta)}.
$$
As $\displaystyle \sum_{p=0}^{M} \sqrt{p} \leq \frac{2}{3} (M+1)^{\frac{3}{2}}$, the result follows.
\end{proof}

We now turn to our final estimate, which addresses the iterated quasi-reversibility method in the presence of noisy data. We recall that, for $\delta>0$, the perturbed data is given by
$(f^\delta,\mathbf{F}^\delta, g^\delta) \in L^2(\Omega) \times L^2(\Omega)^d \times L^2(\omega)$ verifying
$$
\Vert f^\delta - f \Vert_{L^2(\Omega)} + \Vert \mathbf{F}^\delta - \mathbf{F} \Vert_{L^2(\Omega)^d} + \Vert g^\delta - g \Vert_{L^2(\omega)} \leq \delta.
$$
For $\varepsilon>0$, we consider the sequence $X_{\varepsilon,h}^{M,\delta} = (u_{\varepsilon,h}^{M,\delta},\lambda_{\varepsilon,h}^{M,\delta}) \in \mathcal{X}_h$ defined by \eqref{eq_itQRDAapproxnoisy_core}.

\begin{proposition} \label{prop_errorestimate_noise}
There exists a constant $c>0$ such that, for any $\delta >0$, for any $h>0$, for any $\varepsilon>0$ and for any $M\in \N$,
$$
\Vert u_{\varepsilon,h}^{M,\delta} - u_{\varepsilon,h}^M \Vert_{H^1(\Omega)} + \Vert \lambda_{\varepsilon,h}^{M,\delta} - \lambda_{\varepsilon,h}^M \Vert_{H^1(\Omega)}
\leq c \sqrt{\frac{M}{\varepsilon}} \delta.
$$
\end{proposition}

\begin{proof}
The definitions of $\lambda^\delta$ and $\lambda_s$ imply
$$
\int_\Omega \nabla (\lambda^\delta- \lambda_s) \cdot \nabla  \mu dx = \int_\Omega (f^\delta - f ) \mu dx - \int_\Omega (\mathbf{F}^\delta - \mathbf{F}) \cdot \nabla \mu dx,
\quad \forall \mu \in H^1_0(\Omega).
$$
It follows  that $\Vert \lambda^\delta - \lambda_s \Vert_{H^1(\Omega)} \lesssim \delta$. Now, setting
$$
L : (v,\mu) \in \mathcal{X} \mapsto \int_\Omega \nabla \lambda_s \cdot \nabla v dx + \int_\omega g v dx, \quad
 L^\delta : (v,\mu) \in \mathcal{X} \mapsto \int_\Omega \nabla \lambda^\delta \cdot \nabla v dx + \int_\omega g^\delta v dx,
$$
we have
$$
\Vert L^\delta - L \Vert \lesssim   \Vert \lambda^\delta - \lambda_s \Vert_{H^1(\Omega)} + \Vert g^\delta - g \Vert_{L^2(\omega)}\lesssim  \delta.
$$
The result follows from  \autoref{prop_abstracterrorestimatenoise}.
\end{proof}

\subsection{\texorpdfstring{$L^2$}{L2}-convergence}

By adapting the classical Aubin--Nitsche approach to our problem, we derive an improved convergence estimate in the $L^2(\Omega)$-norm, given by:


\begin{proposition}
	\label{prop_superconvergence_l2}
	Let $\sigma \in \left(\frac{1}{2},1\right]$ be given by \autoref{thm_errorestimatediscretization}, and $\eps_0 > 0$. For any $\eps \in (0, \eps_0]$ and any $M \in \N$, there exists a constant $c > 0$ such that for any $h > 0$,
	\[\norm{u_{\eps,h}^M - u_\eps^M}[L^2(\Omega)] + \norm{\lambda_{\eps,h}^M - \lambda_\eps^M}[L^2(\Omega)] \leq c h^{2\sigma} \norm{u_s}[H^1(\Omega,\Delta)].\]
\end{proposition}

We emphasize that, unlike all previous results, our analysis does not track the dependence of the constant $c$ on $\varepsilon$ and $M$. This choice is motivated by two considerations. First, the computations are already cumbersome, and explicitly tracking these dependencies would further complicate them. Second, the dependence we obtain through our proof scheme is likely suboptimal, particularly when compared to the dependencies derived in earlier works (see \autoref{rk_dependancy_constant_L2} below).

In order to simplify the notation, let us introduce the norm
\[\norm{(u, \lambda)}[L]^2 \coloneqq \norm{u}[L^2(\Omega)]^2 + \norm{\lambda}[L^2(\Omega)]^2, \quad \forall (u, \lambda) \in \mcX.\]
We denote $\scalprod{}{}[L]$ the associated scalar product.

Recall that, for any $\eps > 0$, the form $Q_\eps$ satisfies the assumptions of Lax-Milgram's theorem. Hence, for all $M \in \N$, we can define $W_\eps^M = (w_\eps^M, \ell_\eps^M) \in \mcX$ such that
\begin{equation}
	\label{pb_aubin}
	Q_\eps(x, W_\eps^M) = \scalprod{X_\eps^M - X_{\eps,h}^M}{x}[L], \quad \forall x \in \mcX,
\end{equation}
alongside its discretized counterpart $W_{\eps,h}^M = (w_{\eps,h}^M, \ell_{\eps,h}^M) \in \mcX_h$ such that
\begin{equation}
	\label{pb_aubin_h}
	Q_\eps(x_h, W_{\eps,h}^M) = \scalprod{X_\eps^M - X_{\eps,h}^M}{x_h}[L], \quad \forall x_h \in \mcX_h.
\end{equation}
Note that \eqref{pb_aubin} reads
	\begin{equation}
		\label{pb_aubin_system}
		\left\{
			\begin{aligned}
				\int_\omega w_\eps^M v - \int_\Omega \nabla \ell_\eps^M \cdot \nabla v + \eps \int_\Omega \nabla w_\eps^M \cdot \nabla v &= \int_\Omega (u_\eps^M - u_{\eps,h}^M) v, && \forall v \in H^1(\Omega), \\
				\int_\Omega \nabla w_\eps^M \cdot \nabla \mu + \int_\Omega \nabla \ell_\eps^M \cdot \nabla \mu &= \int_\Omega (\lambda_\eps^M - \lambda_{\eps,h}^M) \mu, && \forall \mu \in H^1_0(\Omega).
			\end{aligned}
		\right.
	\end{equation}

\begin{lemma}
	\label{lemma_err_Wh}
	Let $\sigma \in \interval[open left, scaled]{\frac{1}{2}}{1}$ be given by \autoref{thm_errorestimatediscretization}. For any $\eps > 0$, there exists a constant $c > 0$ such that, for any $M \in \N$ and $h > 0$:
	\[\norm{W_\eps^M - W_{\eps,h}^M}[Q_\eps] \leq c h^\sigma \norm{X_\eps^M - X_{\eps,h}^M}[L].\]
\end{lemma}

\begin{proof}
%
	By combining the variational formulation \eqref{pb_aubin_system} with the regularity results \autoref{thm_regularityresults_Diri} and \autoref{thm_regularityresults_Neum}, we deduce --  using the same reasoning as in \autoref{subsection_regularity_estimates} -- that both $w_\eps^M$ and $\ell_\eps^M$ belong to $H^{1 + \sigma}(\Omega)$. Furthermore, there exists a constant $c$ depending only on $\varepsilon$ such that
	$$
	\Vert w_{\eps}^M\Vert_{H^{1+\sigma}(\Omega)} + \Vert \ell_{\eps}^M \Vert_{H^{1+\sigma}(\Omega)} \leq c \norm{X_\eps^M - X_{\eps,h}^M}[L].
	$$
The desired result then follows directly from Galerkin orthogonality and \autoref{thm_globalapprox_FE}.
\end{proof}

\begin{lemma}
	\label{lemma_integrationparts_W}
	For any $\eps > 0$ and $M \in \N$, $w_\eps^M \in H^1(\Omega, \Delta)$, $\ell_\eps^M \in H^1(\Omega, \Delta)$ and:
	\[\forall v \in H^1(\Omega), \quad \eps \int_\Omega (\Delta w_\eps^M) v + \eps \int_\Omega \nabla w_\eps^M \cdot \nabla v = \int_\Omega (\Delta \ell_\eps^M) v + \int_\Omega \nabla \ell_\eps^M \cdot \nabla v.\]
	Moreover, there exist two constants $c_w > 0$ and $c_\ell > 0$, depending only on $\eps$, such that:
	\[\norm{\Delta w_\eps^M}[L^2(\Omega)] \leq c_w \norm{X_\eps^M - X_{\eps,h}^M}[L] \quad \text{and} \quad \norm{\Delta \ell_\eps^M}[L^2(\Omega)] \leq c_\ell \norm{X_\eps^M - X_{\eps,h}^M}[L].\]
\end{lemma}

\begin{proof}
	Considering $v = \mu = \varphi \in \mcC^\infty_0(\Omega)$ in \eqref{pb_aubin_system} yields
	\begin{equation}
		\label{pb_aubin_system_strong}
		\left\{
			\begin{aligned}
				\chi_\omega w_\eps^M + \Delta \ell_\eps^M - \eps \Delta w_\eps^M &= u_\eps^M - u_{\eps,h}^M, \\
				-\Delta w_\eps^M - \Delta \ell_\eps^M &= \lambda_\eps^M - \lambda_{\eps,h}^M,
			\end{aligned}
		\right.
	\end{equation}
	which is equivalent to
	$$
	\begin{cases}
	-(1 + \eps) \Delta w_\eps^M = u_\eps^M - u_{\eps,h}^M + \lambda_\eps^M - \lambda_{\eps,h}^M - \chi_\omega w_\eps^M, \\
	\Delta \ell_\eps^M = -\Delta w_\eps^M - (\lambda_\eps^M - \lambda_{\eps,h}^M).
	\end{cases}
	$$
	We directly conclude that both $\Delta w_\eps^M$ and $\Delta \ell_\eps^M$ belong to $L^2(\Omega)$, and we derive the corresponding $L^2$-estimates using \autoref{lemma_err_Wh}.
	Finally, multiplying  \eqref{pb_aubin_system_strong}$_1$ by  $v \in H^1(\Omega)$, and integrating over $\Omega$, leads to
	\[\int_\Omega (\Delta \ell_\eps^M) v - \eps \int_\Omega (\Delta w_\eps^M) v = \int_\Omega (u_\eps^M - u_{\eps,h}^M) v - \int_\omega w_\eps^M v.\]
	The claimed identity follows from \eqref{pb_aubin_system}$_1$.
\end{proof}

\begin{proof}[Proof of \autoref{prop_superconvergence_l2}]
Writing \eqref{pb_aubin} at step $M + 1$ with $x = X_\eps^{M + 1} - X_{\eps,h}^{M + 1}$ leads to
\begin{align*}
	\norm{X_\eps^{M + 1} - X_{\eps,h}^{M + 1}}[L]^2 &= Q_\eps(X_\eps^{M + 1} - X_{\eps,h}^{M + 1}, W_\eps^{M + 1}) \\
	&= Q_\eps(X_\eps^{M + 1} - X_{\eps,h}^{M + 1}, W_\eps^{M + 1} - W_{\eps,h}^{M + 1}) + Q_\eps(X_\eps^{M + 1} - X_{\eps,h}^{M + 1}, W_{\eps,h}^{M + 1}).
\end{align*}
As $W_{\eps,h}^{M + 1} \in \mcX_h$, the definition of $X_\eps^{M + 1}$ and $X_{\eps,h}^{M + 1}$ directly yields
\begin{align*}
	Q_\eps(X_\eps^{M + 1} - X_{\eps,h}^{M + 1}, W_{\eps,h}^{M + 1}) &= \eps b(X_\eps^M - X_{\eps,h}^M, W_{\eps,h}^{M + 1}) \\
	&= \eps b(X_\eps^M - X_{\eps,h}^M, W_{\eps,h}^{M + 1} - W_\eps^{M + 1}) + \eps b(X_\eps^M - X_{\eps,h}^M, W_\eps^{M + 1}).
\end{align*}
We use \autoref{lemma_integrationparts_W} to rewrite:
\begin{multline*}
	\eps b(X_\eps^M - X_{\eps,h}^M, W_\eps^{M + 1}) = \eps \int_\Omega \nabla (u_\eps^M - u_{\eps,h}^M) \cdot \nabla w_\eps^{M + 1} \\
	= \int_\Omega (\Delta \ell_\eps^{M + 1}) (u_\eps^M - u_{\eps,h}^M) - \eps \int_\Omega (\Delta w_\eps^{M + 1}) (u_\eps^M - u_{\eps,h}^M) + \int_\Omega \nabla \ell_\eps^{M + 1} \cdot \nabla (u_\eps^M - u_{\eps,h}^M).
\end{multline*}
The last integral can be rewritten as
\[\int_\Omega \nabla \ell_\eps^{M + 1} \cdot \nabla (u_\eps^M - u_{\eps,h}^M) = \int_\Omega \nabla (\ell_\eps^{M + 1} - \ell_{\eps,h}^{M + 1}) \cdot \nabla (u_\eps^M - u_{\eps,h}^M) + \int_\Omega \nabla \ell_{\eps,h}^{M + 1} \cdot \nabla (u_\eps^M - u_{\eps,h}^M).\]
As $\ell_{\eps,h}^{M + 1} \in V_{h,0}$, the definition of $u_\eps^M$ and $\lambda_\eps^M$ allows to write
\[\int_\Omega \nabla \ell_{\eps,h}^{M + 1} \cdot \nabla (u_\eps^M - u_{\eps,h}^M) = \int_\Omega \nabla \ell_{\eps,h}^{M + 1} \cdot \nabla (\lambda_\eps^M - \lambda_{\eps,h}^M).\]
Once again, we split this integral:
\[\int_\Omega \nabla \ell_{\eps,h}^{M + 1} \cdot \nabla (\lambda_\eps^M - \lambda_{\eps,h}^M) = \int_\Omega \nabla (\ell_{\eps,h}^{M + 1} - \ell_\eps^{M + 1}) \cdot \nabla (\lambda_\eps^M - \lambda_{\eps,h}^M) + \int_\Omega \nabla \ell_\eps^{M + 1} \cdot \nabla (\lambda_\eps^M - \lambda_{\eps,h}^M).\]
As $\lambda_\eps^M$ and $\lambda_{\eps,h}^M$ belong in $H^1_0(\Omega)$, and $\ell_\eps^{M + 1} \in H^1(\Omega, \Delta)$, an integration by parts gives
\[\int_\Omega \nabla \ell_\eps^{M + 1} \cdot \nabla (\lambda_\eps^M - \lambda_{\eps,h}^M) = -\int_\Omega (\Delta \ell_\eps^{M + 1}) (\lambda_\eps^M - \lambda_{\eps,h}^M).\]
We finally obtained
\begin{multline*}
	\norm{X_\eps^{M + 1} - X_{\eps,h}^{M + 1}}[L]^2 = Q_\eps(X_\eps^{M + 1} - X_{\eps,h}^{M + 1}, W_\eps^{M + 1} - W_{\eps,h}^{M + 1}) + \eps b(X_\eps^M - X_{\eps,h}^M, W_{\eps,h}^{M + 1} - W_\eps^{M + 1}) \\
	+ \int_\Omega \nabla (\ell_\eps^{M + 1} - \ell_{\eps,h}^{M + 1}) \cdot \nabla ((u_\eps^M - u_{\eps,h}^M) - (\lambda_\eps^M - \lambda_{\eps,h}^M)) \\
	+ \int_\Omega (\Delta \ell_\eps^{M + 1}) ((u_\eps^M - u_{\eps,h}^M) - (\lambda_\eps^M - \lambda_{\eps,h}^M)) - \eps \int_\Omega (\Delta w_\eps^{M + 1}) (u_\eps^M - u_{\eps,h}^M).
\end{multline*}

As in the proof of \autoref{prop_abstract_Galerkin_estimation}, one can easily show that:
\[\forall (x, \tilde{x}) \in \mcX^2, \quad Q_\eps(x, \tilde{x}) \lesssim \left(1 + \max\left(1, \frac{1}{\eps}\right)\right) \norm{x}[Q_\eps] \norm{\tilde{x}}[Q_\eps].\]
Hence:
\[Q_\eps(X_\eps^{M + 1} - X_{\eps,h}^{M + 1}, W_\eps^{M + 1} - W_{\eps,h}^{M + 1}) \lesssim \left(1 + \max\left(1, \frac{1}{\eps}\right)\right) \norm{X_\eps^{M + 1} - X_{\eps,h}^{M + 1}}[Q_\eps] \norm{W_\eps^{M + 1} - W_{\eps,h}^{M + 1}}[Q_\eps].\]
Using \autoref{lemma_err_Wh} and \autoref{thm_errorestimatediscretization}, we then find a constant $c > 0$ such that
\[Q_\eps(X_\eps^{M + 1} - X_{\eps,h}^{M + 1}, W_\eps^{M + 1} - W_{\eps,h}^{M + 1}) \leq c h^{2\sigma} \norm{u_s}[H^1(\Omega, \Delta)] \norm{X_\eps^{M + 1} - X_{\eps,h}^{M + 1}}[L].\]
Similarly, as
\[\forall (x, \tilde{x}) \in \mcX^2, \quad \eps b(x, \tilde{x}) \leq \norm{x}[Q_\eps] \norm{\tilde{x}}[Q_\eps],\]
there exists another constant $c > 0$ such that
\[\eps b(X_\eps^M - X_{\eps,h}^M, W_{\eps,h}^{M + 1} - W_\eps^{M + 1}) \leq c h^{2\sigma} \norm{u_s}[H^1(\Omega, \Delta)] \norm{X_\eps^{M + 1} - X_{\eps,h}^{M + 1}}[L].\]
Now, by Cauchy-Schwarz inequality:
\begin{align*}
	\int_\Omega \nabla (\ell_\eps^{M + 1} - \ell_{\eps,h}^{M + 1}) &\cdot \nabla ((u_\eps^M - u_{\eps,h}^M) - (\lambda_\eps^M - \lambda_{\eps,h}^M)) \\
	&\leq \norm{\nabla (\ell_\eps^{M + 1} - \ell_{\eps,h}^{M + 1})}[L^2(\Omega)] (\norm{\nabla (u_\eps^M - u_{\eps,h}^M)}[L^2(\Omega)] + \norm{\nabla (\lambda_\eps^M - \lambda_{\eps,h}^M)}[L^2(\Omega)]) \\
	&\leq \norm{W_\eps^{M + 1} - W_{\eps,h}^{M + 1}}[a] (\norm{X_\eps^{M + 1} - X_{\eps,h}^{M + 1}}[b] + \norm{X_\eps^{M + 1} - X_{\eps,h}^{M + 1}}[a]) \\
	&\leq \max(1, \eps) \norm{W_\eps^{M + 1} - W_{\eps,h}^{M + 1}}[Q_\eps] \norm{X_\eps^{M + 1} - X_{\eps,h}^{M + 1}}[Q_\eps].
\end{align*}
Once again, we find a constant $c > 0$ such that:
\[\int_\Omega \nabla (\ell_\eps^{M + 1} - \ell_{\eps,h}^{M + 1}) \cdot \nabla ((u_\eps^M - u_{\eps,h}^M) - (\lambda_\eps^M - \lambda_{\eps,h}^M)) \leq c h^{2\sigma} \norm{u_s}[H^1(\Omega, \Delta)] \norm{X_\eps^{M + 1} - X_{\eps,h}^{M + 1}}[L].\]
The last two integrals can be controlled in a similar way. We first write
\[\int_\Omega (\Delta w_\eps^{M + 1}) (u_\eps^M - u_{\eps,h}^M) \leq \norm{\Delta w_\eps^{M + 1}}[L^2(\Omega)] \norm{u_\eps^M - u_{\eps,h}^M}[L^2(\Omega)] \leq \norm{\Delta w_\eps^{M + 1}}[L^2(\Omega)] \norm{X_\eps^M - X_{\eps,h}^M}[L].\]
Applying \autoref{lemma_integrationparts_W} and \autoref{thm_errorestimatediscretization}, we then find a constant $c > 0$ such that
\[\int_\Omega (\Delta w_\eps^{M + 1}) (u_\eps^M - u_{\eps,h}^M) \leq c \norm{X_\eps^M - X_{\eps,h}^M}[L] \norm{X_\eps^{M + 1} - X_{\eps,h}^{M + 1}}[L].\]
In the end, we proved the existence of a constant $c > 0$ such that:
\[\norm{X_\eps^{M + 1} - X_{\eps,h}^{M + 1}}[L]^2 \leq c (h^{2\sigma} \norm{u_s}[H^1(\Omega, \Delta)] + \norm{X_\eps^M - X_{\eps,h}^M}[L]) \norm{X_\eps^{M + 1} - X_{\eps,h}^{M + 1}}[L].\]
Then, \autoref{prop_superconvergence_l2} follows by induction.
\end{proof}

\begin{remark} \label{rk_dependancy_constant_L2}
	In contrast to \autoref{subsection_FEM_H1conv}, \autoref{prop_superconvergence_l2} does not provide an explicit dependence of the constant on the parameters $\varepsilon$ and $M$.

Although it is possible to track this dependence in the steps outlined above, doing so results in a constant that grows exponentially with respect to $M$. Additionally, the resulting constant is not homogeneous with respect to the ratio $\frac{M}{\varepsilon}$. This suggests that, while the chosen approach suffices to establish the $L^2$ superconvergence of the finite element method, it does not yield a reasonable dependence on $M$ and $\varepsilon$.
\end{remark}

\section{Numerical results and validation} \label{sec_NumIll}

This section presents numerical results illustrating the different estimates that have been obtained in the study. The meshes used in the simulations have all been generated with Gmsh~\cite{geuzaineremacle09}. The complete code is written in Python and is available in an online repository~\cite{heleine26}, where scripts to reproduce the presented results can be found. The finite element matrices are assembled with the library scikit-fem~\cite{gustafssonmcbain20}. We use SymPy~\cite{sympy} to generate source terms from analytic solutions, and we use the libraries NumPy~\cite{numpy} and SciPy~\cite{scipy} for the numerical computations. In particular, we use SciPy's sparse LU decomposition to factorize the matrix of the system, significantly decreasing the computation time of the iterates. Finally, we use the library Matplotlib~\cite{hunter07}  to produce the figures.


\subsection{Finite element rate of convergence}

As the primary objective of the study is to derive convergence estimates for the finite element approximation of the iterated quasi-reversibility sequence,
we first examine the finite element convergence on its own, \textit{i.e.} \autoref{thm_errorestimatediscretization}. More precisely, for a fixed choice of the parameters
$\varepsilon$ and $M$, \autoref{thm_errorestimatediscretization} gives
\[\norm{u_{\eps,h}^M - u_\eps^M}[H^1(\Omega)] + \norm{\lambda_{\eps,h}^M - \lambda_\eps^M}[H^1(\Omega)] \lesssim h^\sigma,\]
where $\sigma \in \left(\frac{1}{2}, 1\right]$ depends only on $\Omega$.
Furthermore, \autoref{prop_superconvergence_l2} gives, for the same $\sigma$,
\[\norm{u_{\eps,h}^M - u_\eps^M}[L^2(\Omega)] + \norm{\lambda_{\eps,h}^M - \lambda_\eps^M}[L^2(\Omega)] \lesssim h^{2\sigma}.\]
This is exactly the estimates we seek to illustrate numerically. However, we face a difficulty: the continuous iterated quasi-reversibility solution $(u_\varepsilon^M,\lambda_\varepsilon^M)$ is not accessible for any choice of the parameters. To address this, we propose to compare the numerical sequence to a reference sequence, computed on a finer mesh of size $h_0$. Indeed, a simple application of the triangular inequality gives
\[\norm{u_{\eps,h}^M - u_{\eps,h_0}^M}[H^1(\Omega)] + \norm{\lambda_{\eps,h}^M - \lambda_{\eps,h_0}^M}[H^1(\Omega)] \lesssim h^\sigma + h_0^\sigma.\]
By varying  $h$, we can observe the convergence rate governed by $h^\sigma$, the term $h_0^\sigma$ being negligible if $h_0$ is chosen small enough.



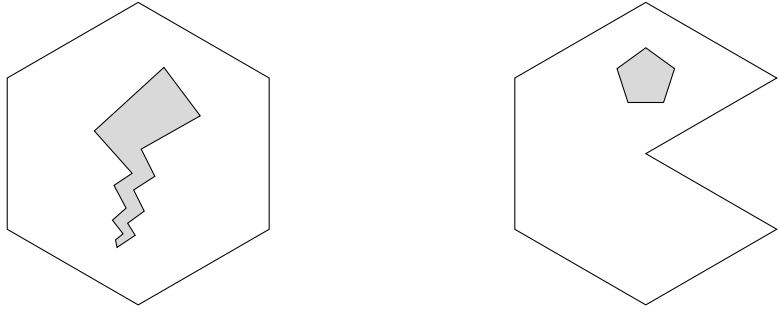
\begin{figure}[hbt]
	\centering

	\begin{tikzpicture}[scale = 2]
		\draw
			(0.8660254037844387, 0.5)
			-- (0,1)
			-- (-0.8660254037844385, 0.5)
			-- (-0.8660254037844388, -0.5)
			-- (0, -1)
			-- (0.8660254037844384, -0.5)
			-- cycle
		;

		\draw[fill = gray, fill opacity = 0.3]
			(0.17, 0.57)
			-- (-0.29, 0.15)
			-- (-0.04, -0.13)
			-- (-0.16, -0.21)
			-- (-0.08, -0.36)
			-- (-0.17, -0.44)
			-- (-0.1, -0.53)
			-- (-0.15, -0.57)
			-- (-0.14, -0.62)
			-- (-0.02, -0.54)
			-- (-0.07, -0.46)
			-- (0.04, -0.38)
			-- (-0.03, -0.24)
			-- (0.11, -0.15)
			-- (0.02, 0.03)
			-- (0.41, 0.25)
			-- cycle
		;
	\end{tikzpicture}
	\hspace{3cm}
	\begin{tikzpicture}[scale = 2]
		\draw
			(0, 0)
			-- (0.8660254037844387, 0.5)
			-- (0, 1)
			-- (-0.8660254037844385, 0.5)
			-- (-0.8660254037844388, -0.5)
			-- (0, -1)
			-- (0.8660254037844384, -0.5)
			--cycle
		;

		\draw[fill = gray, fill opacity = 0.3]
			(0, 0.7)
			-- (-0.1902113032590307, 0.5618033988749895)
			-- (-0.11755705045849466, 0.3381966011250105)
			-- (0.11755705045849459, 0.3381966011250105)
			-- (0.19021130325903074, 0.5618033988749894)
			-- cycle
		;
	\end{tikzpicture}

	\caption{Geometries of interest in the numerical tests. Left: a convex geometry. Right: a non-convex geometry. The data subdomain $\omega$ is highlighted in gray in both cases.}
	\label{fig_geometries}
\end{figure}

We use two different geometries. In the first one, $\Omega$ is a regular hexagon, hence a convex polygon. In the second case, we replaced one side by an entring corner, ending up with a non-convex domain. The two geometries are presented in \autoref{fig_geometries}. Note that, in the non-convex geometry, the corner is of radius $\frac{\pi}{3}$ and placed at the origin.

The parameters of the meshes generated from these geometries are listed in \autoref{tab_meshsizes_convex} and \autoref{tab_meshsizes_nonconvex}. The reference mesh is named $M_0$.

\begin{table}[hbt]
	\centering

	\begin{tabular}{|c|c|c|c|c|c|}
		\hline
		Convex & $M_0$ & $M_1$ & $M_2$ & $M_3$ & $M_4$ \\
		\hline
		$h$ & \scnum[2]{0.011059937083712128} & \scnum[2]{0.050748081328974536} & \scnum[2]{0.0806721125562766} & \scnum[2]{0.10612492638395539} & \scnum[2]{0.12529964086141668} \\
		\hline
		Triangles & \num{96972} & \num{4386} & \num{2032} & \num{1242} & \num{808} \\
		\hline
		Vertices & \num{48862} & \num{2269} & \num{1068} & \num{661} & \num{435} \\
		\hline
	\end{tabular}

	\caption{Parameters of the meshes used to test the convergence rates in the case of the convex geometry.}
	\label{tab_meshsizes_convex}
\end{table}

\begin{table}[hbt]
	\centering

	\begin{tabular}{|c|c|c|c|c|c|}
		\hline
		Non-convex & $M_0$ & $M_1$ & $M_2$ & $M_3$ & $M_4$ \\
		\hline
		$h$ & \scnum[2]{0.010778604132880848} & \scnum[2]{0.05420068051002905} & \scnum[2]{0.07698941877111512} & \scnum[2]{0.09831517559219528} & \scnum[2]{0.12670382525626672} \\
		\hline
		Triangles & \num{79361} & \num{3255} & \num{1493} & \num{863} & \num{592} \\
		\hline
		Vertices & \num{40119} & \num{1716} & \num{807} & \num{478} & \num{332} \\
		\hline
	\end{tabular}

	\caption{Parameters of the meshes used to test the convergence rates in the case of the non-convex geometry.}
	\label{tab_meshsizes_nonconvex}
\end{table}

\subsubsection{Smooth exact solution}

In this first test, the data $f$ and $g$ correspond to a smooth exact solution given by
\[u_s \colon (x, y) \mapsto x^2 + y.\]
Its projection onto $H^1_0(\Omega)$, noted $\lambda_s$, is obtained numerically by computing $\Delta u_s$ and then solving
\[\left\{\begin{aligned}
	\Delta \lambda_s &= \Delta u_s, && \text{in } \Omega, \\
	\lambda_s &= 0, && \text{on } \partial\Omega.
\end{aligned}\right.\]
Note that the exact solution is used solely to generate the data and does not appear in the convergence estimates we illustrate.

\begin{figure}[hbt]
	\centering

	\includegraphics[width=0.45\textwidth]{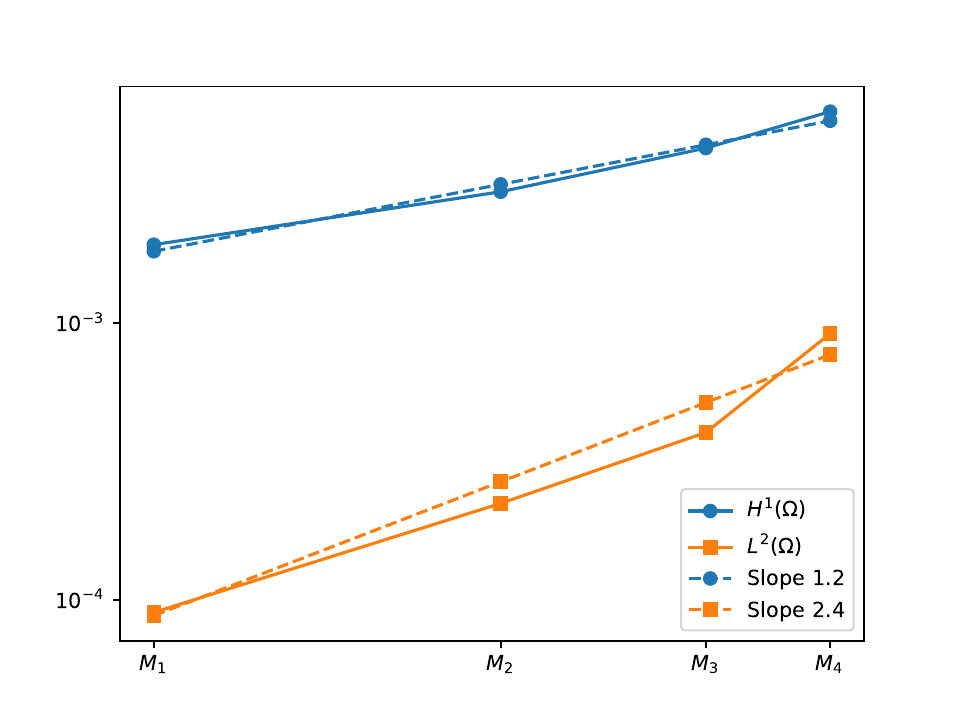}
	\hfill
	\includegraphics[width=0.45\textwidth]{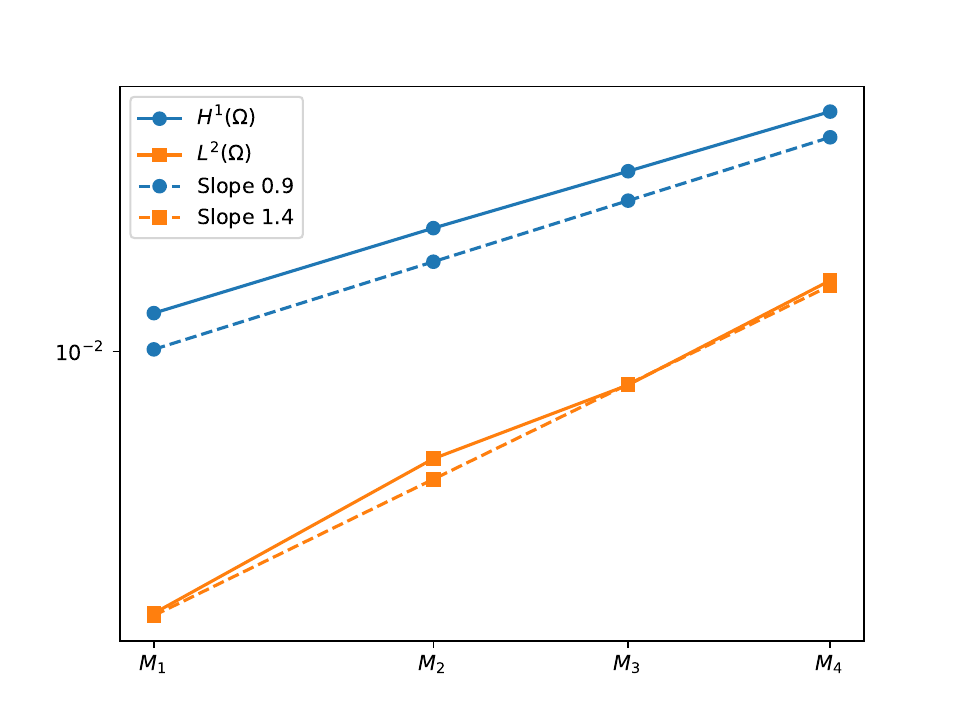}

	\caption{Illustration of the convergence rate of the error with respect to $h$. Left: convex geometry. Right: non-convex geometry.}
	\label{fig_convergence_h}
\end{figure}

The iterates of the quasi-reversibility are first computed on the reference mesh of size $h_0$, with $\eps = 10^{-2}$, until we reach $M = 100$. For each mesh of size $h$, the iterates are computed with the same values of the parameters $\varepsilon$ and $M$. Then, we compute the relative error in $H^1(\Omega)$-norm
\[\frac{\norm{u_{\eps,h}^M - u_{\eps,h_0}^M}[H^1(\Omega)] + \norm{\lambda_{\eps,h}^M - \lambda_{\eps,h_0}^M}[H^1_0(\Omega)]}{\norm{u_{\eps,h_0}^M}[H^1(\Omega)] + \norm{\lambda_{\eps,h_0}^M}[H^1_0(\Omega)]}.\]
In addition, we compute the relative error in $L^2(\Omega)$-norm using the same formula, but replacing all the norms with $\norm{}[L^2(\Omega)]$. In \autoref{fig_convergence_h}, we show the evolution of these errors with respect to $h$.

We observe the expected orders of convergence in the convex geometry, as $\sigma = 1$ in that case. Furthermore, the $L^2(\Omega)$-error decreases faster than the  $H^1(\Omega)$-error, as predicted by \autoref{prop_superconvergence_l2}. Same behaviors are observed in the second geometry, but with a slower rate of convergence, which is expected as  $\sigma < 1$ for non-convex geometry. 

\subsubsection{Singular exact solution}

\begin{figure}[hbt]
	\centering

	\includegraphics[width=0.45\textwidth]{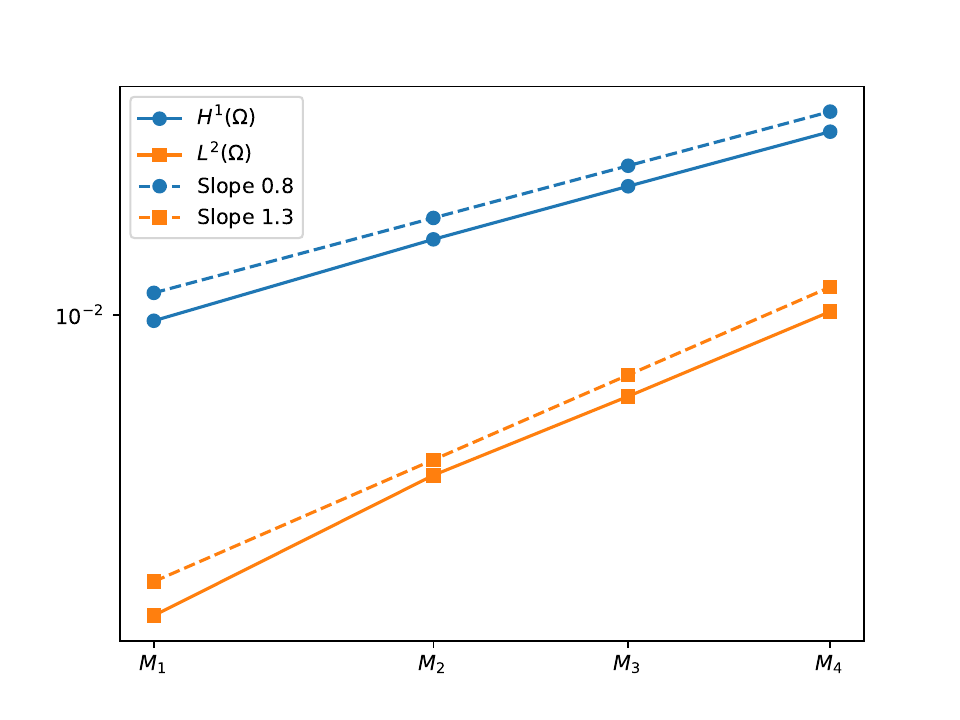}

	\caption{Evolution of the finite element error w.r.t. $h$, with a singular exact solution, in the non-convex geometry.}
	\label{fig_convergence_h_singular}
\end{figure}

We now do a second test in the non-convex geometry, with a different exact solution given by (in polar coordinates)
\[u_s \colon (r, \theta) \mapsto \sin\left(\frac{3}{5} \left(\theta - \frac{\pi}{6}\right)\right) r^{3/5}.\]
This function is harmonic and belongs to $H^1(\Omega)$, but fails to be $H^2(\Omega)$ due to a singularity in its second derivatives at the origin. 
The finite element convergence is presented in \autoref{fig_convergence_h_singular}.

We observe numerical convergence rates similar to those obtained in the previous section.
This is expected, since our analysis shows that despite the now-singular exact solution,
the iterates $u_\varepsilon^M$ still belong to $H^2(\Omega)$. Consequently, the singularity
in $u_s$ introduces no additional loss.

\subsection{Effect of the parameters \texorpdfstring{$\varepsilon$}{ε} and \texorpdfstring{$M$}{M} on the finite element and reconstruction errors}

\begin{figure}[hbt]
	\centering

	\includegraphics[width=0.45\textwidth]{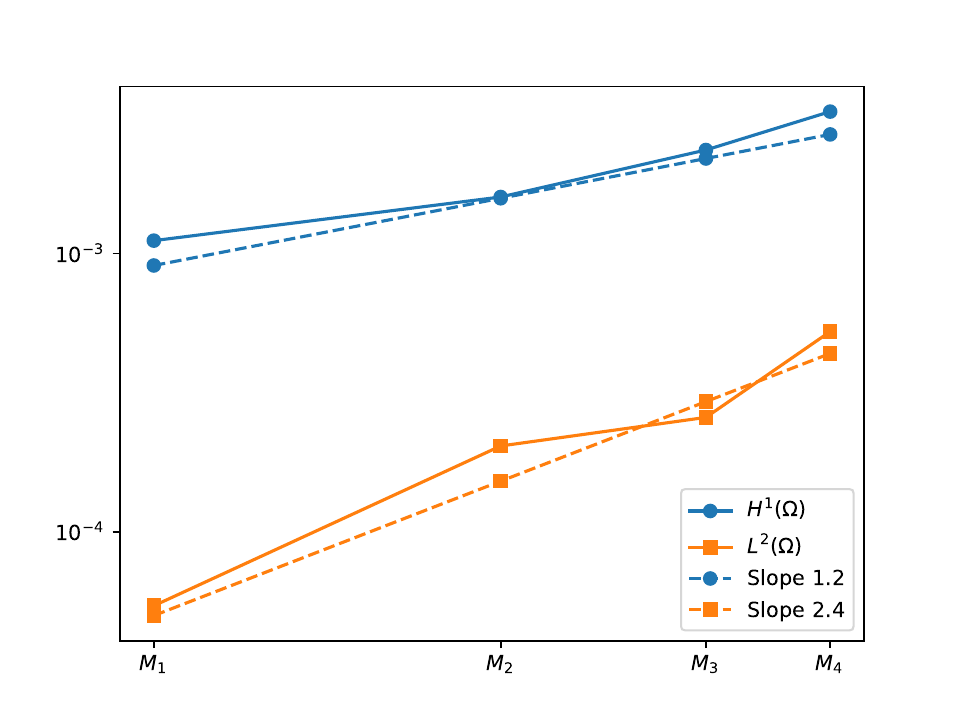}
	\hfill
	\includegraphics[width=0.45\textwidth]{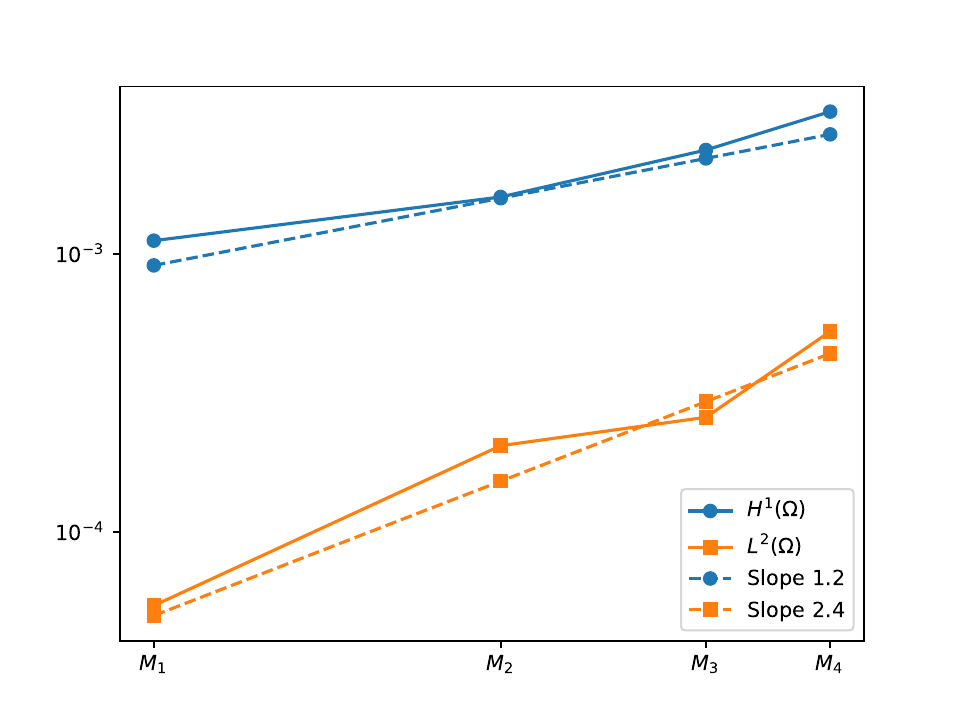}

	\caption{Evolution of the finite element error w.r.t. $h$, convex geometry. Left: $\eps = 10^{-2}$ and $M = 100$. Right: $\eps = 10^{-1}$ and $M = 1000$.}
	\label{fig_convergence_h_exp}
\end{figure}

We now come back to the convex geometry, and consider the function
\[u_s \colon (x, y) \mapsto e^x \sin(y).\]
We focus first on the finite element convergence rate, using the same meshes and procedure as in the previous section. We obtain \autoref{fig_convergence_h_exp}: the left figure corresponds to the choice $M = 100$ and $\varepsilon = 10^{-2}$, the right figure to the choice $M = 1000$ and $\eps = 10^{-1}$.

The results are very similar in both cases. Note that in both tests, the ratio $\frac{M}{\varepsilon}$ equals $10^5$, and \autoref{thm_errorestimatediscretization} shows that the estimates depend precisely on this ratio. Hence, the similarity in the results is expected; however, their near-perfect superposition is somewhat surprising.

\begin{figure}[hbt]
	\centering

	\includegraphics[width=0.6\textwidth]{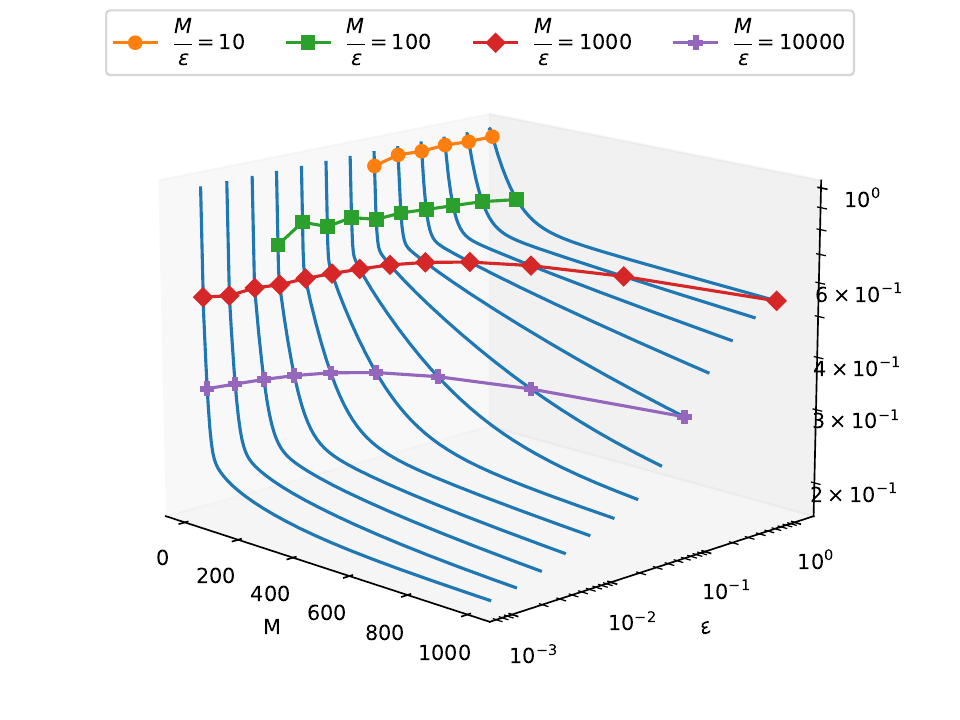}

	\caption{Illustration of the influence of the ratio $\frac{M}{\eps}$ over the error.}
	\label{fig_convergence_M}
\end{figure}

To investigate further, we now focus on the discrepency between the finite element approximation of the iterated quasi-reversibility method and the exact solution itself.
Indeed, our main estimate \autoref{thm_intro_second} highlights again the central role of the ratio $\frac{M}{\varepsilon}$.
\autoref{fig_convergence_M} shows the evolution of the error
\[\frac{\norm{u_{\eps,h}^M - u_s}[H^1(\Omega)] + \norm{\lambda_{\eps,h}^M - \lambda_s}[H^1_0(\Omega)]}{\norm{u_s}[H^1(\Omega)] + \norm{\lambda_s}[H^1_0(\Omega)]}\]
as a function of $M$, for different values of $\eps$, obtained on the mesh $M_1$. On top of these, we add some curves corresponding to fixed ratios $\frac{M}{\eps}$. This figure confirms the previous behavior: the error remains stable along these curves. In practice, there is a balance to find: to achieve a given precision, larger $\varepsilon$ requires larger $M$, whereas smaller $\varepsilon$ requires fewer iterations.


\subsection{Admissible strategies to ensure convergence in presence of noise}

While it is never easy to determine a "good" choice of parameters to ensure a satisfying reconstruction, 
the estimate
\[\norm{u_{\eps,h}^{M,\delta} - u_s}[H^1(\Omega)] \lesssim \sqrt{\frac{M}{\eps}} \delta + \left(\frac{M}{\eps}\right)^{3/2} h^\sigma + \norm{u_\eps^M - u_s}[H^1(\Omega)] + \norm{\lambda_\eps^M - \lambda_s}[H^1(\Omega)],\]
obtained in \autoref{thm_intro_second}, allows to define \emph{stategies of choice} to guarantee global convergence of the method.

First, as already mentioned in the introduction, we assign $\varepsilon$ the role of \emph{regularization parameter}. As such, we assume its value is fixed a priori, typically moderately small, and independent of all other parameters.
With that in mind, $M$ plays the role of the \emph{convergence parameter}, at least at the continuous level. This is evident in the last two terms of the estimate, that
go to zero as $M$ goes to infinity.
But we want the global error to go to zero, hence we need the first two terms in the error estimate to go to zero as well. In other words, working in the convex geometry for which we have $\sigma = 1$, we want $M$, $h$ and $\delta$ to satisfy
\[\sqrt{M} \delta \xrightarrow[]{} 0 \quad \text{and} \quad M^{3/2} h \xrightarrow[]{} 0.\]

In practical applications, the data is polluted by some noise, the amplitude of which is not controlled. In our setting, this means that the imposed parameter is $\delta$,
and we must fixed $M$ and $h$ as functions of $\delta$ to ensure the global convergence. 
In view of the estimate, a reasonable choice is to select $M$ and $h$ such that
 \[\sqrt{M} \delta \approx M^{3/2} h \xrightarrow{\delta \rightarrow 0} 0.\]
 Such a choice gives $\delta \approx M h$. As we want at the same time $M \xrightarrow{\delta \rightarrow 0} + \infty$ and $ \sqrt{M} \delta \xrightarrow{\delta \rightarrow 0} 0$, we deduce that we must have $\delta^3 \ll h \ll \delta$.

As such, in the following we set $h = c_h \delta^2$, and $M = \ceil*{c_M \delta^{-1}}$, where $c_h$ and $c_M$ are two constants chosen to obtain reasonable sizes of mesh and number of iterations. Note that these constants does not affect the asymptotics of the estimate. As in our context of numerical testing, it is easier to control $\delta$, as it can be fixed as we see fit, than to precisely control the mesh size $h$, in the numerics we proceed as follows: for one mesh $M_k$, we compute the exact value of $h$, deduce the corresponding amplitude of noise $\delta$  to apply to the data, and finally compute corresponding number of iterations $M$.

\begin{figure}[hbt]
	\centering

	\includegraphics[width=0.45\textwidth]{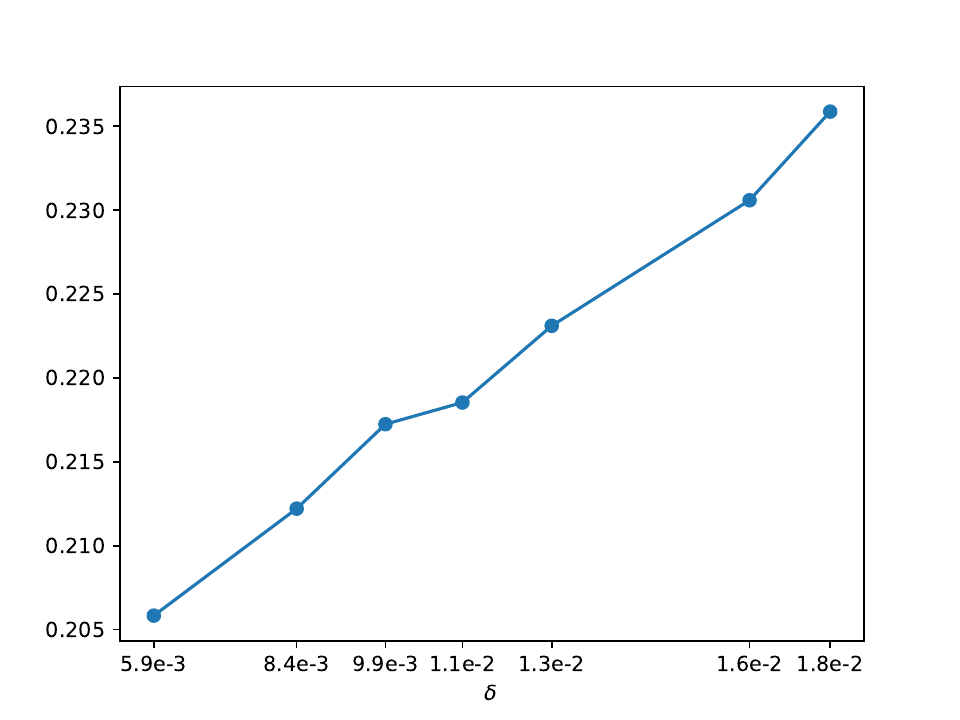}

	\caption{Illustration of the error decrease with respect to the amplitde of noise.}
	\label{fig_convergence_delta}
\end{figure}

For these simulations, we apply the following procedure to generate noise in the data. Once the exact solution $u_s$ is given, we compute the degrees of freedom of $f \coloneqq \Delta u_s$ and $g \coloneqq \restriction{u_s}{\omega}$. We then define two vectors of randomly generated numbers $f^\text{pert}$ and $g^\text{pert}$. This perturbation is normalized to become
\[\tilde{f}^\text{pert} \coloneqq \frac{f^\text{pert}}{\sqrt{\norm{f^\text{pert}}[L^2(\Omega)]^2 + \norm{g^\text{pert}}[L^2(\omega)]^2}}, \quad \tilde{g}^\text{pert} \coloneqq \frac{g^\text{pert}}{\sqrt{\norm{f^\text{pert}}[L^2(\Omega)]^2 + \norm{g^\text{pert}}[L^2(\omega)]^2}}.\]
Finally, for a given level of noise $\delta$, we work with the data
\[f^\delta \coloneqq f + \delta \tilde{f}^\text{pert}, \quad g^\delta \coloneqq g + \delta \tilde{g}^\text{pert},\]
and the rest of the numerical procedure remains the same as before. In \autoref{fig_convergence_delta}, we show the evolution of the relative error
\[\frac{\norm{u_{\eps,h}^{M,\delta} - u_s}[H^1(\Omega)]}{\norm{u_s}[H^1(\Omega)]}\]
with respect to $\delta$. As expected, the error decreases as $\delta$ goes to zero.

\section{Further comments} \label{section_further_comments}

\subsection{A feature of the proposed quasi-reversibility method}

Let \( u_s \in H^1(\Omega) \) denote the exact solution of the data assimilation problem that we aim to reconstruct. Let \( (f^\delta, \mathbf{F}^\delta, g^\delta) \in L^2(\Omega) \times L^2(\Omega)^d \times L^2(\omega) \) represent perturbed data in the sense of \eqref{eq_def_noise_Laplace}, and let \( \lambda^\delta \in H^1_0(\Omega) \) be defined by \eqref{eq_def_lambda}.

\begin{proposition}
Let $(u_\varepsilon^{M,\delta},\lambda_\varepsilon^{M,\delta})$ be the iterated quasi-reversibility sequence defined by \eqref{eq_defQRiterated}. For all $M\in \N_*$,
$$
\int_\omega u_\varepsilon^{M,\delta} dx = \int_\omega g^\delta dx.
$$
\end{proposition}

\begin{proof}
Choosing $(v,\mu) = (1,0)$ in \eqref{eq_defQRiterated} directly gives the result.
\end{proof}

At first glance, this result may appear anecdotal. However, it is important to recall that \( g^\delta \) represents the noisy data corresponding to \( u_s \) restricted to \( \omega \).
Since noise, if unbiased, acts as a high-frequency perturbation, it is reasonable to assume that the mean value of \( g^\delta \) over \( \omega \) remains close to the mean value of the exact solution over \( \omega \). Consequently, this property naturally extends to the approximated solution:
\[
\int_\omega u_\varepsilon^{M,\delta} \, dx = \int_\omega g^\delta \, dx \approx \int_\omega u_s \, dx.
\]
In other words, the mean value of the exact solution over \(\omega\) is a quantity that is reliably reconstructed by the method, regardless of the noise amplitude affecting the data.

\subsection{General second order elliptic equation}

This study primarily examines data assimilation for the Laplacian as a simplified case, though the methodology generalizes naturally to second-order elliptic operators satisfying the unique continuation property. As a typical example, consider a smooth, uniformly elliptic, symmetric matrix-valued function  $A : \overline{\Omega} \mapsto \mathcal{M}_d(\R)$  ($d=2$ or 3), \textit{i.e.}  $A \in C^\infty(\overline{\Omega};\mathcal{M}^d(\R))$, $A^T = A$ a.e. in $\Omega$, and there exists $c>0$ such that
$$
A \xi \cdot \xi \geq c \Vert \xi \Vert^2, \quad \forall \xi \in \R^d, \text{a.e. in } \Omega.
$$
Given $f^\delta\in L^2(\Omega)$ and $g^\delta \in L^2(\omega)$, we consider the following data assimilation problem for the operator $\nabla \cdot A \nabla$: find $u \in H^1(\Omega)$ such that
\begin{equation} \label{pbl_dataassimilation_divAgrad} \tag{$\mathcal{P}_{DA}^A$}
- \nabla\cdot A \nabla u  = f^\delta \text{ in } \Omega, \quad u = g^\delta \text{ in } \omega.
\end{equation}
Such problem shares the same properties as the data assimilation \eqref{eq_dataassimpb} \cite{lerousseaulebeaurobbiano22,vessella25}.

We define $\lambda^\delta$ as the unique element of $H^1_0(\Omega)$ satisfying
$$
\int_\Omega A \nabla \lambda^\delta \cdot \nabla \mu \, dx = \int_\Omega f^\delta \mu \, dx, \quad \forall \mu \in H^1_0(\Omega).
$$
Then, we define $a$,  continuous positive bilinear form
on $H^1(\Omega) \times H^1_0(\Omega)$, by
$$
a \big((u,\lambda),(v,\mu) \big) = \int_\Omega A \nabla \lambda\cdot \nabla v\, dx + \int_\omega u v \, dx +\int_\Omega A\nabla \lambda\cdot \nabla \mu dx - \int_\Omega A \nabla u \cdot \nabla \mu dx,
$$
$L$, continuous linear form on $H^1(\Omega) \times H^1_0(\Omega)$ , by
$$
L \big( (u,\lambda) \big) = \int_\Omega A \nabla \lambda^\delta \cdot \nabla v \, dx + \int_\omega g^\delta v \, dx,
$$
and set $b$ as in  \autoref{sec_contQRregres},
$$
b \big((u,\lambda), (v,\mu) \big) = \int_\Omega \nabla u \cdot \nabla v\, dx.
$$
In this setting, the abstract framework developed in  \autoref{sec_abstractQR} applies directly to \eqref{pbl_dataassimilation_divAgrad}.

The corresponding iterated quasi-reversibility sequence is initialized with $(u_\varepsilon^{\delta,0},\lambda_\varepsilon^{\delta,0}) = (0,0)$. For any iteration index
$M$ in $\N$,  the subsequent pair $(u_\varepsilon^{M+1,\delta},\lambda_\varepsilon^{M+1,\delta})$ is defined as the unique solution in $H^1(\Omega)\times H^1_0(\Omega)$ to the following coupled system:
\begin{equation*}
\begin{cases}
\displaystyle \int_\Omega A \nabla \lambda_\varepsilon^{M+1,\delta} \cdot \nabla v \, dx + \int_\omega u_\varepsilon^{M+1,\delta} v \, dx + \varepsilon \int_\Omega \nabla u_\varepsilon^{M+1,\delta} \cdot \nabla v \, dx
= \int_\Omega \nabla \lambda^\delta \cdot \nabla v \, dx + \int_\omega g^\delta v \, dx \\
\hfill \displaystyle + \varepsilon \int_\Omega  \nabla u_\varepsilon^{M,\delta} \cdot \nabla v \, dx, \quad \forall v \in H^1(\Omega),
\\
\\ \displaystyle \int_\Omega A \nabla \lambda_\varepsilon^{M+1,\delta} \cdot \nabla \mu \, dx - \int_\Omega A \nabla u_\varepsilon^{M+1,\delta} \cdot \nabla \mu \, dx = 0, \quad \forall \mu \in H^1_0(\Omega).
\end{cases}
\end{equation*}

It is straightforward to verify that all the estimates derived for the Laplacian operator remain valid in this more general setting. A crucial observation is that the elliptic regularity results for the Laplace operator presented in \autoref{appendix_miscres} continue to hold for the operator $\nabla \cdot A \nabla$, with the understanding that the constants $\sigma_D$ and $\sigma_N$ now depend not only on $\Omega$ but also on $A$. For further details on elliptic results for general elliptic operators in polygonal domains, we refer interested readers to Chapter~31.4 of~\cite{ernguermond2021} (see also~\cite{dauge88,grisvard85,grisvard92}).

As a result, there exists  $\sigma \in \left[ \frac{1}{2}, 1 \right)$ such that for every $\varepsilon_0 > 0$, there exists a constant $c > 0$ such that, for any $\delta >0$, for any $\varepsilon \in (0, \varepsilon_0]$ and any $M \in \mathbb{N}$,

\begin{multline*}
\Vert u_{\varepsilon,h}^{M,\delta} - u_s \Vert_{H^1(\Omega)} \leq c \bigg( \sqrt{\frac{M}{\varepsilon}} \delta
+ \left( \frac{M}{\varepsilon}\right)^{\frac{3}{2}} h^\sigma \sqrt{\Vert u_s \Vert^2_{H^1(\Omega)} + \Vert \nabla \cdot A \nabla u_s \Vert_{L^2(\Omega)}^2}
 \\ + \Vert u_{\varepsilon}^M - u_s \Vert_{H^1(\Omega)} + \Vert \lambda_\varepsilon^M - \lambda_s \Vert_{H^1(\Omega)}\bigg).
\end{multline*}

\subsection{Improved regularity}

To simplify the discussion, we assume here that $\partial \Omega$ is smooth.
As shown in \autoref{sec_contQRregres}, the $H^1 \times H^1$ quasi-reversibility formulation~\eqref{eq_itQRDA_core} used in this study inherently provides a regularity gain, which can be schematized as:
\begin{equation*}
    u_s \in H^1(\Omega, \Delta) \implies (u_\varepsilon^M, \lambda_\varepsilon^M) \in H^2(\Omega) \times H^2(\Omega).
\end{equation*}
However, this gain is inherently limited by the formulation itself: as noted in \autoref{rk_limitation_regularity}, even if the exact solution $u_s$ is smoother, the quasi-reversibility solutions do not exhibit higher regularity. This limitation arises from the terms integrated over $\omega$ in~\eqref{eq_itQRDA_core}.

A minor adjustment to the quasi-reversibility formulation resolves this issue. We begin by introducing a smooth cut-off function $\chi \in C^\infty_c(\Omega)$ such that $\chi = 0$ in $\Omega \setminus \omega$, and $\chi = 1$ in $\omega_0$,
where $\omega_0$ is a non-empty open set compactly supported in $\omega$.

Next, we define a modified continuous positive bilinear form on $H^1(\Omega) \times H^1_0(\Omega)$:
$$
    a_\chi\big((u,\lambda), (v,\mu)\big)
   = \int_\Omega \nabla \lambda \cdot \nabla v \, dx
    + \int_\Omega \chi^2 u v \, dx
    \quad + \int_\Omega \nabla \lambda \cdot \nabla \mu \, dx
    - \int_\Omega \nabla u \cdot \nabla \mu \, dx,
$$
along with a corresponding modified linear form:
\begin{equation*}
    L_\chi(v,\mu)
    = \int_\Omega \nabla \lambda^\delta \cdot \nabla v \, dx
    + \int_\Omega \chi^2 g^\delta \, v \, dx.
\end{equation*}
All results from \autoref{sec_contQRregres} remain valid when replacing $a$ with $a_\chi$ and $L$ with $L_\chi$, while retaining the same symmetric, positive bilinear form $b$.
The modified iterated quasi-reversibility sequence for unperturbed data is defined recursively as follows:
\begin{itemize}
    \item Initialize with $(u_\varepsilon^0, \lambda_\varepsilon^0) = (0, 0)$
    \item For each $M \in \mathbb{N}$, the pair $(u_\varepsilon^{M+1}, \lambda_\varepsilon^{M+1})$ is the unique solution in $H^1(\Omega) \times H^1_0(\Omega)$ satisfying:
\end{itemize}
\begin{equation*}
\begin{cases}
\displaystyle \int_\Omega \nabla \lambda_\varepsilon^{M+1} \cdot \nabla v \, dx + \int_\Omega \chi^2 u_\varepsilon^{M+1} v \, dx + \varepsilon \int_\Omega \nabla u_\varepsilon^{M+1} \cdot \nabla v \, dx
= \int_\Omega \nabla \lambda_s \cdot \nabla v \, dx + \int_\Omega \chi^2 u_s v \, dx \\
\hfill \displaystyle + \varepsilon \int_\Omega \nabla u_\varepsilon^M \cdot \nabla v \, dx, \quad \forall v \in H^1(\Omega),
\\
\\ \displaystyle \int_\Omega \nabla \lambda_\varepsilon^{M+1} \cdot \nabla \mu \, dx - \int_\Omega \nabla u_\varepsilon^{M+1} \cdot \nabla \mu \, dx = 0, \quad \forall \mu \in H^1_0(\Omega).
\end{cases}
\end{equation*}
Applying the same reasoning as in \autoref{rk_limitation_regularity}, we find that $\lambda_\varepsilon^1 = \lambda_\varepsilon$ satisfies:
\begin{equation*}
\Delta \lambda_\varepsilon = \frac{1}{1+\varepsilon} \Delta \lambda_s + \frac{1}{1+\varepsilon} \chi^2 (u_\varepsilon - u_s) \quad \text{in } \Omega.
\end{equation*}
Therefore, if $u_s$ belongs to $H^3(\Omega)$, so does $\lambda_s$, which gives now $\Delta \lambda_\varepsilon \in H^1(\Omega)$. As $\lambda_\varepsilon$ is in $H^1_0(\Omega)$, we deduce that $\lambda_\varepsilon$ belongs to $H^3(\Omega)$: we have gain one order of regularity  in comparaison to the standard formulation. The same gain propagate to the full sequence $(u_\varepsilon^M,\lambda_\varepsilon^M)$, which is now as smooth as the exact solution of the problem.

\medskip

In polyhedral domains, regularity is inherently limited near reentrant corners. However, with the modified quasi-reversibility method, the elliptic regularity property holds away from the boundary: $u_s \in H^3(\Omega) \Rightarrow u_\varepsilon^{M} \in H^3_{loc}(\Omega)$.

\medskip

The main divergence in the study of this new formulation emerges at the computational level. Specifically, the method now requires numerical evaluation of integrals of the form
\[
\int_\Omega \chi^2 u_h v_h \, dx,
\]
where \( u_h \) and \( v_h \) are basis functions in the finite element space \( V_h \). The study's underlying assumption, that \emph{all appearing integrals can be computed with exact numerical precision}, is no longer realistic in practice. Consequently, the errors arising from approximate quadrature must now be explicitly addressed. A rigorous analysis of this quasi-reversibility variant, including these computational considerations, will therefore be pursued in subsequent work.

\subsection{Equivalence between \texorpdfstring{\eqref{eq_defQRiterated}}{(\ref{eq_defQRiterated})} and a symmetric variational problem}

Interestingly, one can recast the proposed non-symmetrical iterated quasi-reversibility method studied in the present article into a symmetrical formulation. To do so, we first define the linear continuous operator $B : u \in H^1(\Omega) \mapsto Bu \in H^1_0(\Omega)$, where $Bu$ is the unique element of $H^1_0(\Omega)$
such that
$$
\int_\Omega \nabla Bu \cdot \nabla \mu dx = \int_\Omega \nabla u \cdot \nabla \mu \,dx, \quad \forall \mu \in H^1_0(\Omega).
$$
In other terms, $Bu$ is the projection of $u$ onto $H^1_0(\Omega)$.

We now define an alternative quasi-reversibility sequence: $w_\varepsilon^{0,\delta} = 0$ and for any $M$ in $\N$, $w_\varepsilon^{M+1,\delta}$
is the unique element of $H^1(\Omega)$ such that
\begin{multline} \label{eq_def_alternative_symQR}
\int_\Omega \nabla B w_\varepsilon^{M+1,\delta}\cdot \nabla B v \, dx + \int_\omega w_\varepsilon^{M+1,\delta} v dx
+ \varepsilon \int_\Omega \nabla w_{\varepsilon}^{M+1,\delta} \cdot \nabla v dx = \int_\Omega \nabla \lambda^\delta \cdot\nabla v dx  \\
+ \int_\omega g^\delta v dx + \varepsilon \int_\Omega \nabla w_\varepsilon^{M,\delta} \cdot \nabla v \, dx, \quad \forall v \in H^1(\Omega).
\end{multline}
As before, $\lambda^\delta$ is defined by \eqref{eq_def_lambda}, and $g^\delta \in L^2(\omega)$ is given in \eqref{eq_dataassimpb}.
It is not complicated to obtain that the sequence $(w_\varepsilon^{M,\delta})_{M\in \N}$ is well-defined. Note that the bilinear form appearing on the left-hand  side of Equation  \eqref{eq_def_alternative_symQR} is symmetric, yet it is equivalent to the method introduce in the present paper:

\begin{proposition}
For all $\varepsilon>0$ and all $M\in \N$, we have $(u_\varepsilon^{M,\delta},\lambda_\varepsilon^{M,\delta}) = (w_\varepsilon^{M,\delta},B w_\varepsilon^{M,\delta})$.
\end{proposition}

\begin{proof}
The second equation in~\eqref{eq_defQRiterated} gives
$
\lambda_\varepsilon^{M,\delta} = B u_\varepsilon^{M,\delta}
$
for any integer $M$.
Consequently, to demonstrate equivalence between the two formulations, it suffices to show that
$
u_\varepsilon^{M,\delta} = w_\varepsilon^{M,\delta}
$
for all iteration indices $M$.

Assume this equality holds for some fixed integer $M$. For any test function $v \in H^1(\Omega)$, we have:
$$
\int_\Omega \nabla \lambda_\varepsilon^{M+1,\delta} \cdot \nabla v \, dx
= \int_\Omega \nabla \lambda_\varepsilon^{M+1,\delta} \cdot \nabla B v \, dx
= \int_\Omega \nabla B u_\varepsilon^{M+1,\delta} \cdot \nabla B v \, dx.
$$
Substituting this identity into the first equation of~\eqref{eq_defQRiterated} and invoking our inductive hypothesis that
$
u_\varepsilon^{M,\delta} = w_\varepsilon^{M,\delta},
$
we deduce that $u_\varepsilon^{M+1,\delta}$ satisfies~\eqref{eq_def_alternative_symQR}. By the uniqueness of solutions to this problem, we conclude that $
u_\varepsilon^{M+1,\delta} = w_\varepsilon^{M+1,\delta}.
$

Since the base case
$u_\varepsilon^{0,\delta} = 0 = w_\varepsilon^{0,\delta}
$ holds by definition, the result follows by induction.
\end{proof}

At first glance, the second version of the iterated quasi-reversibility method \eqref{eq_def_alternative_symQR} appears simpler than the original formulation \eqref{eq_defQRiterated}. It reduces the problem to a single unknown and exhibits symmetry, properties that, notably, allow the direct application of the theoretical framework developed in~\cite{darde2016}.

While this partially holds at the continuous level, the explicit presence of the operator \( B \) on the left-hand side of \eqref{eq_def_alternative_symQR} complicates its discretization. The only practical approach to discretizing \eqref{eq_def_alternative_symQR} with finite elements is to revert to the original form \eqref{eq_defQRiterated}. For this reason, it is most logical to focus the analysis on \eqref{eq_defQRiterated} from the outset.


\appendix

\section{Miscellaneous results} \label{appendix_miscres}

\subsection*{Elliptic regularity results}

We present the elliptic regularity results used in our study.
Let $\Theta$ be a bounded open set of $\R^d$, with a Lipschitz boundary.
Let $f$ belong to $L^2(\Theta)$.
We consider the following two problems: find $u_D \in H^1_0(\Theta)$ such that
for all $v \in H^1_0(\Theta)$, we have
\begin{equation} \label{eq_Prob_DiriHom} \tag{$\mathcal{P}_{diri}$}
\int_\Theta \nabla u_D \cdot \nabla v dx = \int_\Theta f \, vdx.
\end{equation}
and find $u_N \in H^1(\Theta)$ such that for
all $v \in H^1(\Theta)$,
\begin{equation} \label{eq_Prob_NeumHom} \tag{$\mathcal{P}_{neum}$}
\int_\Theta \left(\nabla u_N \cdot \nabla v + u_N\, v\right) dx  = \int_\Theta f v dx.
\end{equation}
Problem \eqref{eq_Prob_DiriHom} corresponds to the weak  formulation of the second-order problem $-\Delta u = f$ in $\Theta$, with homogeneous Dirichlet condition on $\partial \Theta$,
whereas problem \eqref{eq_Prob_NeumHom} corresponds to the weak formulation of the second-order problem $-\Delta u + u = f $ posed in $\Theta$, with homogeneous Neumann boundary conditions on $\partial \Theta$.

Below, we give the main regularity results we use in the study. They can all be found in \cite[Section 1.5.2]{ciarletluneville23}

\begin{theorem} \label{thm_regularityresults_Diri}
	Problem \eqref{eq_Prob_DiriHom} admits a unique solution $u_D \in H^1_0(\Theta)$. Furthermore,
	\begin{itemize}
		\item suppose $\partial \Theta$ is of class $C^{1,1}$, or $d= 2$ and $\Theta$ is a convex polygon, or $d=3$ and $\Theta$ is a convex polyhedron. Then $u_D$ belongs to $H^2(\Theta)$, and there exists $c$, depending only on $\Theta$, such that
		$$
		\Vert u_D \Vert_{H^2(\Theta)} \leq c \Vert f \Vert_{L^2(\Theta)}.
		$$
		\item suppose $d=2$ and $\Theta$ is a polygon, or $d=3$ and $\Theta$ is a polyhedron. Then there exists $\sigma_D \in \left]1/2,1 \right]$ and $c$, depending only on $\Theta$, such that $u_D$ belongs to $H^{1+\sigma_D}(\Theta)$ and
		$$
		\Vert u_D \Vert_{H^{1+\sigma_D}(\Theta)} \leq c \Vert f \Vert_{L^2(\Theta)}.
		$$
	\end{itemize}
\end{theorem}

\begin{proof}
	The existence and uniqueness is a standard application of Lax-Milgramm theorem and Poincar\'e's inequality. The regularity results for domains with $C^{1,1}$ boundary is proven in \cite[Theorem 2.2.2.3 and Corollary 2.2.2.4]{grisvard85}. The case of a two-dimensional polygon (convex or not) is covered in \cite[Theorem 2.2.3 and Theorem 2.4.3]{grisvard92}.
	The result for three-dimensional polyhedron is given in \cite[Corollaries 2.6.7 and 2.6.8]{grisvard92}.
\end{proof}

\begin{theorem} \label{thm_regularityresults_Neum}
	Problem \eqref{eq_Prob_NeumHom} admits a unique solution $u_N \in H^1(\Theta)$. Furthermore,
	\begin{itemize}
		\item suppose $\partial \Theta$ is of class $C^{1,1}$, or $d= 2$ and $\Theta$ is a convex polygon, or $d=3$ and $\Theta$ is a convex polyhedron. Then $u_N$ belongs to $H^2(\Theta)$, and there exists $c$, depending only on $\Theta$, such that
		$$
		\Vert u_N \Vert_{H^2(\Theta)} \leq c \Vert f \Vert_{L^2(\Theta)}.
		$$
		\item suppose $d=2$ and $\Theta$ is a polygon, or $d=3$ and $\Theta$ is a polyhedron. Then there exists $\sigma_N \in \left]1/2,1 \right]$ and $c$, depending only on $\Theta$, such that $u_N$ belongs to $H^{1+\sigma_N}(\Theta)$ and
		$$
		\Vert u_N \Vert_{H^{1+\sigma_N
		}(\Theta)} \leq c \Vert f \Vert_{L^2(\Theta)}.
		$$
	\end{itemize}
\end{theorem}

\begin{proof}
	The existence and uniqueness is a standard application of Lax-Milgramm theorem. The regularity results for domains with $C^{1,1}$ boundary is proven in \cite[Theorem 2.2.2.5 and Corollary 2.2.2.6]{grisvard85}. The case of a two-dimensional polygon (convex or not) is an adaptation of \cite[Theorem 2.2.3 and Theorem 2.4.3]{grisvard92} to the case of homogeneous Neumann condition. The result for three-dimensional polyhedron is given in \cite[Corollary 2.6.7 and Remark 2.6.9]{grisvard92}.
\end{proof}

\begin{remark}
In dimension 2, the parameters $\sigma_D$ and $\sigma_N$, which characterize the regularity of the solutions to the Dirichlet and Neumann problems, respectively, are actually equal. This equality does not generally hold in dimension 3; see \cite{ciarletluneville23,grisvard92} for further details.
\end{remark}

\subsection*{Conditional stability estimates for the data assimilation problem}

In the study, the $L^2$ convergence estimates for the iterated quasi-reversibilty method are obtained via the following stability estimate from \cite{alessandrinirondirossetvessella2009}.

\begin{theorem}[{\cite[Theorem 5.3]{alessandrinirondirossetvessella2009}}] \label{thm_GlobalStability_Alessandrini}
	Let $\Omega$ be a bounded open set of $\R^d$ with Lipschitz boundary, and $\omega$ be a non-empty open set compactly contained in $\Omega$. There exists $c>0$ and $\mu \in (0,1)$ such that, for all $u \in H^1(\Omega)$ satisfying $\Delta u = f \in L^2(\Omega)$ in $\Omega$, $u_{|\omega} = g \in L^2(\omega)$ with the following a priori bounds
	$$
	\Vert f \Vert_{L^2(\Omega)} \leq \eta_f, \quad \Vert g \Vert_{L^2(\omega)} \leq \eta_g, \quad \Vert u \Vert_{H^1(\Omega)} \leq \mathcal K,
	$$
	where $\mathcal K$, $\eta_g$ and $\eta_f$ satisfy
	$$
	\frac{\eta_g + \eta_f}{\mathcal K + \eta_f} < 1,
	$$
	one has
	$$
	\Vert u \Vert_{L^2(\Omega)} \leq c \left(\mathcal K + \eta_f \right) \ln \left(\frac{\mathcal K + \eta_f}{\eta_g + \eta_f}\right)^{-\mu}.
	$$
\end{theorem}

\subsection*{Global approximation properties of $P_1$-Lagrange finite elements}

As introduced page \pageref{subsection_intro_discretization}, $\mathcal{T}_h$ is an admissible, shape regular partition of $\Omega$ into triangular elements in dimension 2, and tetrahedral elements in dimension 3. We recall that $V_h$ refers to the standard conforming $P_1$-Lagrange finite element space, and $V_{h,0} = V_h \cap H^1_0(\Omega)$.

Our convergence analysis relies on the following well-known finite element approximation result.

\begin{theorem} \label{thm_globalapprox_FE}
Let $\eta \in [1,2]$. There exists $c>0$ such that for all $w \in H^\eta(\Omega)$ and $v \in H^\eta \cap H^1_0(\Omega)$,
$$
\inf_{w_h\in V_h} \Vert w - w_h \Vert_{L^2(\Omega)} + h \inf_{w_h \in V_h} \Vert w - w_h \Vert_{H^1(\Omega)} \leq c h^{\eta} \Vert w \Vert_{H^\eta(\Omega)},
$$
$$
\inf_{v_h\in V_{h,0}} \Vert v - v_h \Vert_{L^2(\Omega)} + h \inf_{v_h \in V_{h,0}} \Vert v - v_h \Vert_{H^1(\Omega)} \leq c h^{\eta} \Vert v \Vert_{H^\eta(\Omega)}.
$$
\end{theorem}

The $H^s(\Omega)$ case is proven in {\cite[Theorem 9.10]{steinbach08}}, and the $H^\eta\cap H^1_0(\Omega)$ case is discussed in \cite[Section 6.3.1]{grisvard92}.
We also refer to \cite{ciarletluneville23} and \cite{grisvard85} on related topics.

\subsection*{Technical lemmas}

\begin{lemma} \label{lemma_estimate_series}
	Let $(a_n)_{n\in \N_*}$ be a non-increasing summable sequence of positive real numbers. Set $\displaystyle c  = \sum_{n\in \N_*} a_n$.
	Then, for all $n \in \N_*$, we have $ n\, a_n \leq c$.
\end{lemma}

\begin{proof}
	We simply note that, because $0\leq a_n \leq a_k$ for all $n\in \N$ and all $k \in \left\lbrace1,\ldots, n\right\rbrace$,
	$$
	n \, a_n = \sum_{k=1}^n a_n \leq \sum_{k=1}^n a_k \leq \sum_{k\in \N_*} a_k = c.
	$$
\end{proof}

\begin{lemma} \label{lemma_wm_explicit} Let $\mathcal B$ be a Banach space.
	Let $q\neq 1$, $\kappa \in \mathcal B$ and $\eta \in \mathcal B^\N$. Let $w \in \mathcal B^\N$ be defined recursively by, for all integer $m$,  $w_{m+1} = q w_m + \kappa + \eta_m$ . Then
	$$
	w_m = q^m \,w_0 + \frac{1-q^m}{1-q} \kappa + \sum_{k=0}^{m-1} q^k \eta_{m-1-k}.
	$$
\end{lemma}

\begin{proof} The formula holds for $q= 0$. Let $q \notin \left\lbrace 0,1\right\rbrace$, and set
	 $a_m = w_m q^{-m}$. Then $a_0 = w_0$, and, for all integer $m$, $a_{m+1} = a_m + q^{-(m+1)} \kappa + q^{-(m+1)} \eta_m$. As a consequence, for all $m\in \N$,
	$$
	a_m = a_0 \sum_{k=0}^{m-1} (a_{k+1} - a_k) = a_0 + \kappa \sum_{k=0}^{m-1} q^{-(k+1)} + \sum_{k=0}^{m-1} \eta_k q^{-(k+1)}
	= a_0 + q^{-m} \frac{1-q^m}{1-q} \kappa + q^{-1} \sum_{k=0}^{m-1} \frac{\eta_k}{q^k}.
	$$
	The result follows.
\end{proof}

\begin{lemma} \label{lemma_bound_sum}
There exists $f \in C^0([0,1); \R_+)$ such that for all $m\in \N_*$, for all $q \in [0,1)$,
$$
\sum_{p=0}^{m-1} \frac{q^p}{\sqrt{m-p}} \leq \frac{f(q)}{(1-q)\sqrt m }.
$$
\end{lemma}

\begin{proof} Let  $\eta$ be in $(0,1)$. We have
$$
(1-q) \sqrt m \sum_{p=0}^{\lfloor \eta m \rfloor}  \frac{q^p}{\sqrt{m-p}}
 \leq \frac{1}{\sqrt{1-\eta}} (1-q) \sum_{p=0}^{\lfloor \eta m \rfloor} q^p
 \leq \frac{1}{\sqrt{1-\eta}},
$$
and, setting $g : x >0 \mapsto  x q^{2\eta x}$
$$
(1-q) \sqrt m \sum_{p=\lfloor \eta m \rfloor+1}^{m-1} \frac{q^p}{\sqrt{m-p}}
=\sqrt m\, q^{\lfloor \eta m \rfloor+1}\,
(1-q) \sum_{p=0}^{m-\lfloor \eta m \rfloor} q^p
 \leq \sqrt m q^{\eta m} \leq \sqrt{ \sup_{x>0} g(x)}.
$$
A direct calculation gives $\displaystyle \sup_{x>0} g(x) = (2e\, \eta\, \vert \ln(q) \vert)^{-1}$. Consequently, for all $\eta \in (0,1)$,
$$
(1-q) \sqrt m \sum_{p=0}^{m-1} \frac{q^p}{\sqrt{m-p}} \leq \frac{1}{\sqrt{1-\eta}} + \frac{1}{\sqrt{2 e \eta \vert \ln (q) \vert}}.
$$
The right-hand side is minimal for $\eta = \sqrt[3]{1+ 2 e \vert \ln (q)\vert}$, which gives
$$
(1-q) \sqrt m \sum_{p=0}^{m-1} \frac{q^p}{\sqrt{m-p}} \leq \frac{\left( 1+ \sqrt[3]{2 e \vert \ln(q) \vert}\right)^{\frac{3}{2}}}{\sqrt{2 e \vert \ln q\vert}}.
$$
We can therefore set
$$
f : q \in [0,1) \mapsto 1 \text{ if } q=0,\ \frac{\left( 1+ \sqrt[3]{2 e \vert \ln(q) \vert}\right)^{\frac{3}{2}}}{\sqrt{2 e \vert \ln q\vert}} \text{ otherwise},
$$
which ends the proof.

\end{proof}

\end{document}